\documentclass[11pt,a4paper]{article}

\usepackage[left=2cm,right=2cm,top=2cm,bottom=2cm]{geometry}
\usepackage{graphicx}
\usepackage{multirow}
\usepackage{amsmath,amssymb,amsfonts}
\usepackage{amsthm}
\usepackage{mathrsfs}
\usepackage[title]{appendix}
\usepackage{xcolor}
\usepackage{textcomp}
\usepackage{manyfoot}
\usepackage{booktabs}
\usepackage{algorithm}
\usepackage{algorithmicx}
\usepackage{algpseudocode}
\usepackage{listings}
\ifdefined\debug
  \usepackage[inline]{showlabels}
  \showlabels{ref}
  \showlabels{cite}
  \usepackage{fancyhdr}
\fi
\usepackage{bm}
\usepackage{hyperref}

\theoremstyle{plain}
\newtheorem{theorem}{Theorem}[section]
\newtheorem{proposition}[theorem]{Proposition}%
\newtheorem{lemma}[theorem]{Lemma}

\theoremstyle{definition}
\newtheorem{definition}[theorem]{Definition}
\newtheorem{remark}[theorem]{Remark}

\newtheorem{question}[theorem]{Question}
\newtheorem*{notation}{Notation}
\theoremstyle{plain}
\newtheorem{maintheorem}{Theorem}
\renewcommand{\themaintheorem}{\Alph{maintheorem}}
\newcommand{\wiptext}[1]{\textcolor{blue}{#1}}
\ifdefined\debug
  \newcommand{\debugtext}[1]{\textcolor{red}{#1}}
  \renewcommand{\headrulewidth}{0pt}
\else
  \newcommand{\debugtext}[1]{\ignorespaces}
\fi

\counterwithin{figure}{section}
\counterwithin{table}{section}

\numberwithin{equation}{section}

\newcommand{\isMaster}{}

\begin{document}

\title{
  \textbf{
    Heterodimensional cycles derived from homoclinic tangencies
    via Hopf bifurcations
  }
}

\author{Shuntaro Tomizawa
  \footnote{tomizawa-s@g.ecc.u-tokyo.ac.jp,
    Graduate School of Mathematical Sciences,
    The University of Tokyo, 
    3-8-1 Komaba, Meguro, Tokyo, 153-8914, Japan
}}

\maketitle
\ifdefined\debug
  \thispagestyle{fancy}
\fi

\noindent \textbf{Abstract.}
We analyze three-dimensional $C^r$ diffeomorphisms ($r \geq 5$) exhibiting a quadratic
focus–saddle homoclinic tangency whose multipliers satisfy
$|\lambda \gamma| = 1$.  
For a proper unfolding family with three-parameters
that split the tangency,
vary the argument of the stable multipliers,
and control the modulus $|\lambda \gamma|$,
we show that a Hopf bifurcation occurs on this curve
and that a homoclinic point to the bifurcating periodic orbit is present.  
As a consequence,
the original map $f$ can be $C^r$-approximated by a diffeomorphism
exhibiting a coindex one heterodimensional cycle in the saddle case.

\medskip


\medskip


\tableofcontents

\ifdefined\isMaster
\else
  \documentclass[11pt,a4paper]{article}

\usepackage[left=2cm,right=2cm,top=2cm,bottom=2cm]{geometry}
\usepackage{graphicx}
\usepackage{multirow}
\usepackage{amsmath,amssymb,amsfonts}
\usepackage{amsthm}
\usepackage{mathrsfs}
\usepackage[title]{appendix}
\usepackage{xcolor}
\usepackage{textcomp}
\usepackage{manyfoot}
\usepackage{booktabs}
\usepackage{algorithm}
\usepackage{algorithmicx}
\usepackage{algpseudocode}
\usepackage{listings}
\ifdefined\debug
  \usepackage[inline]{showlabels}
  \showlabels{ref}
  \showlabels{cite}
  \usepackage{fancyhdr}
\fi
\usepackage{bm}
\usepackage{hyperref}

\theoremstyle{plain}
\newtheorem{theorem}{Theorem}[section]
\newtheorem{proposition}[theorem]{Proposition}%
\newtheorem{lemma}[theorem]{Lemma}
\newtheorem{corollary}[theorem]{Corollary}
\newtheorem{claim}[theorem]{Claim}
\theoremstyle{definition}
\newtheorem{definition}[theorem]{Definition}
\newtheorem{remark}[theorem]{Remark}
\newtheorem{example}[theorem]{Example}
\newtheorem{question}[theorem]{Question}
\newtheorem*{notation}{Notation}
\theoremstyle{plain}
\newtheorem{maintheorem}{Theorem}
\renewcommand{\themaintheorem}{\Alph{maintheorem}}
\newcommand{\wiptext}[1]{\textcolor{blue}{#1}}
\ifdefined\debug
  \newcommand{\debugtext}[1]{\textcolor{red}{#1}}
  \pagestyle{fancy} 
  \fancyhead[L]{}
  \fancyhead[C]{}
  \fancyhead[R]{\textcolor{red}{[DEBUG] Debug mode is activated.}}
  \renewcommand{\headrulewidth}{0pt}
\else
  \newcommand{\debugtext}[1]{\ignorespaces}
\fi

\counterwithin{figure}{section}
\counterwithin{table}{section}

\raggedbottom

\numberwithin{equation}{section}

  \begin{document}
\fi

\section{Introduction}
In smooth dynamical systems,
complicated behavior often appears when the system is not uniformly hyperbolic.
Two important phenomena that cause such behavior are \textit{homoclinic tangencies}
and \textit{heterodimensional cycles}.

A homoclinic tangency means that the stable and unstable manifolds of
a hyperbolic periodic point intersect in a non-transversal way.
This kind of intersection can produce complicated dynamics,
such as infinitely many sinks or sources,
or strange attractors.
The phenomenon of homoclinic tangency was first observed in \cite{N1970}.
Later studies revealed its deep connection with the so-called \textit{Newhouse domain},
where persistent homoclinic tangencies and infinitely many sinks can coexist;
for instance,
the studies of the domain are \cite{N1979,GTS1993,R1995,GTS1997,GST2008,BD2012,L2024}.
Furthermore, homoclinic tangencies have been studied in connection with
non-hyperbolic properties,
including the occurrence of zero Lyapunov exponents \cite{DG2009},
the divergence of Birkhoff averages \cite{KS2017,B2022,BB2023,KNS2023},
the emergence of infinitely many sinks \cite{B2016,N1974,PV1994,GST2008},
and the complexity of bifurcation structures \cite{T2010,T2015}.

A heterodimensional cycle is a situation where two hyperbolic periodic points
have different unstable indices
(that is, different dimensions of their unstable manifolds),
and their invariant manifolds intersect in both directions.
Such a cycle were discovered in \cite{AS1970,S1972}.
Later, Bonatti and Díaz identified regions,
now called the \textit{Bonatti–Díaz domains},
where such cycles occur robustly \cite{BD1996}.
Subsequent studies have explored the dynamical complexity
within these domains \cite{BD2008,BDK2012,LT2024,L2024_2},
as well as other forms of rich behavior arising from heterodimensional structures
\cite{AST2017,AST2021,GIKN2005}.

In recent studies,  
researchers have found a strong connection between a homoclinic tangency  
and a heterodimensional cycle,  
Many studies have investigated this relationship in depth
\cite{DU1994,DG2009,L2017,LT2017,BDP2022,LLST2022}.
Understanding this connection is important for studying non-hyperbolic dynamics.

Another important bifurcation related to non-hyperbolic dynamics is the Hopf bifurcation.
This bifurcation occurs when a fixed point of a nonlinear system loses its stability,
and a limit cycle appears or disappears.
In continuous-time systems of dimension two or higher,
the Hopf bifurcation plays a key role in the emergence of oscillatory behavior,
such as nonlinear or self-excited vibrations \cite{H1943}.
Neimark and Sacker extended the Hopf bifurcation to discrete-time systems.
The bifurcation is now known as the Neimark–Sacker bifurcation \cite{N1959,S2009}.
This discrete analogue also creates invariant closed curves from fixed points,
and is fundamental in the study of bifurcations in maps.
It has been observed that  
homoclinic tangencies and heterodimensional cycles
can occur near Neimark–Sacker bifurcations,  
especially when the system exhibits a Hopf-homoclinic cycle \cite{M1999,T2019}.  

\subsection{Previous work and our approach}
\label{s11-VS}
The prior work related to our research is the study of the relationship  
between homoclinic tangencies and heterodimensional cycles \cite{LLST2022}.  
In this work,
they study homoclinic tangencies in a manifold $M_\mathrm{ph}$
with $\dim M_\mathrm{ph} \geq 3$.
Let $\Gamma$ be an orbit of a homoclinic tangency to a hyperbolic periodic point
$O^*$ of a $C^r$, $r \in \mathbb{Z}_{> 0} \cup \{\infty, \omega\}$,
diffeomorphism $f$.
Here, we write $\mathbb{Z}_{> 0} := \{1, 2, \cdots\}$.
We denote the multipliers of $O^*$, which are the eigenvalues of
$D(f^{\mathrm{per}(O^*)})_{O^*}$, by
\begin{align*}
  \lambda^*_{d^s}, \,
  \lambda^*_{d^s - 1}, \,
  \cdots, \,
  \lambda^*_1, \,
  \gamma^*_1, \,
  \gamma^*_2, \,
  \cdots, \,
  \gamma^*_{d^u}
\end{align*}
with
\begin{align*}
  |\lambda^*_{d^s}| \leq |\lambda^*_{d^s - 1}| \leq \cdots \leq |\lambda^*_1|
  < 1 <
  |\gamma^*_1| \leq |\gamma^*_2| \leq \cdots \leq |\gamma^*_{d^u}|,
\end{align*}
where $d^s$ and $d^u$ indicate the stable and unstable index of $O^*$, respectively
and $\mathrm{per}(O^*)$ is the period of $O^*$.
The \textit{center-stable and center-unstable multipliers} of $O^*$ are the ones
closest to the unit circle, with the former just inside and the latter just outside.
By an arbitrarily $C^r$ small perturbation, we may assume that the central
multipliers are just $\lambda^*_1$ and $\gamma^*_1$ and their complex conjugates, if any.
Such a generic orbit of a homoclinic tangency has several classes:
\begin{itemize}
  \item \textbf{Saddle~(1, 1)}: $\lambda^*_1, \gamma^*_1 \in \mathbb{R}$.
  \item \textbf{Saddle-Focus~(1, 2)}: $\lambda^*_1 \in \mathbb{R}$,
    and $\gamma^*_1 = \overline{\gamma^*_2} = \gamma^* \mathrm{e}^{i \omega^*}$
    for some $\gamma^*$ with $|\gamma^*| > 1$ and $\omega^* \in (0, \pi)$.
  \item \textbf{Focus-Saddle~(2, 1)}:
    $\lambda^*_1 = \overline{\lambda^*_2} = \lambda^* \mathrm{e}^{i \omega^*}$
    for some $\lambda^* \in (0, 1)$ and $\omega^* \in (0, \pi)$,
    and $\gamma^*_1 \in \mathbb{R}$.
  \item \textbf{Bi-Focus~(2, 2)}:
    $\lambda^*_1 = \overline{\lambda^*_2} = \lambda^* \mathrm{e}^{i \omega^*_1}$
    for some $\lambda^* \in (0, 1)$ and $\omega^*_1 \in (0, \pi)$,
    and $\gamma^*_1 = \overline{\gamma^*_2} = \gamma^* \mathrm{e}^{i \omega^*_2}$
    for some $\gamma^*$ with $|\gamma^*| > 1$ and $\omega^*_2 \in (0, \pi)$,
\end{itemize}
where $\mathrm{e}$ is the base of the natural logarithm.
The above terminologies are based on papers \cite{GST2008,LLST2022}.

\begin{figure}[h]
  \centering
  \includegraphics{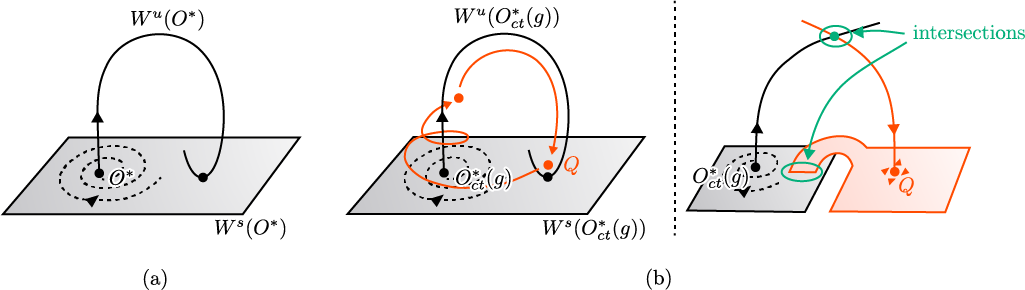}
  \caption{
    (a) The phase portrait of $f$ in the Focus-Saddle~(2, 1) class when
    $\dim M_\mathrm{ph} = 3$. (b) The phase portrait of $g$ having a heterodimensional
    cycle involving $O^*(g)$ and $Q$.
    The new hyperbolic periodic point $Q$ arises near the orbit of
    the homoclinic tangency.
    The right picture indicates the cycle in a topological view.}
    \label{fig-phase_LLST}
  \end{figure}

Depending on the product $|\lambda^* \gamma^*|$, we can generally consider cases shown
in Table~\ref{tab-cases_of_HT}.
\begin{table}[h]
  \centering
  \caption{Generic cases of homoclinic tangencies}\label{tab-cases_of_HT}
  \begin{tabular}{l|lc}
    \midrule
    \textbf{Name} & \textbf{Class} & \textbf{The product $|\lambda^* \gamma^*|$} \\
    \midrule
    Case~(1, 1)-Sm & \multirow{2}{*}{Saddle~(1, 1)} & $< 1$ \\
    Case~(1, 1)-Lg & & $> 1$ \\
    \midrule
    Case~(1, 2)-Sm & \multirow{2}{*}{Saddle-Focus~(1, 2)} & $< 1$ \\
    Case~(1, 2)-Lg & & $> 1$ \\
    \midrule
    Case~(2, 1)-Sm & \multirow{2}{*}{Focus-Saddle~(2, 1)} & $< 1$ \\
    Case~(2, 1)-Lg & & $> 1$ \\
    \midrule
    Case~(2, 2)-Sm & \multirow{2}{*}{Bi-Focus~(2, 2)} & $< 1$ \\
    Case~(2, 2)-Lg & & $> 1$ \\
    \midrule
  \end{tabular}
\end{table}
We now focus on the Focus-Saddle~(2, 1) class studied in \cite{LLST2022}.
Regarding the Case~(2, 1)-Lg, they showed that $f$ is $C^r$-approximated by a
diffeomorphism $g$ having a heterodimensional cycle involving the continuation
$O^*_{ct}(g)$ of $O^*$ and a new hyperbolic periodic point $Q$ of $g$;
see Figure~\ref{fig-phase_LLST},
where the definition of the heterodimensional cycle is done later,
see before Theorem~\ref{thm-main} for details.
The continuation $O^*_{ct}$ of $O^*$ refers to a $C^r$ map  
from a small neighborhood $\mathcal{U}$ of $f$ in $\mathrm{Diff}^r(M_\mathrm{ph})$  
to a small neighborhood of $O^*$ in $M_\mathrm{ph}$,  
which assigns to each $g \in \mathcal{U}$
the hyperbolic periodic point $O^*_{ct}(g)$ of $g$,  
satisfying $O^*_{ct}(f) = O^*$,
where $\mathrm{Diff}^r(M_\mathrm{ph})$ denotes the set of all $C^r$ diffeomorphisms  
from $M_\mathrm{ph}$ to itself.  
Hereafter, whenever a continuation is naturally determined  
and does not cause confusion,  
we will omit the detailed definition of such a continuation.  

\begin{figure}[h]
  \centering
  \includegraphics{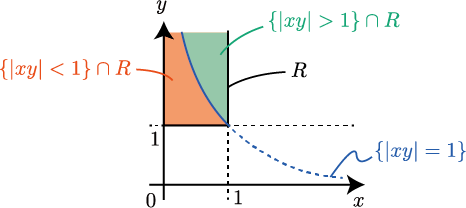}
  \caption{
    The $xy$-plane. The main region is
    $R = \{0 < x < 1, y > 1\}$.
  } \label{fig-lambda_gamma_plane}
\end{figure}

In their result, the assumption $|\lambda^* \gamma^*| > 1$ is essential to create a
hyperbolic periodic point $Q$ whose unstable index is $d^u + 1$.
In the case of $|\lambda^* \gamma^*| < 1$, the unstable index of $Q$ becomes $d^u$,
and so a heterodimensional cycle would not occur.
As a result, the $g$ which has a heterodimensional cycle is also in the region
$\{|\lambda^* \gamma^*| > 1\}$.
Let us explain it more precisely.
We may assume $g$ has a hyperbolic periodic point
$O^*_{ct}(g)$ which is the continuation
of $O^*$ and we can consider the continuations
$\lambda^*_{ct}(g)$ and $\gamma^*_{ct}(g)$
of $\lambda^*$ and $\gamma^*$, respectively.
Consider the region
$R = \{(x, y) \:|\: 0 < x < 1, \, y > 1\}$
in $xy$-plane, see the Figure~\ref{fig-lambda_gamma_plane}.
The pair $(|\lambda^*_{ct}(g)|, |\gamma^*_{ct}(g)|)$ is always in the region
$\{ |xy| > 1\} \cap R$.
This was stated as `$g$ is in the region $\{|\lambda^* \gamma^*| > 1\}$'.

As mentioned before, if the original $f$ is in the region $\{|\lambda^* \gamma^*| < 1\}$,
then we can't find a hyperbolic periodic point whose unstable index is
$d^u + 1$ and we can't find a heterodimensional cycle.
On the other hand, if $f$ is in the curve $\{|\lambda^* \gamma^*| = 1\}$,
then we may find a hyperbolic periodic point whose unstable index is $d^u + 1$.
Our research is to analyze such a $f$ and to extend their result.
Note that $f$ with $|\lambda^* \gamma^*| = 1$ can be $C^r$-approximated by a diffeomorphism
in the region $\{|\lambda^* \gamma^*| > 1\}$,
and hence $f$ is $C^r$-approximated by a diffeomorphism $g$
having a heterodimensional cycle.
Thus, the central question of our study is as follows, and this paper addresses
the following question.
\begin{question}[Central question]\label{q-CQ}
  Let $f \in \{ |\lambda^* \gamma^*| = 1 \}$ be a $C^r$, $r \geq 1$, diffeomorphism.
  Can we get $g$ arbitrarily $C^r$-close to $f$,
  having a heterodimensional cycle so that $g \in \{|\lambda^* \gamma^*| < 1 \}$?
\end{question}\unskip

\medskip

Let us discuss our results.
Suppose that the whole manifold has a dimension three: $\dim M_\mathrm{ph} = 3$.
We denote the multipliers of
a hyperbolic periodic point $O^*$ by
$\lambda_1^*$, $\lambda_2^*$, $\gamma^*$.
Assume
\begin{itemize}
  \item $|\lambda_1^*| = |\lambda_2^*| < 1 < |\gamma^*|$;
  \item The $\lambda_1^*$ and $\lambda_2^*$ are complex conjugate:
    $\lambda_1^* = \lambda^* \mathrm{e}^{i \omega^*}$
    and $\lambda_2^* = \lambda^* \mathrm{e}^{-i \omega^*}$
    for some $\lambda^* \in (0, 1)$ and $\omega^* \in (0, \pi)$;
  \item $|\lambda^* \gamma^*| = 1$.
\end{itemize}
The phase portrait of $f$ is like Figure~\ref{fig-phase_LLST}~(a).
Note that the last assumption makes differences between the inspired paper
\cite{LLST2022} and this paper.
For such a diffeomorphism $f$, we find the following result.
The expanding condition \textbf{(EC)} for a pair $(f, \Gamma)$ is given later,
see Section~\ref{s23-resfam}.
Roughly speaking,
\textbf{(EC)} guarantees that a global map on the orbit of $\Gamma$ has
the area expansion property,
see Remark~\ref{rem-meaningEC} for more details.
In the following,
$\mathrm{u\text{-}index}(X)$ indicates the unstable index of
a hyperbolic periodic orbit $X$ and
$\# Y$ is the number of elements in a finite set $Y$.
We say that diffeomorphism $g$ has a \textit{heterodimensional cycle} involving
two hyperbolic periodic orbits $L_1$ and $L_2$ if
\begin{align*}
  \mathrm{u\text{-}index}(L_1) \neq \mathrm{u\text{-}index}(L_2), \quad
  W^u(L_1) \cap W^s(L_2) \neq \emptyset, \quad
  W^u(L_2) \cap W^s(L_1) \neq \emptyset,
\end{align*}
where $W^u(X)$ and $W^s(X)$ denote the unstable and stable manifolds
of a hyperbolic periodic point, orbit, or set $X$,  
respectively, of the hyperbolic periodic orbit $X$.  
\begin{maintheorem}[Main theorem]\label{thm-main}
  For the above three-dimensional $C^r$, $r \geq 1$, diffeomorphism $f$ with
  $|\lambda^* \gamma^*| = 1$, there exists a $C^r$ diffeomorphism $g$ arbitrarily
  $C^r$-close to $f$ such that $g$ has a heterodimensional cycle involving
  two hyperbolic periodic orbits $L_1$ and $L_2$ of saddles satisfying
  \begin{align*}
    \# L_1 = \# L_2, \quad
    \mathrm{u\text{-}index}(L_1) = 1, \quad
    \mathrm{u\text{-}index}(L_2) = 2.
  \end{align*}
  Moreover, if the pair $(f, \Gamma)$ satisfies the expanding condition \textbf{(EC)},
  then the $g$ can be chosen so that
  \begin{align*}
    |\lambda^*_{ct}(g) \gamma^*_{ct}(g)| < 1,
  \end{align*}
  where $\lambda^*_{ct}(g)$ and $\gamma^*_{ct}(g)$ are the continuations of
  $\lambda^*$ and $\gamma^*$ for $g$, respectively.
\end{maintheorem}

\begin{remark}
  \begin{itemize}
    \item We can create a $C^r$ diffeomorphism $f$ on 3-sphere $S^3$ satisfying
      the assumptions of Theorem~\ref{thm-main} and the expanding condition
      \textbf{(EC)}, so the above theorem gives an affirmative answer
      to Question~\ref{q-CQ}, see the appendix for the construction.
      In fact, the set of diffeomorphism satisfying \textbf{(EC)}
      contains at least an open set in the space of diffeomorphisms having
      the homoclinic tangency, and hence the second half of Theorem~\ref{thm-main}
      can be applied to a lot of diffeomorphisms.
    \item It looks like the first half of the above theorem is the same as
      the main result in \cite{LLST2022}.
      The heterodimensional cycle in their result involves the continuation of $O^*$
      and periodic point $Q$ whose periods basically never coincide.
      Moreover, the heterodimensional cycle we have found  
      is in the Saddle case in the terminology of \cite[Section~2.1]{LT2024},  
      whereas the cycle discovered by them is related to $O^*$  
      and hence does not belong to the Saddle case.        
      Thus, our result differs from theirs in these aspects.
    \item The heterodimensional cycle in our result has coindex one,  
    and hence $f$ can be $C^r$-approximated by a diffeomorphism $g$  
    having a $C^1$-robust heterodimensional dynamics  
    by \cite[Theorem~A]{LT2024},  
    where “$C^1$-robust heterodimensional dynamics” is the terminology  
    defined in that paper.  
    On the other hand, in \cite{LLST2022},  
    they perturb the original $f$ within a generic two-parameter family  
    to obtain a heterodimensional cycle not in the Saddle case,  
    and then show that this cycle is stabilized within the same parameter space.  
    We also perturb the original $f$ within a generic three-parameter family  
    (Theorem~\ref{thm-sec}),  
    but then further perturb the system using another result  
    (Theorem~\ref{thm-tomizawa})  
    to obtain a heterodimensional cycle in the Saddle case.  
    Therefore, it remains an open question  
    whether a heterodimensional cycle in the Saddle case  
    can be obtained within the initial three-parameter family,  
    and if so, whether it can be stabilized.  
    Nevertheless, we conjecture that both questions can be affirmatively answered.  
    
  \end{itemize}
\end{remark}

\subsection{Plan of proof of Theorem~\ref{thm-main}}
\label{s11-plan}
In this section, we give a plan of proof of the main theorem (Theorem~\ref{thm-main}).
It will be reduced to Theorem~\ref{thm-sec}.
First, we review the Hopf bifurcation and related topics to assert Theorem~\ref{thm-sec}.
In this section, we assume $r \geq 4$ unless otherwise noted.
We always allow $r = \infty, \omega$ throughout this paper.

\medskip

Let $g$ be a $C^r$ diffeomorphism having a periodic point $Q$
with period $\mathrm{per}(Q)$,
where the dimension of the whole manifold is greater than or equal to 2:
$\dim M_\mathrm{ph} \geq 2$.
Assume the differential $D(g^{\mathrm{per}(Q)})_Q$ has complex eigenvalues $\nu$ and
$\bar{\nu}$ such that
\begin{align*}
  \nu = \cos\psi + \mathrm{i} \sin\psi,
  \quad \text{and} \quad
  \bar{\nu} = \cos\psi - \mathrm{i} \sin\psi,
\end{align*}
where $\psi \in (0, \pi)$ and $\mathrm{i}$ is the imaginary unit.
By the small perturbation, we may suppose
\begin{itemize}
  \item for any eigenvalue $\tau$ of $D(g^{\mathrm{per}(Q)})_Q$,
    if $\tau$ is different from neither $\nu$ nor $\bar{\nu}$,
    then $|\tau| \neq 1$,
  \item and $\psi \in \Psi_\mathrm{reg}$,
\end{itemize}
where
\begin{align}
  \Psi_\mathrm{reg} :=
    \left\{ \psi \in (0, \pi) \:\middle|\: \psi \notin \frac{2\pi}{j} \mathbb{Z}
    \quad \text{for any} \quad
    j \in \{1, 2, 3, 4\} \right\}
    = (0, \pi) \setminus \left\{ \frac{\pi}{2}, \frac{2\pi}{3} \right\}.
    \label{e12-Preg}
\end{align}
By the assumption and the center manifold theorem \cite{K1967}
and \cite[Section~5A]{HPS1970},
there exist a two-dimensional local center manifold
$W^c_\mathrm{loc}(Q)$ of $Q$.
The smoothness of $W^c_\mathrm{loc}(Q)$ is at least $C^4$  
since $r \geq 4$.  
Note that when $r \in \{\infty, \omega\}$,  
the smoothness does not become $C^r$;  
see, e.g., \cite[Section~5.10.2]{R1998}.  
By \cite[Section~7, 8]{RT1971},
\cite[Chapter~III]{I1979},
\cite[Section~6, 6A]{MM2012},
or \cite[Section~2.8]{D2018},
there exists a neighborhood of $Q$ in $W^c_\mathrm{loc}(Q)$
having $C^4$ complex coordinates $w \in \mathbb{C}$ such that
$g^{\mathrm{per}(Q)}|_{W^c_\mathrm{loc}(Q)}: w \mapsto \tilde{w}$
with $\tilde{w} \in W^c_\mathrm{loc}(Q)$
has the form
\begin{equation}
  \tilde{w}
    = \nu w + \alpha w^2 \bar{w} + O(|w|^4) \label{eq-01-Cform}
\end{equation}
for some constant $\alpha \in \mathbb{C}$,
where $O(|w|^4)$ is a term of fifth order or higher.
From this, we have
\begin{align*}
  |\tilde{w}|
    = |w| \sqrt{1 + 2\Re(\bar{\nu} \alpha)|w|^2 + O(|w|^3)}
    = |w| + \Re(\bar{\nu} \alpha) |w|^3 + O(|w|^4),
\end{align*}
where $O(|w|^n)$ is a term of $n$-th order or higher
for any $n \geq 1$
and $\Re(X)$ denotes the real part of the complex number $X$.
This implies when
\begin{align*}
  \mathrm{LC}(Q) = \mathrm{LC}(Q; w)
    := -\Re(\bar{\nu} \alpha)
\end{align*}
is negative,
$Q$ is weakly repelling on $W^c_\mathrm{loc}(Q)$,
and when $\mathrm{LC}(Q)$ is positive, $Q$ is attracting on $W^c_\mathrm{loc}(Q)$.
Therefore, the sign of $\mathrm{LC}(Q)$ is determined independently
of the way the coordinates giving the canonical form \eqref{eq-01-Cform} are taken.
We call $\mathrm{LC}(Q)$ the \textit{first Lyapunov coefficient},
or simply the \textit{Lyapunov coefficient} of a generic point $Q$.
\begin{definition}[Generic Hopf point]\label{def-Hopf_point}
  We say that $Q$ is a generic Hopf point of a $C^r$, $r \geq 4$,
  diffeomorphism $g$ if the Lyapunov coefficient $\mathrm{LC}(Q)$ is not zero:
  $\mathrm{LC}(Q) \neq 0$.
\end{definition}

\medskip

\begin{figure}[h]
  \centering
  \includegraphics{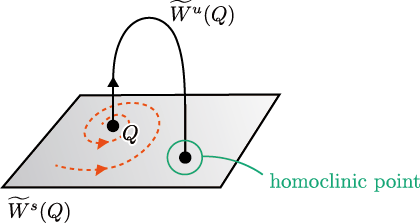}
  \caption{
    The phase portrait of $g$ when $\dim M_\mathrm{ph} = 3$ and $\mathrm{LC}(Q) > 0$.
    In this setting, $\widetilde{W}^s(Q)$ contains $W^c_\mathrm{loc}(Q)$ and
    $\widetilde{W}^u(Q)$ is one-dimensional manifold.
  } \label{fig-HHC}
\end{figure}

Assume $Q$ is a generic Hopf point.
We define
\begin{align}
\begin{aligned}
  \widetilde{W}^s(Q)
    &:= \{ M \in M_\mathrm{ph} \:|\:
    \lim_{n \to \infty} \mathrm{dist}(f^n(M), f^n(Q)) = 0 \}, \\
  \widetilde{W}^u(Q)
    &:= \{ M \in M_\mathrm{ph} \:|\:
    \lim_{n \to \infty} \mathrm{dist}(f^{-n}(M), f^{-n}(Q)) = 0 \},
\end{aligned} \label{e12-gmfds}
\end{align}
where the $\mathrm{dist}$ is the metric that defines the same topology as $M_\mathrm{ph}$.
Though $Q$ is non-hyperbolic,
$\widetilde{W}^s(Q)$ and $\widetilde{W}^u(Q)$ are immersed submanifolds
since $Q$ is determined to attract or repel on the local central manifold
$W^c_\mathrm{loc}(Q)$.
We call $\widetilde{W}^s(Q)$ and $\widetilde{W}^u(Q)$  
as the generalized stable manifold and  
generalized unstable manifold of $Q$, respectively.  
\begin{definition}[Hopf-homoclinic cycle]\label{def-HHC}
  We say that $C^r$, $r \geq 4$, diffeomorphism $g$ has a Hopf-homoclinic cycle
  of a generic Hopf point $Q$ if
  \begin{align*}
    \left( \widetilde{W}^s(Q) \cap \widetilde{W}^u(Q) \right) \setminus \{Q\}
      \neq \emptyset.
  \end{align*}
\end{definition}
See Figure~\ref{fig-HHC} to understand how the cycle looks.
If $g$ has a Hopf-homoclinic cycle,
then $\dim M_\mathrm{ph} \geq 3$,
since $\dim \widetilde{W}^s(Q)$ and $\dim \widetilde{W}^u(Q)$
are greater than or equal to 1
and $\dim \widetilde{W}^s(Q)$ or $\dim \widetilde{W}^u(Q)$
is greater than or equal to 2.

\medskip

Now, we assert our secondary theorem.
In the following theorem, we assume $r \geq 5$  
in order to ensure the boundedness of the partial derivatives up to third order  
of the functions $q_k^{(i)}$, $i \in \{1, 2, 3\}$,  
which appear in Section~\ref{s312-r}.  

\begin{maintheorem}[Secondary theorem]\label{thm-sec}
  For the three-dimensional $C^r$, $r \geq 5$, diffeomorphism $f$
  in Theorem~\ref{thm-main} with $|\lambda^* \gamma^*| = 1$,
  there exists a $C^r$ diffeomorphism $g$ arbitrarily $C^r$-close to $f$
  such that $g$ has a Hopf-homoclinic cycle
  of a generic Hopf point with a negative Lyapunov coefficient.
  Moreover, if the pair $(f, \Gamma)$ satisfies the expanding condition \textbf{(EC)},
  then the $g$ can be chosen so that
  \begin{align*}
    |\lambda^*_{ct}(g) \gamma^*_{ct}(g)| < 1.
  \end{align*}
\end{maintheorem}

The above theorem implies our main result (Theorem~\ref{thm-main}),
by using the following result.
In the following theorem, the assumption $r \geq 5$ is made  
to ensure the existence of an invariant circle after the Hopf bifurcation.  

\begin{theorem}[Three-dimensional version of Theorem~1.1 in \cite{T2019}]
\label{thm-tomizawa}
  Let $g$ be a $C^r$, $r \geq 5$, diffeomorphism on a manifold $M_\mathrm{ph}$
  with $\dim M_\mathrm{ph} = 3$, having a Hopf-homoclinic cycle
  of a generic Hopf point.
  Then there exists a $C^r$ diffeomorphism $g'$ arbitrarily $C^r$-close to $g$
  such that $g'$ has a heterodimensional cycle involving
  two hyperbolic periodic orbits $L_1$ and $L_2$ of saddles satisfying
  \begin{align*}
    \# L_1 = \# L_2, \quad
    \mathrm{u\text{-}index}(L_1) = 1,
    \quad \text{and} \quad
    \mathrm{u\text{-}index}(L_2) = 2.
  \end{align*}
\end{theorem}

\subsubsection{Idea of proof of Theorem~\ref{thm-sec}}
Now, our main objective is reduced to prove Theorem~\ref{thm-sec}.
Let us explain the idea of proof of this theorem.

\medskip

First, we can find the periodic points $Q$ as shown in Figure~\ref{fig-phase_LLST}~(b).
Since we assumed $|\lambda^* \gamma^*| = 1$, $Q$ would be non-hyperbolic.
The existence of $Q$ is achieved in Proposition~\ref{p34-Qk}.

Second, we compute the Lyapunov coefficient in detail
and observe in Proposition~\ref{p42-LC} that $Q$ is actually a generic Hopf point.
In fact, the Lyapunov coefficient of $Q$ can be taken to be always negative.
In other words, $Q$ is always weakly repelling, see Proposition~\ref{p42-LC}.
Therefore, we can also obtain a situation like the one shown on the right side
of Figure~\ref{fig-phase_LLST}~(b) in our settings.

Finally, since $\widetilde{W}^u(Q)$ rotates and approaches $W^u(O^*_{ct}(g))$,
a small perturbation allows a Hopf-homoclinic cycle of $Q$ to be found.
This observation is confirmed by Proposition~\ref{p51-Tint} and \ref{p52-HHC}.
This is a summary of the proof of Theorem~\ref{thm-sec}.

\medskip

The construction of this paper is as follows.
In Section~\ref{s2-fam}, we state Theorem~\ref{thm-third},
which is a detailed version of Theorem~\ref{thm-sec}
using the term parameter family.
Hence, our goal will be to prove Theorem~\ref{thm-third}.
In Section~\ref{s3-nhyp}, we find the non-hyperbolic periodic point $Q$
that is mentioned the above.
In Section~\ref{s4-Hopf}, we verify that $Q$ is a generic Hopf point and
$Q$ can be taken to be always weakly repelling.
In Section~\ref{s5-HHC}, we find the Hopf-homoclinic cycle of $Q$ and
complete the proof of Theorem~\ref{thm-third}.

\ifdefined\isMaster
\else
  \end{document}
\fi
\ifdefined\isMaster
\else
  \documentclass[11pt,a4paper]{article}

\usepackage[left=2cm,right=2cm,top=2cm,bottom=2cm]{geometry}
\usepackage{graphicx}
\usepackage{multirow}
\usepackage{amsmath,amssymb,amsfonts}
\usepackage{amsthm}
\usepackage{mathrsfs}
\usepackage[title]{appendix}
\usepackage{xcolor}
\usepackage{textcomp}
\usepackage{manyfoot}
\usepackage{booktabs}
\usepackage{algorithm}
\usepackage{algorithmicx}
\usepackage{algpseudocode}
\usepackage{listings}
\ifdefined\debug
  \usepackage[inline]{showlabels}
  \showlabels{ref}
  \showlabels{cite}
  \usepackage{fancyhdr}
\fi
\usepackage{bm}
\usepackage{hyperref}

\theoremstyle{plain}
\newtheorem{theorem}{Theorem}[section]
\newtheorem{proposition}[theorem]{Proposition}%
\newtheorem{lemma}[theorem]{Lemma}
\newtheorem{corollary}[theorem]{Corollary}
\newtheorem{claim}[theorem]{Claim}
\theoremstyle{definition}
\newtheorem{definition}[theorem]{Definition}
\newtheorem{remark}[theorem]{Remark}
\newtheorem{example}[theorem]{Example}
\newtheorem{question}[theorem]{Question}
\newtheorem*{notation}{Notation}
\theoremstyle{plain}
\newtheorem{maintheorem}{Theorem}
\renewcommand{\themaintheorem}{\Alph{maintheorem}}
\newcommand{\wiptext}[1]{\textcolor{blue}{#1}}
\ifdefined\debug
  \newcommand{\debugtext}[1]{\textcolor{red}{#1}}
  \pagestyle{fancy} 
  \fancyhead[L]{}
  \fancyhead[C]{}
  \fancyhead[R]{\textcolor{red}{[DEBUG] Debug mode is activated.}}
  \renewcommand{\headrulewidth}{0pt}
\else
  \newcommand{\debugtext}[1]{\ignorespaces}
\fi

\counterwithin{figure}{section}
\counterwithin{table}{section}

\raggedbottom

\numberwithin{equation}{section}

  \begin{document}
\fi

\section{Three-parameter family of diffeomorphisms}\label{s2-fam}
Our goal is now to prove Theorem~\ref{thm-sec},
which will be reduced to Theorem~\ref{thm-third}.
This theorem is given in Section~\ref{s23-resfam}
by using terms of a parameter family.
In Section~\ref{s21-PU},
we define such a three-parameter family of diffeomorphisms.
In Section~\ref{s23-resfam},
we give several conditions on a pair $(f, \Gamma)$ to assert Theorem~\ref{thm-third},
where $\Gamma$ is an orbit of a homoclinic tangency
and state Theorem~\ref{thm-third}.

In the remaining sections below,
except where explicitly stated,
we always assume $f$ is a three-dimensional $C^r$, $r \geq 5$, diffeomorphism
in Theorem~\ref{thm-sec}.

\subsection{Proper Unfolding}\label{s21-PU}
In this section, we will define the three parameters $\mu$, $\omega$, and $\rho$,
and define proper unfolding,
which is the terminology in \cite{LLST2022}.
Some preparations are made before giving the definition.

\medskip

Let $\mathrm{per}(O^*)$ be the period of $O^*$ and define the local map
\begin{align}
  T_0^*
  = T_0^*(f, U_0)
  := f^{\mathrm{per}(O^*)}|_{U_0 \cap f^{-\mathrm{per}(O^*)}(U_0)}
  \label{eq-def_T00}
\end{align}
on a small neighborhood $U_0$ of $O^*$.
We define the local stable manifold
$W^s_\mathrm{loc}(O^*) = W^s_\mathrm{loc}(O^*; f, U_0)$
of $O^*$ by the connected component
of $W^s(O^*) \cap U_0$ that contains $O^*$.
The local unstable manifold $W^u_\mathrm{loc}(O^*) = W^u_\mathrm{loc}(O^*; f, U_0)$
is defined in the same way.
Pick two base points $M^-_0 \in W^u_\mathrm{loc}(O^*) \cap \Gamma$ and
$M^+_0 \in W^s_\mathrm{loc}(O^*) \cap \Gamma$.
There is $n_0 = n_0(M^-_0, M^+_0) \in \mathbb{Z}_{>0}$ such that
$(f^{\mathrm{per}(O^*)})^{n_0}(M^-_0) = M^+_0$.
We define the global map by
\begin{align}
  T_1^*
  = T_1^*(f, \Gamma, U_0, M^-_0, M^+_0)
  := (f^{\mathrm{per}(O^*)})^{n_0}.
  \label{eq-02-def_T10}
\end{align}
By the assumption,
the image $T_1^*(W^u_\mathrm{loc}(O^*))$ is tangent to $W^s_\mathrm{loc}(O^*)$
at $M^+_0$.

We define the quadratic condition by
\begin{itemize}
  \item \textbf{(QC)} The $(f, \Gamma)$ satisfies \textbf{(QC)}
    if there exist $U_0$,
    $M^-_0$, and $M^+_0$ such that the tangency between
    $T_1^*(W^u_\mathrm{loc}(O^*))$ and
    $W^s_\mathrm{loc}(O^*)$ at $M^+_0$ is quadratic.
\end{itemize}
Here, we say that the tangency $M$ between
embedded $C^r$ submanifolds $\mathcal{N}^u$,
$\mathcal{N}^s \subset M_\mathrm{ph}$ with $\dim \mathcal{N}^u = 1$ and
$\dim \mathcal{N}^s = 2$ is quadratic
if there exists a small neighborhood of $M$ having
$C^r$ coordinates $(u_1, u_2, v)$ such that
\begin{align*}
  \mathcal{N}^s = \{ v = 0 \} \quad \text{and} \quad
  \mathcal{N}^u = \{ v = h(u_1), u_2 = 0 \}
\end{align*}
for some $C^r$ function $h$ satisfying $h(0) = 0$, $h'(0) = 0$, and $h''(0) \neq 0$.
\begin{remark}
  \begin{itemize}
    \item We can verify
      if $(f, \Gamma)$ holds \textbf{(QC)}, then the tangency between
      $T_1^*(W^u_\mathrm{loc}(O^*))$ and $W^s_\mathrm{loc}(O^*)$
      is quadratic for any $U_0$,
      $M^-_0$, and $M^+_0$ since $f$ is a diffeomorphism.
    \item Even if $(f, \Gamma)$ does not satisfy \textbf{(QC)}, then there exist a
      diffeomorphism $g$ arbitrarily $C^r$-close to $f$ and an orbit $\Gamma'$ of a
      homoclinic tangency to the continuation $O(g)$ of $O^*$
      such that $(g, \Gamma')$ holds
      \textbf{(QC)}, see \cite{NPT1983}.
  \end{itemize}
\end{remark}

There are $C^r$ coordinates $(s_1, s_2, t)$ on $U_0$ such that
\begin{equation}
  W^s_\mathrm{loc}(O^*) = \{ t = 0 \}, \quad
  W^u_\mathrm{loc}(O^*) = \{ s_1 = 0,\, s_2 = 0\}
  \label{eq-02-CCWsWu}
\end{equation}
and $T_0^*: (s_1, s_2, t)
\mapsto (\hat{s}_1, \hat{s}_2, \hat{t})$ has the form
\begin{equation}
  \begin{aligned}
    \hat{s}_1
      &= \lambda^* s_1 \cos\omega^*
      - \lambda^* s_2 \sin\omega^*
      + p_1^*(s_1, s_2, t), \\
    \hat{s}_2
      &= \lambda^* s_1 \sin\omega^*
      + \lambda^* s_2 \cos\omega^*
      + p_2^*(s_1, s_2, t), \\
    \hat{t}
      &= \gamma^* t + p_3^*(s_1, s_2, t),
  \end{aligned}\label{eq-02-CCDiag}
\end{equation}
where $p_1^*$, $p_2^*$ and $p_3^*$ are $C^r$ maps with
\begin{equation}
  \begin{aligned}
    p_i^*(0, 0, 0) = 0, \quad
    p_{i, s_j}^*(0, 0, 0), \,
    p_{i, t}^*(0, 0, 0) = 0, \quad
    p_l^*(0, 0, t) \equiv 0, \quad
    p_3^*(s_1, s_2, 0) \equiv 0
  \end{aligned}\label{eq-02-CCpcondi}
\end{equation}
for any $i \in \{1, 2, 3\}$ and $j$, $l \in \{1, 2\}$.
Using these coordinates, we write the two base points $M^-_0$ and $M^+_0$ by
\begin{align}
  M^-_0 = (0, 0, t^-), \quad
  M^+_0 = (s^+_1, s^+_2, 0).
  \label{eq-02-Mp0Mm0}
\end{align}
By using the coordinates $(s_1, s_2, t)$, consider the small cube
\begin{align}
  \Pi^-_0 := [-\delta^-_0, \delta^-_0]^3 + M^-_0,
  \quad \delta^-_0 > 0,
  \label{eq-02-Pim0}
\end{align}
centered at $M^-_0$,
where $X + a = \{x + a \:|\: x \in X\}$
for any $X \subset \mathbb{R}^n$, $a \in \mathbb{R}^n$, and $n \geq 1$.

\medskip

Now, we define the parameters.
Assume $(f, \Gamma)$ holds the quadratic condition \textbf{(QC)}.
For any $g$ sufficiently $C^r$-close to $f$,
there is the continuation
\begin{align*}
  W^u(g) := W^u_\mathrm{loc}(O^*_{ct}(g)) \cap \Pi^-_0
\end{align*}
of the segment
\begin{align*}
  W^u_0 := W^u_\mathrm{loc}(O^*) \cap \Pi^-_0,
\end{align*}
where $O^*_{ct}(g)$ is the continuation of $O^*$ and we define
$W^\sigma_\mathrm{loc}(O^*_{ct}(g))$, $\sigma \in \{s, u\}$,
by the connected component of
$W^u(O^*_{ct}(g)) \cap U_0$
that contains $O^*_{ct}(g)$.
Let $\mu$ be the $C^r$ functional
from a small neighborhood $\mathcal{U}$ of $f$ to $\mathbb{R}$ such that
$\mu(g)$ assigns signed distance between
$(g^{\mathrm{per}(O^*)})^{n_0}(W^u(g))$
and $W^s_\mathrm{loc}(O^*_{ct}(g))$ for any $g \in \mathcal{U}$,
where the direction of the sign is arbitrary.

Recall that we wrote the argument of the stable multiplier $\lambda_1^*$
of $O^*$ by $\omega^* \in (0, \pi)$.
The $\omega$ is the $C^{r - 1}$ functional
from a small neighborhood of $f$ to $(0, \pi)$ such that
$\omega$ is the continuation of $\omega^*$.
The reason why the smoothness of $\omega$ is $r - 1$ is that,
in general, the eigenvalues are solutions of equations
with first partial derivatives as coefficients,
so the smoothness is one lower.

Remember that $\lambda^* = |\lambda_1^*|$
and $\gamma^* \in \mathbb{R}$ is the unstable multiplier.
We define the $C^{r - 1}$ functional $\rho$ by
\begin{align}
  \rho(g) := \log | \lambda^*_{ct}(g) \gamma^*_{ct}(g) | \label{eq-02-rhoD}
\end{align}
for any $g$ sufficiently $C^r$-close to $f$,
where $\lambda^*_{ct}(g)$ and $\gamma^*_{ct}(g)$ are the continuations of
$\lambda^*$ and $\gamma^*$, respectively.
The reason why the smoothness is $r - 1$ is the same as for $\omega$,
and the reason for taking the logarithm is to ensure that
$\rho(f) = \log | \lambda^* \gamma^* | = 0$.

\medskip

Let $\{ f_\varepsilon \}_{\varepsilon \in R_\mathrm{prm}^*}$
be a three-parameter family of $C^r$ diffeomorphisms
with $f_{\varepsilon^*} = f$,
where we assume $\varepsilon$ runs in three-dimensional open ball
$R_\mathrm{prm}^* \subset \mathbb{R}^3$ centered at $\varepsilon^*$.
We always assume that the smoothness with respect to the parameters is also $C^r$;
the $f:M_\mathrm{ph} \times R_\mathrm{prm}^*
\ni (M, \varepsilon) \mapsto f_\varepsilon(M) \in M_\mathrm{ph}$ is $C^r$ and
$f_\varepsilon: M_\mathrm{ph} \to M_\mathrm{ph}$
is a $C^r$ diffeomorphism for any $\varepsilon \in R_\mathrm{prm}^*$.

\begin{definition}[Proper unfolding]\label{def-proper_unfolding}
  We say that $\{ f_\varepsilon \}_{\varepsilon \in R_\mathrm{prm}^*}$
  \textit{unfolds properly}
  at $\varepsilon = \varepsilon^*$ with respect to $\Gamma$
  (or simply that $\{ f_\varepsilon \}_{\varepsilon \in R_\mathrm{prm}^*}$
  \textit{unfolds properly}) if
  \begin{align*}
    \mathrm{det} \,
      \frac{\partial (\mu(f_\varepsilon), \omega(f_\varepsilon), \rho(f_\varepsilon))}
      {\partial \varepsilon}(\varepsilon^*) \neq 0,
  \end{align*}
  where the expression inside $\det$ is the $3 \times 3$ Jacobian matrix.  
\end{definition}

\begin{remark}\label{rmk-param}
  For a proper unfolding family  
  $\{ f_\varepsilon \}_{\varepsilon \in R_\mathrm{prm}^*}$,  
  the inverse function theorem guarantees that $\varepsilon$  
  and $(\mu, \omega, \rho)$ correspond one-to-one  
  via some $C^{r - 1}$ map  
  \begin{align*}
    R_\mathrm{prm}^* \ni \varepsilon  
    \mapsto (\mu(f_\varepsilon), \omega(f_\varepsilon), \rho(f_\varepsilon))  
  \end{align*}  
  and its inverse,  
  by replacing $R_\mathrm{prm}^*$ with a smaller one if necessary.  
  In the following, we identify $\varepsilon$ with $(\mu, \omega, \rho)$  
  via the above map.  
  Thus, we write $\varepsilon = (\mu, \omega, \rho)$  
  and $\varepsilon^* = (0, \omega^*, 0)$.    
\end{remark}

\subsection{Our result in a three-parameter family}
\label{s23-resfam}
In this section,
we first define the accompanying condition \textbf{(AC)}.
Next, we state the expanding condition \textbf{(EC)},
and then describe Theorem~\ref{thm-third}.

\medskip

Recall the coordinates $(s_1, s_2, t)$ satisfying
\eqref{eq-02-CCWsWu} and \eqref{eq-02-CCDiag} with \eqref{eq-02-CCpcondi},
and the point $t^-$ in \eqref{eq-02-Mp0Mm0}.
We put the pair of $U_0$ and $(s_1, s_2, t)$ by
\begin{align*}
  \mathbb{U}_0^* := (U_0; s_1, s_2, t).
\end{align*}
We define the accompanying condition as follows.
\begin{itemize}
  \item \textbf{(AC)} The $(f, \Gamma)$ satisfies \textbf{(AC)}
  if there exist $\mathbb{U}_0^*$,
  $M^-_0$, points of transverse intersection $\{(0, 0, t_i)\}$
  between $W^s(O^*)$ and $W^u_\mathrm{loc}(O^*)$
  such that $\{t_i\}$ converges $t^-$
  from the both sides as $i \to \infty$.
\end{itemize}
\begin{remark}
  \begin{itemize}
    \item Although the situation of \textbf{(AC)} does not seem to occur in general,
      but it happen all the time,
      see Proposition~\ref{p23-AC},
      where the proof is completely based on \cite{LLST2022}.
    \item Whether \textbf{(AC)} is satisfied or not
      does not depend on the choice of
      $U_0$, the coordinates $(s_1, s_2, t)$, and $M^-_0$
      due to the invariance of $W^s(O^*)$ and $W^u_\mathrm{loc}(O^*)$.
    \item Let $\theta_i \in (0, \pi/2]$ be the angle between $W^s(O^*)$ and
      $W^u_\mathrm{loc}(O^*)$ at $(0, 0, t_i)$.
      In general,
      $\theta_i \to 0$ as $i \to \infty$,
      and so even if there is a
      one-dimensional $C^r$ disks $\{\ell_k\}$ $C^1$-converging to the small
      neighborhood of $t^-$ in $W^u_\mathrm{loc}(O^*)$ as $k \to \infty$,
      then we may not be able to find the intersection between $\ell_k$ and $W^s(O^*)$
      when the length of $\ell_k$ converges to 0 as $k \to \infty$.
      On the other hand, if the length of $\ell_k$ is bounded away from 0,
      then we can find the transversal intersection.
      A similar observation will be used in the proof of Proposition~\ref{p51-Tint}
      in Section~\ref{s51-TInt}.
  \end{itemize}
\end{remark}

Recall the neighborhood $\Pi^-_0$ of $M^-_0$ in \eqref{eq-02-Pim0}.
Replacing $\delta_0^- > 0$ with a smaller one if necessary,
we may suppose $T_1^*(\Pi_0^-) \subset U_0$.
Using the coordinates $(s_1, s_2, t)$,
we express the global map
$T_1^*:\Pi^-_0 \ni (\tilde{s}_1, \tilde{s}_2, \tilde{t}) \mapsto
(\bar{s}_1, \bar{s}_2, \bar{t}) \in U_0$ as
\begin{equation}
  \begin{aligned}
    \bar{s}_1 - s_1^+
      &= a_{11}^* \tilde{s}_1 + a_{12}^* \tilde{s}_2
      + b_1^* (\tilde{t} - t^-)
      + O(\|(\tilde{s}_1, \tilde{s}_2, \tilde{t} - t^-)\|^2), \\
    \bar{s}_2 - s_2^+
      &= a_{21}^* \tilde{s}_1 + a_{22}^* \tilde{s}_2
      + b_2^* (\tilde{t} - t^-)
      + O(\|(\tilde{s}_1, \tilde{s}_2, \tilde{t} - t^-)\|^2), \\
    \bar{t}
      &= c_1^* \tilde{s}_1 + c_2^* \tilde{s}_2
      + d^* (\tilde{t} - t^-)^2
      + O(\|(\tilde{s}_1, \tilde{s}_2, \tilde{t} - t^-)\|^2),
  \end{aligned} \label{eq-02-ExpT10}
\end{equation}
where
$O(\cdot)$ are terms of second order or higher of the Taylor expansion, excluding
the explicitly stated terms.
Note that the coefficient of $(\tilde{t} - t^-)$ in $\bar{t}$ vanishes
since $M^-_0$, $M^+_0 \in \Gamma$
and the quadratic condition \textbf{(QC)} says $d^* \neq 0$.
Note also that the above coefficients depend on
$f$, $\Gamma$, $\mathbb{U}_0^* = (U_0; s_1, s_2, t)$,
$M^-_0$, and $M^+_0$ by \eqref{eq-02-def_T10}:
\begin{align*}
  &a_{ij}^* = a_{ij}^*(f, \Gamma, \mathbb{U}_0^*, M^-_0, M^+_0), \quad
  b_i^* = b_i^*(f, \Gamma, \mathbb{U}_0^*, M^-_0, M^+_0), \\
  &c_i^* = c_i^*(f, \Gamma, \mathbb{U}_0^*, M^-_0, M^+_0), \quad
  d^* = d^*(f, \Gamma, \mathbb{U}_0^*, M^-_0, M^+_0)
\end{align*}
for any $i$, $j \in \{1, 2\}$.

We consider the following quantity:
\begin{align}
  \mathcal{E}
    = \mathcal{E}(f, \Gamma, \mathbb{U}_0^*, M^-_0, M^+_0)
    := \sqrt{(b_1^*)^2 + (b_2^*)^2} \sqrt{(c_1^*)^2 + (c_2^*)^2}.
  \label{eq-02-mathcalE}
\end{align}
We define the expanding condition \textbf{(EC)} as follows.
\begin{itemize}
  \item \textbf{(EC)} The $(f, \Gamma)$ satisfies \textbf{(EC)}
    if there exist $\mathbb{U}_0^*$, $M^-_0$, and $M^+_0$ such that
    \begin{align*}
      \mathcal{E}(f, \Gamma, \mathbb{U}_0^*, M^-_0, M^+_0) > 1.
    \end{align*}
\end{itemize}

\begin{remark} \label{rem-meaningEC}
  For each point $M \in U_0$,
  define a basis of the tangent space at $M$, denoted $T_M U_0$, 
  by taking the natural basis associated with the coordinate system $(s_1, s_2, t)$,
  and denote the basis vectors 
  by $e_M^{(1)}$, $e_M^{(2)}$, and $e_M^{(3)}$.
  For any $v = v_1 e^{(1)}_M + v_2 e^{(2)}_M + v_3 e^{(3)}_M \in T_M U_0$
  ($v_1$, $v_2$, $v_3 \in \mathbb{R}$),
  we define
  \begin{align*}
    \| v \|_0 = \sqrt{v_1^2 + v_2^2 + v_3^2}
  \end{align*}
  and $\mathrm{pr}_M$ is the projection defined by
  \begin{align*}
    \mathrm{pr}_M (v) = v_3 e^{(3)}_M.
  \end{align*}
  The notation $\mathrm{span} \, X$ denotes the space spanned by
  the elements of the subset $X$ of a vector space.
  Since \eqref{eq-02-ExpT10},
  the geometric meaning of $\mathcal{E}$ is as
  \begin{align*}
    \mathcal{E}
    &= \max_{v, w} \left\|
      \left( \mathrm{pr}_{M_0^+} \circ D(T_1^*)_{M_0^-}(v) \right)
      \times D(T_1^*)_{M_0^-}(w) \right\|_0 \\
    &= \max_{v, w} \, \left\{
      \text{The area of the rectangle spanned by }
      \mathrm{pr}_{M_0^+} \circ D(T_1^*)_{M_0^-}(v)
        \text{ and } D(T_1^*)_{M_0^-}(w)
    \right\},
  \end{align*}
  where $v$ and $w$ are assumed to move while satisfying
  \begin{align*}
    v \in \mathrm{span} \, \left\{e^{(1)}_{M_0^-}, e^{(2)}_{M_0^-}\right\}, \quad
      \: \| v \|_0 = 1, \quad
    w \in \mathrm{span} \, \left\{e^{(3)}_{M_0^-}\right\}, \quad
      \| w \|_0 = 1.
  \end{align*}
  Thus, \textbf{(EC)} states that the global map is area expanding
  in the above sense.
\end{remark}
Note that if $(f, \Gamma)$ holds \textbf{(EC)},
then $\mathcal{E}(f, \Gamma, \mathbb{U}_0^*, M^-_0, M^+_0) > 1$
for any $U_0$, $M_0^-$, and $M_0^+$
since on $\{s_1 = 0, \, s_2 = 0\} \cup \{t = 0\}$,
$DT_1^*$ is $\lambda^*$-contracting in the $(s_1, s_2)$-direction
and $\gamma^*$-expanding in the $t$-direction and $|\lambda^* \gamma^*| = 1$.
In fact,
it does not depend on how the coordinates $(s_1, s_2, t)$ are taken:
\begin{proposition}[Independence of \textbf{(EC)}] \label{p23-IndEC}
  If $(f, \Gamma)$ holds \textbf{(EC)},
  then $\mathcal{E}(f, \Gamma, \mathbb{U}_0^*, M^-_0, M^+_0) > 1$
  for any $\mathbb{U}_0^*$, $M_0^-$, and $M_0^+$.
\end{proposition}

\begin{remark}
  From the above lemma,
  to check if \textbf{(EC)} is satisfied,
  we just verify that $\mathcal{E} > 1$ at some coordinates
  and base points,
  where we take the coordinates so that they hold
  \eqref{eq-02-CCWsWu}--\eqref{eq-02-CCpcondi}.
\end{remark}

\begin{proof}[Proof of Proposition~\ref{p23-IndEC}]
  See the appendix.
\end{proof}

\medskip

Recall the quadratic condition \textbf{(QC)}
and the accompanying condition \textbf{(AC)}
defined in this and the previous sections.
Theorem~\ref{thm-sec} is reduced to the following theorem.
\begin{maintheorem}[Third theorem] \label{thm-third}
  Suppose $(f, \Gamma)$ satisfies \textbf{(QC)} and \textbf{(AC)}.
  For any proper unfolding three-parameter family
  $\{ f_\varepsilon \}_{\varepsilon \in R_\mathrm{prm}^*}$ of
  $C^r$ diffeomorphisms with $f_{\varepsilon^*} = f$,
  there exists a sequence $\{ \varepsilon_k \}$ in $R_\mathrm{prm}^*$
  converging to $\varepsilon^*$
  such that $f_{\varepsilon_k}$ has a generic Hopf point $Q_k$
  with a negative Lyapunov coefficient
  and $Q_k$ has a Hopf-homoclinic cycle.
  Moreover,
  if the original $(f, \Gamma)$ holds \textbf{(EC)},
  then we can take the sequence $\{ \varepsilon_k \}$ so that
  $\rho|_{\varepsilon = \varepsilon_k} < 0$
  for all $k$.
\end{maintheorem}

Theorem~\ref{thm-sec} follows immediately from
Theorem~\ref{thm-third} and the following proposition.
The proof of the following proposition can be found in \cite{LLST2022}.
\begin{proposition}[Generality of \textbf{(AC)}]
  \label{p23-AC}
  Suppose $(f, \Gamma)$ satisfies \textbf{(QC)}.
  For any proper unfolding three-parameter family
  $\{ f_\varepsilon \}_{\varepsilon \in R_\mathrm{prm}^*}$ of
  $C^r$ diffeomorphisms with $f_{\varepsilon^*} = f$,
  there exist sequences
  $\{ \mu_j \}$, $\{ \omega_j \}$
  converging to $0$, $\omega^*$, respectively, such that
  \begin{itemize}
    \item the $f_{(\mu_j, \omega_j, 0)}$ has an orbit $\Gamma_j$ of
      a homoclinic tangency to the continuation $O(\mu_j, \omega_j, 0)$ of
      $O^*$ for any $j$,
    \item the pair $(f_{(\mu_j, \omega_j, 0)}, \Gamma_j)$ satisfies
      \textbf{(QC)} and \textbf{(AC)}, and
    \item the $\{ f_\varepsilon \}_{\varepsilon \in R_\mathrm{prm}^*}$
      unfolds properly at
      $\varepsilon = (\mu_j, \omega_j, 0)$ with respect to $\Gamma_j$.
  \end{itemize}
\end{proposition}

Now, our main goal is to prove Theorem~\ref{thm-third}.
In the following sections, we will focus on the proof.

\ifdefined\isMaster
\else
  \end{document}
\fi
\ifdefined\isMaster
\else
  \documentclass[11pt,a4paper]{article}

\usepackage[left=2cm,right=2cm,top=2cm,bottom=2cm]{geometry}
\usepackage{graphicx}
\usepackage{multirow}
\usepackage{amsmath,amssymb,amsfonts}
\usepackage{amsthm}
\usepackage{mathrsfs}
\usepackage[title]{appendix}
\usepackage{xcolor}
\usepackage{textcomp}
\usepackage{manyfoot}
\usepackage{booktabs}
\usepackage{algorithm}
\usepackage{algorithmicx}
\usepackage{algpseudocode}
\usepackage{listings}
\ifdefined\debug
  \usepackage[inline]{showlabels}
  \showlabels{ref}
  \showlabels{cite}
  \usepackage{fancyhdr}
\fi
\usepackage{bm}
\usepackage{hyperref}

\theoremstyle{plain}
\newtheorem{theorem}{Theorem}[section]
\newtheorem{proposition}[theorem]{Proposition}%
\newtheorem{lemma}[theorem]{Lemma}
\newtheorem{corollary}[theorem]{Corollary}
\newtheorem{claim}[theorem]{Claim}
\theoremstyle{definition}
\newtheorem{definition}[theorem]{Definition}
\newtheorem{remark}[theorem]{Remark}
\newtheorem{example}[theorem]{Example}
\newtheorem{question}[theorem]{Question}
\newtheorem*{notation}{Notation}
\theoremstyle{plain}
\newtheorem{maintheorem}{Theorem}
\renewcommand{\themaintheorem}{\Alph{maintheorem}}
\newcommand{\wiptext}[1]{\textcolor{blue}{#1}}
\ifdefined\debug
  \newcommand{\debugtext}[1]{\textcolor{red}{#1}}
  \pagestyle{fancy} 
  \fancyhead[L]{}
  \fancyhead[C]{}
  \fancyhead[R]{\textcolor{red}{[DEBUG] Debug mode is activated.}}
  \renewcommand{\headrulewidth}{0pt}
\else
  \newcommand{\debugtext}[1]{\ignorespaces}
\fi

\counterwithin{figure}{section}
\counterwithin{table}{section}

\raggedbottom

\numberwithin{equation}{section}

  \begin{document}
\fi

\section{Existence of non-hyperbolic periodic points}
\label{s3-nhyp}
This section aims to prove half of Theorem~\ref{thm-third}, 
specifically the existence of the non-hyperbolic fixed point $Q_k$ of
the so-called first-return map $T_k$.
It is accomplished by Proposition~\ref{p34-Qk} in Section~\ref{s34-Qk}.
In Section~\ref{s31-FRM} and \ref{s32-NF},
we define the first-return map $T_k$
and give the $k$-dependent coordinates $(Z, Y, W)$
that bring $T_k$ to the normal form.
In Section~\ref{s33-cone},
we prove the existence of the invariant cone fields
$\mathcal{C}^{ss}$ and $\mathcal{C}^{cu}$ on the domain of $T_k$.

In the following sections,
$(f, \Gamma)$ is assumed to satisfy \textbf{(QC)}
unless otherwise noted.

\subsection{Iterated local map and global map}\label{s31-FRM}
We begin by formally defining the first-return map.
We start by defining the iterated local map and the global map,
and subsequently define the first-return map as their composition.
Then we give the iterated local map formula and the global map formula.

\begin{notation}
  In this paper, we adopt the following convention for notation of derivatives.  
  \begin{itemize}
    \item For a real-valued function of several variables  
      $F(x_1, x_2, \cdots, x_n)$,  
      the partial derivative of $F$ with respect to  
      $\sigma \in \{x_1, x_2, \cdots, x_n\}$ is denoted by
      \begin{align}
        F_\sigma(x_1, x_2, \cdots, x_n).
        \label{e31-Fsig}
      \end{align}
      If the function has a subscript, such as $F_s$,  
      then its partial derivative with respect to $\sigma$  
      is denoted by $F_{s, \sigma}$.
    \item For a tuple of real-valued functions $(F_1, F_2, \cdots, F_m)$ of
      $(x_1, x_2, \cdots, x_n)$,  
      \begin{align*}
        \frac{\partial (F_1, F_2, \cdots, F_m)}{\partial (x_1, x_2, \cdots, x_n)}
        &= \begin{pmatrix}
            F_{1, x_1} & \cdots & F_{1, x_n} \\
            \vdots & \ddots & \vdots \\
            F_{m, x_1} & \cdots & F_{m, x_n}
          \end{pmatrix}
      \end{align*}
      denotes the $m \times n$ Jacobian matrix.
    \item The differential operator $\partial/(\partial \sigma)$  
      is written as $\partial_\sigma$.
  \end{itemize}
\end{notation}

\subsubsection{First-return map}
For the proper unfolding three-parameter family
$\{ f_\varepsilon \}_{\varepsilon \in R_\mathrm{prm}^*}$
defined in Section~\ref{s2-fam},
we define the local map $T_0$ in the same way as \eqref{eq-def_T00}:
\begin{align*}
  T_0 = T_0(\varepsilon;
    f, U_0, \{ f_\varepsilon \}_{\varepsilon \in R_\mathrm{prm}^*})
  := f_\varepsilon^{\mathrm{per}(O^*)}|_{U_0
    \cap f_\varepsilon^{-\mathrm{per}(O^*)}(U_0)},
\end{align*}
where $U_0$ and $\mathrm{per}(O^*)$ do not depend on $\varepsilon$.
Similar to equation \eqref{eq-02-def_T10} we also define the global map $T_1$
\begin{align}
  T_1 = T_1(\varepsilon; f, \Gamma, U_0, M^-_0, M^+_0,
  \{ f_\varepsilon \}_{\varepsilon \in R_\mathrm{prm}^*})
  := (f_\varepsilon^{\mathrm{per}(O^*)})^{n_0},
  \label{eq-03-defT1}
\end{align}
where $n_0 = n_0(M^-_0, M^+_0)$,
$M^-_0$,
and $M^+_0$ do not depend on $\varepsilon$.

Recall that the range over which the parameter $\varepsilon$ moves was
the three-dimensional open ball $R_\mathrm{prm}^*$ centered at $\varepsilon^*$.
Since \cite[Lemma~6]{GST2008},
by taking a smaller three-dimensional open ball
$R_\mathrm{prm} \subset R_\mathrm{prm}^*$ centered at $\varepsilon^*$
(the smaller one only depends on $(f, U_0)$),
there exist $\varepsilon$-dependent $C^r$ coordinates
$\bm{x} = (x_1, x_2, y)$ on $U_0$
such that the local map $T_0: (x_1, x_2, y) \mapsto (\hat{x}_1, \hat{x}_2, \hat{y})$
can be written in the following form by using these coordinates:
\begin{align}
\begin{aligned}
  \hat{x}_1
    &= \lambda(\varepsilon) x_1 \cos\omega
    - \lambda(\varepsilon) x_2 \sin\omega
    + p_1(x_1, x_2, y, \varepsilon), \\
  \hat{x}_2
    &= \lambda(\varepsilon) x_1 \sin\omega
    + \lambda(\varepsilon) x_2 \cos\omega
    + p_2(x_1, x_2, y, \varepsilon), \\
  \hat{y}
    &= \gamma(\varepsilon) y + p_3(x_1, x_2, y, \varepsilon),
\end{aligned}\label{eq-T0prm}
\end{align}
where $\lambda = \lambda(\varepsilon)$
and $\gamma = \gamma(\varepsilon)$ are
the continuations of $\lambda^*$ and $\gamma^*$,
respectively,
and they are $C^{r - 1}$ with respect to $\varepsilon$;
the coordinates $(x_1, x_2, y)$ are $C^{r - 2}$ with respect to parameters;
the $p_i$, $i \in \{1, 2, 3\}$,
are $C^{r - 2}$ with respect to $(x_1, x_2, y, \varepsilon)$;
\debugtext{
  more precisely:
  the $p_i$, $i \in \{1, 2, 3\}$, and their first and second partial derivatives
  with respect to $(x_1, x_2, y)$ are
  $C^{r - 2}$ with respect to $(x_1, x_2, y, \varepsilon)$;
}
the $p_i$ satisfies
\begin{align}
\begin{aligned}
  &p_i(0, 0, y, \varepsilon) \equiv 0, \quad
  p_i(x_1, x_2, 0, \varepsilon) \equiv 0, \quad
  p_{i, x_j}(0, 0, 0, \varepsilon) \equiv 0, \\
  &p_{i, y}(0, 0, 0, \varepsilon) \equiv 0, \quad
  p_{l, x_j}(0, 0, y, \varepsilon) \equiv 0, \quad
  p_{3, y}(x_1, x_2, 0, \varepsilon) \equiv 0
\end{aligned}\label{eq-loc_linear}
\end{align}
for any $i \in \{1, 2, 3\}$ and $j$, $l \in \{1, 2\}$.

\begin{remark} \label{r31-smth}
To describe the smoothness of the coordinates $(x_1, x_2, y)$ in more detail, 
if we take a $C^r$ coordinates $(s_1, s_2, t)$ on $U_0$
that do not depend on $\varepsilon$,
then for the $C^r$ coordinate transformations
$(x_1, x_2, y, \varepsilon) \mapsto (s_1, s_2, t)$,
it and its first and second partial derivatives with respect to $(x_1, x_2, y)$
are $C^{r - 2}$ with respect to $(x_1, x_2, y, \varepsilon)$.
Here, when $r = \infty$ or $\omega$,
we assume $r - k = r$ for any $k < \infty$.
see \cite[Remarks to Lemma~6]{GST2008} for more details.
\end{remark}

Since $(f, \Gamma)$ holds \textbf{(QC)},
by replacing $R_\mathrm{prm}$ with a smaller one
(the new smaller one only depends on
$(f, \Gamma, (U_0; x_1, x_2, y), M_0^-, M_0^+,
\{ f_\varepsilon \}_{\varepsilon \in R_\mathrm{prm}^*})$),
the implicit function theorem
extends $M^-_0$ and $M^+_0$ to depend on $\varepsilon$ as follows:
The $M^-(\varepsilon)$ and $M^+(\varepsilon)$ are $C^{r - 2}$
with respect to $\varepsilon$
\footnote{
  When applying the implicit function theorem,  
  since the equation involves the partial derivative
  of the global map in the $y$-direction,  
  the smoothness of the solution may appear to decrease by one.  
  However, partial differentiating the global map with respect to the spatial variables  
  $(x_1, x_2, y)$ does not affect the smoothness with respect to parameters,  
  which remains $C^{r - 2}$.  
  Therefore, the smoothness of $M^\pm(\varepsilon)$ is $C^{r - 2}$.    
}
such that $M^-(\varepsilon^*) = M^-_0$ and
$M^+(\varepsilon^*) = M^+_0$,
and they can be written 
\begin{align*}
  M^-(\varepsilon) = (0, 0, y^-(\varepsilon)), \quad
  M^+(\varepsilon) = (x_1^+(\varepsilon), x_2^+(\varepsilon), 0)
\end{align*}
by using the above coordinates $\bm{x}$ with
\begin{align*}
  \left.
    \bar{x}_1
  \right|_{\tilde{\bm{x}} = (0, 0, y^-(\varepsilon))}
    = x_1^+(\varepsilon), \quad
  \left.
    \bar{x}_2
  \right|_{\tilde{\bm{x}} = (0, 0, y^-(\varepsilon))}
    = x_2^+(\varepsilon), \quad
  \left.
    (\partial_{\tilde{y}} \bar{y})
  \right|_{\tilde{\bm{x}} = (0, 0, y^-(\varepsilon))} = 0,
\end{align*}
where we write
$T_1: \tilde{\bm{x}} = (\tilde{x}_1, \tilde{x}_2, \tilde{y})
\mapsto (\bar{x}_1, \bar{x}_2, \bar{y})$
by using the coordinates $\bm{x}$.
In the following,
sometimes $(\varepsilon)$ may be dropped:
\begin{align*}
  M^-(\varepsilon) = M^-, \quad
  M^+(\varepsilon) = M^+, \quad
  y^-(\varepsilon) = y^-, \quad
  x^+_1(\varepsilon) = x^+_1, \quad
  x^+_2(\varepsilon) = x^+_2.
\end{align*}

We may assume that $R_\mathrm{prm} \subset R_\mathrm{prm}^*$ is given by
\begin{align*}
  R_\mathrm{prm}
    = I_\mathrm{prm}
    \times (\omega^* + I_\mathrm{prm})
    \times I_\mathrm{prm}, \quad
    I_\mathrm{prm} := (-\delta_\mathrm{prm}, \delta_\mathrm{prm}),\:
    \delta_\mathrm{prm} > 0,
\end{align*}
in the $(\mu, \omega, \rho)$-space,
and we sometimes write $R_\mathrm{prm}$ as $R_\mathrm{prm}(\delta_\mathrm{prm})$.
We define the pair of $U_0$ and $(x_1, x_2, y)$ as
\begin{align*}
  \mathbb{U}_0 := (U_0; x_1, x_2, y)
\end{align*}
and the tuple of the core objects as
\begin{align}
  \mathbb{F} := (f, \Gamma, \mathbb{U}_0, M_0^-, M_0^+,
  \{ f_\varepsilon \}_{\varepsilon \in R_\mathrm{prm}^*})
  \label{eq-03-mathbbF}
\end{align}
to simplify the notation.
There exist small numbers
\begin{align*}
  \hat{\delta}_\mathrm{dom}
    = \hat{\delta}_\mathrm{dom}(\mathbb{F}) > 0, \quad
  \delta_\mathrm{prm}^{new}
    = \delta_\mathrm{prm}^{new}(\mathbb{F})
    \in (0, \delta_\mathrm{prm})
\end{align*}
such that the two cubes
\begin{align*}
  \Pi^- &= \Pi^-(\varepsilon, \delta_\mathrm{dom})
    := [-\delta_\mathrm{dom}, \delta_\mathrm{dom}]^3
    + M^-(\varepsilon), \\
  \Pi^+ &= \Pi^+(\varepsilon, \delta_\mathrm{dom})
    := [-\delta_\mathrm{dom}, \delta_\mathrm{dom}]^3
    + M^+(\varepsilon)
\end{align*}
are disjoint and
\begin{align*}
  \Pi^-(\varepsilon, \delta_\mathrm{dom}), \,
  \Pi^+(\varepsilon, \delta_\mathrm{dom})
    \subset U_0, \quad
  T_1(\Pi^-(\varepsilon, \delta_\mathrm{dom})) \subset U_0
\end{align*}
for any $\varepsilon \in R_\mathrm{prm}(\delta_\mathrm{prm}^{new})$ and
$\delta_\mathrm{dom} \in (0, \hat{\delta}_\mathrm{dom})$.
In the following,
we drop the `new'.
By replacing $\delta_\mathrm{prm} > 0$ with a smaller one
(the new smaller one only depends on $\mathbb{F}$),
there exists
\begin{align*}
  \kappa(\delta_\mathrm{dom})
  = \kappa(\delta_\mathrm{dom}; \mathbb{F}) > 0
\end{align*}
such that
\begin{align*}
  \Pi_k = \Pi_k(\varepsilon, \delta_\mathrm{dom})
    := \Pi^+(\varepsilon, \delta_\mathrm{dom})
    \cap T_0^{-k}(\Pi^-(\varepsilon, \delta_\mathrm{dom}))
\end{align*}
is a nonempty strip-like region for any $k > \kappa(\delta_\mathrm{dom})$
and $\varepsilon \in R_\mathrm{prm}$.
The iterated local map $T_0^k$ is defined on $\Pi_k$.

We define the \textit{first-return map}
$T_k: \Pi_k \to U_0$ by
\begin{align}
  T_k = T_k(\varepsilon, \delta_\mathrm{dom}; \mathbb{F})
    := T_1 \circ T_0^k
  \label{eq-FRmap_def}
\end{align}
for any $\varepsilon \in R_\mathrm{prm}$,
$\delta_\mathrm{dom} \in (0, \hat{\delta}_\mathrm{dom})$,
and $k > \kappa(\delta_\mathrm{dom})$.

\begin{remark} \label{rem-smaller}
  Throughout the rest of the paper,
  we always consider the
  $\hat{\delta}_\mathrm{dom}$ and $\delta_\mathrm{prm}$ are fixed,
  $\delta_\mathrm{dom}$ runs in $(0, \hat{\delta}_\mathrm{dom})$,
  and $\kappa(\delta_\mathrm{dom})$ is a function of $\delta_\mathrm{dom}$.
  By contrast,
  we sometimes replace $\hat{\delta}_\mathrm{dom}$, $\delta_\mathrm{prm}$,
  $\kappa(\delta_\mathrm{dom})$ with smaller ones
  $\hat{\delta}_\mathrm{dom}^{new} \in (0, \hat{\delta}_\mathrm{dom})$,
  $\delta_\mathrm{prm}^{new} \in (0, \delta_\mathrm{prm})$,
  and a bigger one
  $\kappa^{new}(\delta_\mathrm{dom}) > \kappa(\delta_\mathrm{dom})$, respectively.
  However, throughout the rest of the paper,
  we always take them as
  \begin{gather*}
    \hat{\delta}_\mathrm{dom}^{new}
      = \hat{\delta}_\mathrm{dom}^{new}(\mathbb{F}), \quad
    \delta_\mathrm{prm}^{new}
      = \delta_\mathrm{prm}^{new}(\mathbb{F}), \quad
    \kappa^{new}(\delta_\mathrm{dom})
      = \kappa^{new}(\delta_\mathrm{dom}; \mathbb{F}),
  \end{gather*}
  in other words, new ones
  $\hat{\delta}_\mathrm{dom}^{new}$, $\delta_\mathrm{prm}^{new}$,
  $\kappa^{new}(\delta_\mathrm{dom})$ at least just depend on $\mathbb{F}$.
\end{remark}

\subsubsection{Representation of the local map and the global map} \label{s312-r}
From \cite[Lemma~7]{GST2008},
replacing $\delta_\mathrm{prm} > 0$
and $\kappa(\delta_\mathrm{dom})$ with smaller and larger ones
according to Remark~\ref{rem-smaller} yields the following:
if $T_0^k: \Pi_k \ni (x_1, x_2, y)
\mapsto (\tilde{x}_1, \tilde{x}_2, \tilde{y}) \in \Pi^-$,
then
\begin{align}
\begin{aligned}
  \tilde{x}_1
    &= (\lambda(\varepsilon))^k x_1 \cos(k \omega)
    - (\lambda(\varepsilon))^k x_2 \sin(k \omega)
    + \hat{\lambda}^k q_k^{(1)}(x_1, x_2, \tilde{y}, \varepsilon), \\
  \tilde{x}_2
    &= (\lambda(\varepsilon))^k x_1 \sin(k \omega)
    + (\lambda(\varepsilon))^k x_2 \cos(k \omega)
    + \hat{\lambda}^k q_k^{(2)}(x_1, x_2, \tilde{y}, \varepsilon), \\
  y &= (\gamma(\varepsilon))^{-k} \tilde{y}
    + \hat{\gamma}^{-k} q_k^{(3)}(x_1, x_2, \tilde{y}, \varepsilon)
\end{aligned}\label{eq-03-T0kprm}
\end{align}
for any
$k > \kappa(\delta_\mathrm{dom})$ and
$\varepsilon \in R_\mathrm{prm}$,
where $\hat{\lambda} = \hat{\lambda}(\mathbb{F})$
and $\hat{\gamma} = \hat{\gamma}(\mathbb{F})$ are constants such that
$\hat{\lambda} < \lambda(\varepsilon)$ and $\hat{\gamma} > \gamma(\varepsilon)$
for any $\varepsilon \in R_\mathrm{prm}$;
the $q_k^{(i)}$, $i \in \{1, 2, 3\}$,
are $C^{r - 2}$ with respect to $(x_1, x_2, \tilde{y}, \varepsilon)$;
\debugtext{
  more precisely:
  the $q_k^{(i)}$ ($i \in \{1, 2, 3\}$) and their first and second partial derivatives
  with respect to $(x_1, x_2, \tilde{y})$ are
  $C^{r - 2}$ with respect to $(x_1, x_2, \tilde{y}, \varepsilon)$;
}
the $j$-th partial derivatives of $q_k^{(i)}$
with respect to $(x_1, x_2, \tilde{y}, \varepsilon)$
are bounded with respect to $(k, x_1, x_2, \tilde{y}, \varepsilon)$
for any $j \in \{0, 1, \cdots, r-2\}$.
Note that $\hat{\lambda}$, $\hat{\gamma}$
can be taken so that
\begin{align*}
  &(\lambda^*)^2
  = \lambda^* |\gamma^*|^{-1}
  = |\gamma^*|^{-2}
  < \hat{\gamma}^{-1}
  < \hat{\lambda}
  < \lambda^*
  = |\gamma^*|^{-1}.
\end{align*}
Moving $\hat{\gamma}^{-1}$ closer to $(\lambda^*)^2$
from the right side
and $\hat{\lambda}$ closer to $\lambda^*$
from the left side,
\begin{align*}
  \lambda^* < |\gamma^*| \hat{\gamma}^{-1}
  < \hat{\lambda}^{1/2}
\end{align*}
can be further satisfied.
We also take a constant $\hat{\hat{\lambda}} = \hat{\hat{\lambda}}(\mathbb{F})$
with 
\begin{align}
  \hat{\lambda} < \hat{\hat{\lambda}} < \lambda(\varepsilon).
  \label{e31-hhlmd}
\end{align}
Thus,
replacing $\delta_\mathrm{prm} > 0$ with a smaller one
according to Remark~\ref{rem-smaller} if necessary,
we may suppose
\begin{align}
\begin{gathered}
  (\lambda(\varepsilon))^2, \,
  \lambda(\varepsilon) |\gamma(\varepsilon)|^{-1}, \,
  |\gamma(\varepsilon)|^{-2}
  < \hat{\gamma}^{-1}
  < \hat{\lambda}
  < \lambda(\varepsilon), \,
  |\gamma(\varepsilon)|^{-1}, \\
  \lambda(\varepsilon) |\gamma(\varepsilon)| \hat{\gamma}^{-1}
  < \hat{\lambda}, \quad
  \lambda(\varepsilon)
  < |\gamma(\varepsilon)| \hat{\gamma}^{-1}
  < \hat{\lambda}^{1/2}
\end{gathered} \label{eq-03-hatC}
\end{align}
for any $\varepsilon \in R_\mathrm{prm}$.

\medskip

As in equation \eqref{eq-02-ExpT10},
the global map
$T_1: \Pi^- \ni (\tilde{x}_1, \tilde{x}_2, \tilde{y})
\mapsto (\bar{x}_1, \bar{x}_2, \bar{y}) \in U_0$
is written as follows.
\begin{align*}
  \bar{x}_1 - x_1^+
    &= a_{11}' \tilde{x}_1 + a_{12}' \tilde{x}_2
    + b_1' (\tilde{y} - y^-)
    + O(\|(\tilde{x}_1, \tilde{x}_2, \tilde{y} - y^-)\|^2), \\
  \bar{x}_2 - x_2^+
    &= a_{21}' \tilde{x}_1 + a_{22}' \tilde{x}_2
    + b_2' (\tilde{y} - y^-)
    + O(\|(\tilde{x}_1, \tilde{x}_2, \tilde{y} - y^-)\|^2), \\
  \bar{y}
    &= y^+(\varepsilon) + c_1' \tilde{x}_1 + c_2' \tilde{x}_2
    + d' (\tilde{y} - y^-)^2
    + O(\|(\tilde{x}_1, \tilde{x}_2, \tilde{y} - y^-)\|^2),
\end{align*}
where
$O(\cdot)$ are terms of second order or higher of the Taylor expansion,
excluding the explicitly stated terms,
$\|\cdot\|$ denotes the Euclidean norm, and
\begin{align*}
  y^+ = y^+(\varepsilon) :=
    \left.
      \bar{y}
    \right|_{\tilde{\bm{x}} = (0, 0, y^-(\varepsilon))}.
\end{align*}
Moreover,
since the smoothness of the coordinates $(x_1, x_2, y)$ with respect to parameters  
is $C^{r - 2}$,
$O(\cdot)$ are $C^{r - 2}$
with respect to $(\tilde{x}_1, \tilde{x}_2, \tilde{y}, \varepsilon)$
and $y^+(\varepsilon)$ is $C^{r - 2}$ with respect to $\varepsilon$.
\debugtext{
  More precisely:
  $O(\cdot)$ and their first and second partial derivatives
  with respect to $(\tilde{x}_1, \tilde{x}_2, \tilde{y})$ are
  $C^{r - 2}$ with respect to $(\tilde{x}_1, \tilde{x}_2, \tilde{y}, \varepsilon)$.
}
By the definition of the proper unfolding,
\begin{align}
    y^+_\mu(0, \omega^*, 0) \neq 0
  \label{eq-mu_yplus}
\end{align}
see Section~\ref{s21-PU}.
Note that the above coefficients depend on
$\mathbb{F}$ and $\varepsilon$ since \eqref{eq-03-defT1} and
they are $C^{r - 2}$ with respect to $\varepsilon$:
\begin{align}
\begin{aligned}
  a_{ij}' = a_{ij}'(\varepsilon)
    = a_{ij}'(\varepsilon; \mathbb{F}), \quad
  b_i' = b_i'(\varepsilon)
    = b_i'(\varepsilon; \mathbb{F}), \quad
  c_i' = c_i'(\varepsilon)
    = c_i'(\varepsilon; \mathbb{F}), \quad
  d' = d'(\varepsilon)
    = d'(\varepsilon; \mathbb{F})
\end{aligned} \label{eq-03-coepdep}
\end{align}
for any $i$, $j \in \{1, 2\}$.

\begin{remark} \label{rmk-prm}
  Up to this point, we have regarded $\varepsilon = (\mu, \omega, \rho)$  
  as the parameter (see Remark~\ref{rmk-param}).
  However, from this point on, we switch the roles of $\mu$ and $y^+$,  
  and treat $(y^+, \omega, \rho)$ as the parameters.  
  This is justified by \eqref{eq-mu_yplus} and the fact that $y^+(\varepsilon)$  
  is $C^{r - 2}$.  
  To simplify notation, we will write $y^+$ again as $\mu$.  
  Therefore, we continue to write $\varepsilon = (\mu, \omega, \rho)$,  
  but note that from now on, $\mu$ refers to $y^+$.
\end{remark}

Consider the new $\varepsilon$-dependent $C^r$ coordinates
\begin{align*}
  (x_1^{new}, x_2^{new})^\mathsf{T}
    := R(-\arctan_2(\bm{b}')) (x_1, x_2)^\mathsf{T}, \quad
  y^{new} := y,
\end{align*}
where $\bm{b}' = (b_1', b_2')$;
the $\arctan_2(v_1, v_2) \in [0, 2\pi)$ is the angle determined by
\begin{align}
  \cos(\arctan_2((v_1, v_2))) = \frac{v_1}{\sqrt{v_1^2 + v_2^2}}, \quad
  \sin(\arctan_2((v_1, v_2))) = \frac{v_2}{\sqrt{v_1^2 + v_2^2}}
  \label{e31-at2}
\end{align}
for any $(v_1, v_2) \in \mathbb{R}^2 \setminus \{0\}$;
$X^\mathsf{T}$ denotes the transpose of a matrix $X$;
and $R(\varphi)$ denotes the rotation matrix of angle $\varphi$:
\begin{align}
  R(\varphi) =
  \begin{pmatrix}
    \cos\varphi & -\sin\varphi \\
    \sin\varphi & \cos\varphi \\
  \end{pmatrix}.
  \label{eq-03-rotMat}
\end{align}
In particular,
\begin{align*}
  R(-\arctan_2(\bm{b}')) = \frac{1}{\sqrt{(b_1')^2 + (b_2')^2}}
  \begin{pmatrix}
    b_1' & b_2' \\
    -b_2' & b_1' \\
  \end{pmatrix}.
\end{align*}
The smoothness of the new coordinates
$(x_1^{new}, x_2^{new}, y^{new})$  
with respect to the parameter $\varepsilon$ is the same as that of  
the previous coordinates $(x_1, x_2, y)$.  

The coordinates $(x_1^{new}, x_2^{new}, y^{new})$ do not break
equations \eqref{eq-T0prm} and \eqref{eq-loc_linear}, and further rewrite
the global map $T_1$ as follows:
\begin{align}
\begin{aligned}
  \bar{x}_1 - x_1^+
    &= a_{11} \tilde{x}_1 + a_{12} \tilde{x}_2
    + b (\tilde{y} - y^-)
    + O(\|(\tilde{x}_1, \tilde{x}_2, \tilde{y} - y^-)\|^2), \\
  \bar{x}_2 - x_2^+
    &= a_{21} \tilde{x}_1 + a_{22} \tilde{x}_2
    + O(\|(\tilde{x}_1, \tilde{x}_2, \tilde{y} - y^-)\|^2), \\
  \bar{y}
    &= \mu + c_1 \tilde{x}_1 + c_2 \tilde{x}_2
    + d (\tilde{y} - y^-)^2
    + O(\|(\tilde{x}_1, \tilde{x}_2, \tilde{y} - y^-)\|^2),
\end{aligned} \label{eq-03-T1prm}
\end{align}
where $O(\cdot)$ are terms of second order or higher of the Taylor expansion,
excluding the explicitly stated terms,
the `new' is dropped,
and the above coefficients can be written as
\begin{gather*}
  \bm{A} = R(-\arctan_2(\bm{b}'))\bm{A}'
    R(\arctan_2(\bm{b}')), \quad
  (b, 0)^\mathsf{T} = R(-\arctan_2(\bm{b}'))(\bm{b}')^\mathsf{T}, \\
  \bm{c}^\mathsf{T} = R(-\arctan_2(\bm{b}')) (\bm{c}')^\mathsf{T}, \quad
  d = d'.
\end{gather*}
Here, we put
\begin{align*}
  \bm{A}' := 
  \begin{pmatrix}
    a_{11}' & a_{12}' \\
    a_{21}' & a_{22}'
  \end{pmatrix}, \quad
  \bm{A} := 
  \begin{pmatrix}
    a_{11} & a_{12} \\
    a_{21} & a_{22}
  \end{pmatrix}, \quad
  \bm{c}' := (c_1', c_2'), \quad
  \bm{c} := (c_1, c_2).
\end{align*}
Note that the coefficient of $\tilde{y}$
in the $(\bar{x}_2 - x_2^+)$ equation is zero and
\begin{align}
  b = \sqrt{(b_1')^2 + (b_2')^2} \geq C \label{eq-03-b}
\end{align}
for some constant $C = C(\mathbb{F}) > 0$
by replacing $\delta_\mathrm{prm} > 0$ with a smaller one
according to Remark~\ref{rem-smaller}
since the original $f$ is a diffeomorphism.
We sometimes denote these new coefficients as in \eqref{eq-03-coepdep}:
\begin{align}
  \begin{aligned}
    a_{ij} = a_{ij}(\varepsilon)
      = a_{ij}(\varepsilon; \mathbb{F}), \quad
    b = b(\varepsilon)
      = b(\varepsilon; \mathbb{F}), \quad
    c_i = c_i(\varepsilon)
      = c_i(\varepsilon; \mathbb{F}), \quad
    d = d(\varepsilon)
      = d(\varepsilon; \mathbb{F})
  \end{aligned} \label{eq-03-coedep}
\end{align}
for any $i$, $j \in \{1, 2\}$.
Note also that the above functions $a_{ij}$, $b$, $c_i$, and $d$ are $C^{r - 2}$
with respect to $\varepsilon$.

\subsection{Normal form for the first-return map}
\label{s32-NF}
Recall the first-return map $T_k$ in \eqref{eq-FRmap_def}.
The same method as in \cite{LLST2022} is used to put $T_k$ into normal form.
It is done in several steps.
Before stating results concerning the normal form,  
we make some preparations.  

\medskip

Let us recall the coefficients in \eqref{eq-03-coedep}.  
We define  
\begin{align}  
\begin{aligned}  
  &\alpha^* = \alpha^*_k(\varepsilon)  
    = \alpha^*_k(\varepsilon; \mathbb{F})  
    := c_1(\varepsilon)\cos(k \omega)
    + c_2(\varepsilon)\sin(k \omega).
\end{aligned} \label{e32-alp}  
\end{align}  
This quantity can be written as  
\begin{align}  
  \alpha^*_k(\varepsilon) = c(\varepsilon) \sin (k \omega + \eta^*(\varepsilon))  
  \label{e52-alp2}
\end{align}  
where  
\begin{align}  
  c = c(\varepsilon) = c(\varepsilon; \mathbb{F})  
    := \sqrt{(c_1(\varepsilon))^2 + (c_2(\varepsilon))^2}, \quad  
  \eta^* = \eta^*(\varepsilon)
  := \arctan_2 (c_1(\varepsilon), c_2(\varepsilon)) \in [0, 2\pi),
  \label{eq-03-defc}
\end{align}  
and $\arctan_2$ is a function defined by \eqref{e31-at2}.
We define
\begin{align}
  &I_k^\mathrm{bd}
    := \{ \omega \in (\omega^* + I_\mathrm{prm}) \:|\: 
    |\sin(k \omega + \eta^*(0, \omega^*, 0))| > 2e_\mathrm{bd} \},
  \label{e32-Ikbd} \\
  &R_k^\mathrm{bd} := I_\mathrm{prm} \times I_k^\mathrm{bd} \times I_\mathrm{prm}
    \quad (\subset R_\mathrm{prm}),
  \label{e32-Rkbd}
\end{align}
where $e_\mathrm{bd} \in (0, 1)$ is a completely arbitrary number.  
Since $e_\mathrm{bd}$ can be chosen freely,
hereafter, we always set $e_\mathrm{bd} = 1/20$.  
By replacing $\delta_\mathrm{prm} > 0$ with a smaller one
according to Remark~\ref{rem-smaller}, we have
\begin{align}
  |\sin(k \omega + \eta^*(\mu, \omega, \rho))| > e_\mathrm{bd}
  \label{e32-sinE}
\end{align}
for any $\varepsilon \in R_k^\mathrm{bd}$.
Since $f$ is a diffeomorphism,
\eqref{eq-03-T1prm} implies $c(\varepsilon^*) \neq 0$,
and $c(\varepsilon) \geq C$ for some constant $C = C(\mathbb{F}) > 0$
by replacing $\delta_\mathrm{prm} > 0$ with a smaller one.
Thus, as long as $\varepsilon \in R_k^\mathrm{bd}$,  
we have $|\alpha^*_k(\varepsilon)| \geq C$ for some constant $C = C(\mathbb{F}) > 0$.
We further define
\begin{align}
  E_k = E_k(\varepsilon)
    = E_k(\varepsilon; \mathbb{F})
    := -b(\varepsilon) \alpha^*_k(\varepsilon).
  \label{eq-03-Ek}
\end{align}
By \eqref{eq-03-b},
this quantity is also bounded away from zero
when $\varepsilon \in R_k^\mathrm{bd}$.

\begin{notation}
  Throughout the paper, unless otherwise noted,
  for any $F = F(\varepsilon, \delta_\mathrm{dom}, k, M)$
  and $G = G(\varepsilon, \delta_\mathrm{dom}, k, M)$,
  $F = O(G)$ means there exists $C = C(\mathbb{F}) > 0$ such that
  \begin{align*}
    |F| \leq C |G| 
  \end{align*}
  for any $\varepsilon \in R_\mathrm{prm}$,
  $\delta_\mathrm{dom} \in (0, \hat{\delta}_\mathrm{dom})$,
  $k > \kappa(\delta_\mathrm{dom})$,
  and $M \in \Pi_k(\varepsilon, \delta_\mathrm{dom})$.
\end{notation}

\begin{proposition}[Normal form of $T_k$]\label{p32-NF}
  For any
  $k \in \mathbb{Z}_{> \kappa(\delta_\mathrm{dom})}$
  and $\varepsilon \in R_k^\mathrm{bd}$,
  there exist $(\varepsilon, k)$-dependent $C^r$ coordinates $(Z, Y, W)$ on $\Pi_k$
  and $\delta_\mathrm{dom}' > 0$ such that $\Pi_k$ contains
  \begin{align}
    \Pi_k' = \Pi_k'(\varepsilon, \delta_\mathrm{dom}')
      := [-\delta_\mathrm{dom}', \delta_\mathrm{dom}']^3
    \label{eq-03-Pikcirc}
  \end{align}
  in the $(Z, Y, W)$ coordinates
  and the first-return map
  $T_k: \Pi_k' \ni (Z, Y, W) \mapsto (\bar{Z}, \bar{Y}, \bar{W})$
  with $(\bar{Z}, \bar{Y}, \bar{W}) \in \Pi_k$
  can be written by the form
  \begin{align}
  \begin{aligned}
    \bar{Z}
      &= \lambda^k \alpha_1 Z
      - E_k Y
      + \lambda^k \beta_1 W
      + h_1(Z, Y, W, \varepsilon), \\
    \bar{Y}
      &= \hat{\mu}
      + \lambda^k \gamma^k Z
      + d \gamma^k Y^2
      + h_2(Z, Y, W, \varepsilon), \\
    \bar{W}
      &= \lambda^k \alpha_3 Z
      + \lambda^k \beta_3 W
      + h_3(Z, Y, W, \varepsilon),
  \end{aligned} \label{eq-03-NF}
  \end{align}
  where the above quantities are given as follows:
  the smoothness of the coordinates $(Z, Y, W)$
  with respect to parameter $\varepsilon$ is the same as that of $(x_1, x_2, y)$,
  see also Remark~\ref{r31-smth};
  the $\delta_\mathrm{dom}'$ holds
  \begin{align}
    C_1 \delta_\mathrm{dom} \leq \delta_\mathrm{dom}' \leq C_2 \delta_\mathrm{dom}
    \label{e32-ddomp}
  \end{align}
  for some constants $C_1 = C_1(\mathbb{F}) > 0$
  and $C_2 = C_2(\mathbb{F}) > 0$;
  the $\hat{\mu}$ is given as
  \begin{align}
  \begin{aligned}
    \hat{\mu} = 
      \gamma^k \mu - y^-
      + \lambda^k \gamma^k (\alpha^* x_1^+ + \beta^* x_2^+)
      + \gamma^k O(\hat{\lambda}^k)
  \end{aligned} \label{eq-03-muhat}
  \end{align}
  and the above $O(\hat{\lambda}^k)$ is a $C^{r - 2}$ function
  of $\varepsilon$ and its first partial derivatives
  with respect to $\mu$ and $(\omega, \rho)$ have estimates of
  $O(\gamma^k \hat{\gamma}^{-k})$ and $O(\hat{\hat{\lambda}}^k)$, respectively,
  where $\hat{\hat{\lambda}}$ is a constant with \eqref{e31-hhlmd};
  the coefficients
  \begin{align*}
    \alpha_1 = \alpha_1^{(k)}(\varepsilon), \quad
    \beta_1 = \beta_1^{(k)}(\varepsilon), \quad
    \alpha_3 = \alpha_3^{(k)}(\varepsilon), \quad
    \beta_3 = \beta_3^{(k)}(\varepsilon)
  \end{align*}
  are $C^{r - 2}$ with respect to $\varepsilon$ and satisfy
  \begin{align}
    \alpha_i^{(k)}, \,
    \beta_i^{(k)}, \,
    \alpha_{i, \tilde{\sigma}}^{(k)}, \,
    \beta_{i, \tilde{\sigma}}^{(k)}, \,
      = O(1), \quad
    \alpha_{i, \omega}^{(k)}, \,
    \beta_{i, \omega}^{(k)}
      = O(k) \label{eq-03-ab_ev}
  \end{align}
  for any $i \in \{1, 3\}$
  and $\tilde{\sigma} \in \{\mu, \rho\}$.
  the higher order terms are given as
  \begin{align}
  \begin{aligned}
    &h_i = O(\hat{\lambda}^k)(Z + W) + O(\lambda^k)Y + O(Y^2), \\
    &h_2 = \gamma^k \left(
      O(\hat{\lambda}^k)(Z + W) + O(\lambda^k)(Z + W)Y + O(\lambda^k) Y^2 + O(Y^3)
      \right).
  \end{aligned} \label{eq-03-NF_higher1}
  \end{align}
  Furthermore, for the partial derivatives, we have
  \begin{align}
  \begin{aligned}
    &h_{i, X}, \,
    h_{i, X X'}, \,
    h_{i, X X' X''}, \,
    h_{i, X X' Y}, \,
    h_{i, X Y Y}
      = O(\hat{\lambda}^k), \quad
    h_{i, Y}, \,
    h_{i, X Y}
      = O(\lambda^k), \\
    &h_{i, Y Y}, \,
    h_{i, Y Y Y}
      = O(1), \quad
    h_{i, \sigma}, \,
    h_{i, X \sigma}
      = O(\hat{\lambda}^k), \\
    &h_{i, Y \mu}
      = O(\gamma^k \hat{\gamma}^{-k})
      + O(\lambda^k \gamma^{2k} \hat{\gamma}^{-k}), \quad
    h_{i, Y \sigma'}
      = O(k \lambda^k),
  \end{aligned} \label{eq-03-NF_higher2}
  \end{align}
  \begin{align}
  \begin{aligned}
    &h_{2, X}, \,
    h_{2, X X'}, \,
    h_{2, X X' X''}, \,
    h_{2, X X' Y}
      = \gamma^k O(\hat{\lambda}^k), \quad
    h_{2, Y}, \,
    h_{2, X Y}
      = \gamma^k O(\lambda^k), \\
    &h_{2, Y Y}
      = \gamma^k O(\gamma^k \hat{\gamma}^{-k}), \quad 
    h_{2, X Y Y}
      = \gamma^k (O(\gamma^k \hat{\gamma}^{-k})
      + O(\lambda^k \gamma^{2k} \hat{\gamma}^{-k})), \quad
    h_{2, Y Y Y}
      = \gamma^k O(1), \\
    &h_{2, \mu}
      = \gamma^k O(\gamma^k \hat{\gamma}^{-k}), \quad
    h_{2, \sigma'}, \,
    h_{2, X \sigma'}
      = \gamma^k O(\hat{\hat{\lambda}}^k), \\
    &h_{2, X \mu}, \,
    h_{2, Y \mu}
      = \gamma^k (O(\gamma^k \hat{\gamma}^{-k})
      + O(\lambda^k \gamma^{2k} \hat{\gamma}^{-k})), \quad
    h_{2, Y \sigma'}
      = \gamma^k O(k \lambda^k)
  \end{aligned} \label{eq-03-NF_higher3}
  \end{align}
  for any $Y$ with $|Y| \leq O(\lambda^k)$,
  $\mu$ with $|\mu| \leq O(\lambda^k) + O(\gamma^{-k})$,
  $i \in \{1, 3\}$,
  $X$, $X'$, $X'' \in \{Z, W\}$,
  $\sigma \in \{\mu, \omega, \rho\}$,
  and $\sigma' \in \{\omega, \rho\}$.
\end{proposition}

\begin{remark} \label{rem-proof}
  The essential part of the proof follows \cite{LLST2022}.  
  In Section~\ref{s4-Hopf}, it becomes necessary to estimate  
  the higher-order partial derivatives of the remainder terms $h_i$  
  ($i \in \{1, 2, 3\}$).  
  As a new element, we have incorporated these estimates into the proof.    
\end{remark}

\begin{proof}[Proof of Proposition~\ref{p32-NF}]
We divide the proof into several steps.

\medskip

\textbf{(1) Composition.}\quad
Substituting \eqref{eq-03-T0kprm} into \eqref{eq-03-T1prm},
we get
\begin{align}
\begin{aligned}
  \bar{x}_1 - x_1^+ &= \lambda^k \hat{\alpha}_1 x_1
    + \lambda^k \hat{\beta}_1 x_2
    + b (\tilde{y} - y^-)
    + O(\lambda^k) x_1 (\tilde{y} - y^-)
    + O(\lambda^k) x_2 (\tilde{y} - y^-) \\
    &\quad + O((\tilde{y} - y^-)^2)
    + O(\hat{\lambda}^k), \\
  \bar{x}_2 - x_2^+ &= \lambda^k \hat{\alpha}_2 x_1
    + \lambda^k \hat{\beta}_2 x_2
    + O(\lambda^k) x_1 (\tilde{y} - y^-)
    + O(\lambda^k) x_2 (\tilde{y} - y^-) \\
    &\quad + O((\tilde{y} - y^-)^2)
    + O(\hat{\lambda}^k), \\
  \bar{y} &= \mu
    + \lambda^k \alpha^* x_1
    + \lambda^k \beta^* x_2
    + O(\lambda^k) x_1 (\tilde{y} - y^-)
    + O(\lambda^k) x_2 (\tilde{y} - y^-) \\
    &\quad + d(\tilde{y} - y^-)^2
    + O((\tilde{y} - y^-)^3)
    + O(\hat{\lambda}^k),
\end{aligned} \label{eq-03-comp}
\end{align}
where
\begin{align}
  \begin{aligned}
    &\beta^* = \beta^*_k(\varepsilon)
    = \beta^*_k(\varepsilon; \mathbb{F})
      := -c_1(\varepsilon)\sin(k \omega)
        + c_2(\varepsilon)\cos(k \omega), \\
    &\hat{\alpha}_i = \hat{\alpha}_i^{(k)}(\varepsilon)
      = \hat{\alpha}_i^{(k)}(\varepsilon; \mathbb{F})
      := a_{i1}(\varepsilon)\cos(k \omega)
        + a_{i2}(\varepsilon)\sin(k \omega), \\
    &\hat{\beta}_i = \hat{\beta}_i^{(k)}(\varepsilon)
    = \hat{\beta}_i^{(k)}(\varepsilon; \mathbb{F})
      := -a_{i1}(\varepsilon)\sin(k \omega)
        + a_{i2}(\varepsilon)\cos(k \omega)
\end{aligned} \label{e32-bsetc}
\end{align}
for any $i \in \{1, 2\}$;
the $O(\lambda^k)$ are at least
$C^{r - 2}$ functions of $\varepsilon$;
the $O(\hat{\lambda}^k)$ are
$C^{r - 2}$ functions of $(x_1, x_2, \tilde{y}, \varepsilon)$;
the $O((\tilde{y} - y^-)^i)$ ($i \in \{2, 3\}$) are
$C^{r - 2}$ functions of $(\tilde{y}, \varepsilon)$.
Here, we used \eqref{eq-03-hatC} to sort $O(\cdot)$ terms.
\debugtext{
  Hint:
  \begin{align*}
    \bar{x}_1 - x_1^+ 
      &= a_{11} \left(
        \lambda^k x_1 \cos(k \omega)
        - \lambda^k x_2 \sin(k \omega)
        + \hat{\lambda}^k q_k^{(1)}
      \right)
      + a_{12} \left(
        \lambda^k x_1 \sin(k \omega)
        + \lambda^k x_2 \cos(k \omega)
        + \hat{\lambda}^k q_k^{(2)}
      \right) \\
      &\quad + b (\tilde{y} - y^-)
      + O(\lambda^{2k}) + O((\tilde{y} - y^-)^2) \\
      &\quad + O(1) \left(
        \lambda^k x_1 \cos(k \omega)
        - \lambda^k x_2 \sin(k \omega)
        + \hat{\lambda}^k q_k^{(1)}
      \right) (\tilde{y} - y^-) \\
      &\quad + O(1) \left(
        \lambda^k x_1 \sin(k \omega)
        + \lambda^k x_2 \cos(k \omega)
        + \hat{\lambda}^k q_k^{(2)}
      \right) (\tilde{y} - y^-) \\
      &= \lambda^k \hat{\alpha}_1 x_1 + \lambda^k \hat{\beta}_1 x_2
        + b (\tilde{y} - y^-)
        + O(\lambda^k) x_1 (\tilde{y} - y^-)
        + O(\lambda^k) x_2 (\tilde{y} - y^-)
        + O((\tilde{y} - y^-)^2) \\
      &\quad + a_{11} \hat{\lambda}^k q_k^{(1)}
      + a_{12} \hat{\lambda}^k q_k^{(2)}
      + O(\lambda^{2k})
      + O(1) \hat{\lambda}^k q_k^{(1)} (\tilde{y} - y^-)
      + O(1) \hat{\lambda}^k q_k^{(2)} (\tilde{y} - y^-), \\
    \bar{x}_2 - x_2^+
      &= a_{21} \left(
        \lambda^k x_1 \cos(k \omega)
        - \lambda^k x_2 \sin(k \omega)
        + \hat{\lambda}^k q_k^{(1)}
      \right)
      + a_{22} \left(
        \lambda^k x_1 \sin(k \omega)
        + \lambda^k x_2 \cos(k \omega)
        + \hat{\lambda}^k q_k^{(2)}
      \right) \\
      &\quad + O(\lambda^{2k}) + O((\tilde{y} - y^-)^2) \\
      &\quad + O(1) \left(
        \lambda^k x_1 \cos(k \omega)
        - \lambda^k x_2 \sin(k \omega)
        + \hat{\lambda}^k q_k^{(1)}
      \right) (\tilde{y} - y^-) \\
      &\quad + O(1) \left(
        \lambda^k x_1 \sin(k \omega)
        + \lambda^k x_2 \cos(k \omega)
        + \hat{\lambda}^k q_k^{(2)}
      \right) (\tilde{y} - y^-) \\
      &= \lambda^k \hat{\alpha}_2 x_1 + \lambda^k \hat{\beta}_2 x_2
        + O(\lambda^k) x_1 (\tilde{y} - y^-)
        + O(\lambda^k) x_2 (\tilde{y} - y^-)
        + O((\tilde{y} - y^-)^2) \\
      &\quad + a_{21} \hat{\lambda}^k q_k^{(1)}
      + a_{22} \hat{\lambda}^k q_k^{(2)}
      + O(\lambda^{2k})
      + O(1) \hat{\lambda}^k q_k^{(1)} (\tilde{y} - y^-)
      + O(1) \hat{\lambda}^k q_k^{(2)} (\tilde{y} - y^-), \\
    \bar{y}
      &= \mu + c_1 \left(
        \lambda^k x_1 \cos(k \omega)
        - \lambda^k x_2 \sin(k \omega)
        + \hat{\lambda}^k q_k^{(1)}
      \right)
      + c_2 \left(
        \lambda^k x_1 \sin(k \omega)
        + \lambda^k x_2 \cos(k \omega)
        + \hat{\lambda}^k q_k^{(2)}
      \right) \\
      &\quad + d (\tilde{y} - y^-)^2
      + O(\lambda^{2k}) + O((\tilde{y} - y^-)^3) \\
      &\quad + O(1) \left(
        \lambda^k x_1 \cos(k \omega)
        - \lambda^k x_2 \sin(k \omega)
        + \hat{\lambda}^k q_k^{(1)}
      \right) (\tilde{y} - y^-) \\
      &\quad + O(1) \left(
        \lambda^k x_1 \sin(k \omega)
        + \lambda^k x_2 \cos(k \omega)
        + \hat{\lambda}^k q_k^{(2)}
      \right) (\tilde{y} - y^-) \\
      &= \mu + \lambda^k \alpha^* x_1 + \lambda^k \beta^* x_2
      + O(\lambda^k) x_1 (\tilde{y} - y^-) + O(\lambda^k) x_2 (\tilde{y} - y^-)
      + d(\tilde{y} - y^-)^2
      + O((\tilde{y} - y^-)^3) \\
      &\quad + c_1 \hat{\lambda}^k q_k^{(1)}
      + c_2 \hat{\lambda}^k q_k^{(2)}
      + O(\lambda^{2k})
      + O(1) \hat{\lambda}^k q_k^{(1)} (\tilde{y} - y^-)
      + O(1) \hat{\lambda}^k q_k^{(2)} (\tilde{y} - y^-).
  \end{align*}
}

Since the trigonometric functions are multiplied by $(\lambda(\varepsilon))^k$
in \eqref{eq-03-T0kprm},
$j$-th partial derivatives of
$O(\lambda^k)$ with respect to $\varepsilon$
have estimate of $O(k^j \lambda^k)$
for any $j \in \{1, 2, \cdots, r - 2\}$.
Since the partial derivatives of $q_k^{(i)}$ up to order $r-2$ are uniformly bounded,
$j$-th partial derivatives of
$O(\hat{\lambda}^k)$ with respect to $(x_1, x_2, \tilde{y}, \varepsilon)$
have estimate of $O(\hat{\lambda}^k)$
for any $j \in \{1, 2, \cdots, r - 2\}$.
The $j$-th partial derivatives of
$O((\tilde{y} - y^-)^i)$ ($i \in \{2, 3\}$)
with respect to $\varepsilon$
have estimate of $O((\tilde{y} - y^-)^i)$
for any $j \in \{1, 2, \cdots, r - 2\}$.
The $j$-th partial derivatives of
$O((\tilde{y} - y^-)^i)$ ($i \in \{2, 3\}$)
with respect to $\tilde{y}$
have estimate of $O((\tilde{y} - y^-)^{\max\{i - j, 0\}})$
for any $j \in \{1, 2, \cdots, r - 2\}$.
Note that these partial derivatives are uniformly bounded with respect to $k$,
see the definition of $O(\cdot)$ terms.

\medskip

\textbf{(2) Shilnikov coordinates.}\quad
Similar to \cite{GS1972},
we introduce the following `Shilnikov coordinates' on $\Pi_k$
(this terminology is from \cite{LLST2022}):
\begin{align}
  X_1 := x_1 - x_1^+, \quad
  X_2 := x_2 - x_2^+, \quad
  Y := \tilde{y} - y^-.
  \label{eq-03-Shilnikov}
\end{align}
The $\Pi_k$ is written as
\begin{align*}
    \Pi_k = [-\delta_\mathrm{dom}, \delta_\mathrm{dom}]^3
\end{align*}
in $(X_1, X_2, Y)$ coordinates.

We write $T_k: \Pi_k \ni (X_1, X_2, Y)
\mapsto (\bar{X}_1, \bar{X}_2, \bar{Y})$.
Applying \eqref{eq-03-Shilnikov} to \eqref{eq-03-comp},
we have
\begin{align*}
  \bar{X}_1 &= \lambda^k \hat{\alpha}_1 (X_1 + x_1^+)
    + \lambda^k \hat{\beta}_1 (X_2 + x_2^+)
    + b Y \\
    &\quad + O(\lambda^k) (X_1 + x_1^+) Y
    + O(\lambda^k) (X_2 + x_2^+) Y
    + O(Y^2)
    + O(\hat{\lambda}^k), \\
  \bar{X}_2 &= \lambda^k \hat{\alpha}_2 (X_1 + x_1^+)
    + \lambda^k \hat{\beta}_2 (X_2 + x_2^+) \\
    &\quad + O(\lambda^k) (X_1 + x_1^+) Y
    + O(\lambda^k) (X_2 + x_2^+) Y
    + O(Y^2)
    + O(\hat{\lambda}^k).
\end{align*}
From this,
we obtain
\begin{align}
\begin{aligned}
  \bar{X}_1 &= \lambda^k \hat{\alpha}_1 (X_1 + x_1^+)
    + \lambda^k \hat{\beta}_1 (X_2 + x_2^+)
    + b Y + \hat{\hat{h}}_1(X_1, X_2, Y, \varepsilon), \\
  \bar{X}_2 &= \lambda^k \hat{\alpha}_2 (X_1 + x_1^+)
    + \lambda^k \hat{\beta}_2 (X_2 + x_2^+)
    + \hat{\hat{h}}_2(X_1, X_2, Y, \varepsilon),
\end{aligned} \label{eq-03-SForm2.1}
\end{align}
where
\begin{align}
  \hat{\hat{h}}_1, \hat{\hat{h}}_2 = O(\lambda^k)(1 + X_1 + X_2)Y
    + O(Y^2) + O(\hat{\lambda}^k).
  \label{eq-03-evalhhh12org}
\end{align}
Note the estimate of the partial derivatives of
$O(\lambda^k)$, $O(Y^2)$, and $O(\hat{\lambda}^k)$
described in the end of the step~(1).
By partial differentiation of the above equation,
we obtain
\begin{align}
\begin{aligned}
  &\hat{\hat{h}}_{i, X_j}, \,
  \hat{\hat{h}}_{i, X_j X_l}, \,
  \hat{\hat{h}}_{i, X_j X_l X_m}, \,
  \hat{\hat{h}}_{i, X_j X_l Y}, \,
  \hat{\hat{h}}_{i, X_j Y Y}
  = O(\hat{\lambda}^k), \quad
  \hat{\hat{h}}_{i, Y}, \,
  \hat{\hat{h}}_{i, X_j Y}
  = O(\lambda^k), \\
  &\hat{\hat{h}}_{i, Y Y}, \,
  \hat{\hat{h}}_{i, Y Y Y}
  = O(1), \quad
  \hat{\hat{h}}_{i, \sigma}, \,
  \hat{\hat{h}}_{i, X_j \sigma}
  = O(\hat{\lambda}^k), \quad
  \hat{\hat{h}}_{i, Y \sigma}
  = O(k \lambda^k)
\end{aligned} \label{eq-03-evalhhh12}
\end{align}
for any $Y$ with $|Y| \leq O(\lambda^k)$,
$i$, $j$, $l$, $m \in \{1, 2\}$,
and $\sigma \in \{\mu, \omega, \rho\}$,
where we used \eqref{eq-03-hatC} to sort $O(\cdot)$ terms.

By \eqref{eq-03-T0kprm} and \eqref{eq-03-Shilnikov},
we get
\begin{align*}
  \bar{y} = \gamma^{-k}(\bar{Y} + y^-)
    + \hat{\gamma}^{-k}q_k^{(3)}(\bar{X}_1 + x_1^+,
      \bar{X}_2 + x_2^+, \bar{Y} + y^-, \varepsilon)
\end{align*}
in \eqref{eq-03-comp}.
Applying \eqref{eq-03-Shilnikov} to \eqref{eq-03-comp} again,
\begin{align}
\begin{aligned}
  \bar{Y} &= \gamma^k \mu - y^-
    + \lambda^k \gamma^k \alpha^* (X_1 + x_1^+)
    + \lambda^k \gamma^k \beta^* (X_2 + x_2^+) \\
    &\quad + \gamma^k O(\lambda^k) (X_1 + x_1^+) Y
    + \gamma^k O(\lambda^k) (X_2 + x_2^+) Y
    + \gamma^k d Y^2
    + \gamma^k O(Y^3) \\
    &\quad + \gamma^k O(\hat{\lambda}^k)
    - \gamma^k \hat{\gamma}^{-k}q_k^{(3)}(\bar{X}_1 + x_1^+,
    \bar{X}_2 + x_2^+, \bar{Y} + y^-, \varepsilon).
\end{aligned} \label{eq-03-SForm1.2}
\end{align}
Substituting \eqref{eq-03-SForm2.1} into the last term,
we can think of it as a function of $(X_1, X_2, Y, \bar{Y}, \varepsilon)$.
Using Proposition~\ref{pc1-sgl}
(see the appendix for the proof),
we can solve such a equation
with respect to $\bar{Y}$ as a function of $(X_1, X_2, Y, \varepsilon)$
since for the last term
$-\gamma^k \hat{\gamma}^{-k} q_k^{(3)}$
of the right-hand side of \eqref{eq-03-SForm1.2},
\begin{align*}
  -\gamma^k \hat{\gamma}^{-k} q_k^{(3)}, \quad
  \partial_{\bar{Y}}(-\gamma^k \hat{\gamma}^{-k} q_k^{(3)})
  = -\gamma^k \hat{\gamma}^{-k} q_{k, \tilde{y}}^{(3)}
\end{align*}
uniformly converge to 0 as $k \to \infty$.
Then, we obtain the solution
\begin{align*}
  \bar{Y} = \gamma^k \mu - y^-
  + \lambda^k \gamma^k \alpha^* (X_1 + x_1^+)
  + \lambda^k \gamma^k \beta^* (X_2 + x_2^+)
  + \gamma^k d Y^2
  + \hat{\hat{h}}_3(X_1, X_2, Y, \varepsilon),
\end{align*}
where
\begin{align}
  \hat{\hat{h}}_3 = \gamma^k \left(
    O(\lambda^k)(1 + X_1 + X_2)Y
    + O(Y^3) + O(\hat{\lambda}^k) \right)
    + \hat{\hat{h}}_3'(X_1, X_2, Y, \varepsilon), \quad
    \hat{\hat{h}}_3' = \gamma^k O(\hat{\gamma}^{-k}).
  \label{eq-03-hhh3}
\end{align}
Here, the estimate of the partial derivatives of
$O(\lambda^k)$, $O(Y^3)$, and $O(\hat{\lambda}^k)$
described in the end of the step~(1).
However, for the last term $\hat{\hat{h}}_3'$ in the above equation,
its partial derivatives up to order three have estimates of
\begin{align}
\begin{aligned}
  &\hat{\hat{h}}_{3, X_i}', \,
  \hat{\hat{h}}_{3, Y}', \,
  \hat{\hat{h}}_{3, X_i X_j}', \,
  \hat{\hat{h}}_{3, X_i Y}', \,
  \hat{\hat{h}}_{3, X_i X_j X_l}', \,
  \hat{\hat{h}}_{3, X_i X_j Y}'
  = \gamma^k O(\hat{\lambda}^k), \\
  &\hat{\hat{h}}_{3, Y Y}'
  = \gamma^k O(\gamma^k \hat{\gamma}^{-k}), \quad
  \hat{\hat{h}}_{3, X_i Y Y}', \,
  \hat{\hat{h}}_{3, Y Y Y}'
  = \gamma^k (O(\gamma^k \hat{\gamma}^{-k})
  + O(\lambda^k \gamma^{2k} \hat{\gamma}^{-k})), \\
  &\hat{\hat{h}}_{3, \sigma'}', \,
  \hat{\hat{h}}_{3, X_i \sigma'}', \,
  \hat{\hat{h}}_{3, Y \sigma'}'
  = \gamma^k O(\hat{\lambda}^k), \quad
  \hat{\hat{h}}_{3, \mu}'
  = \gamma^k O(\gamma^k \hat{\gamma}^{-k}), \\
  &\hat{\hat{h}}_{3, X_i \mu}', \,
  \hat{\hat{h}}_{3, Y \mu}'
  = \gamma^k (O(\gamma^k \hat{\gamma}^{-k})
  + O(\lambda^k \gamma^{2k} \hat{\gamma}^{-k}))
\end{aligned} \label{eq-03-evalGH}
\end{align}
for any $Y$ with $|Y| \leq O(\lambda^k)$,
$\mu$ with $|\mu| \leq O(\lambda^k) + O(\gamma^{-k})$,
$i$, $j$, $l \in \{1, 2\}$,
and $\sigma' \in \{\omega, \rho\}$,
see Remark~\ref{rem-hhathat3p} for more details.
Hence, by partial differentiation of the former equation in \eqref{eq-03-hhh3},
we have
\begin{align}
\begin{aligned}
  &\hat{\hat{h}}_{3, X_i}, \,
  \hat{\hat{h}}_{3, X_i X_j}, \,
  \hat{\hat{h}}_{3, X_i X_j X_l}, \,
  \hat{\hat{h}}_{3, X_i X_j Y}
  = \gamma^k O(\hat{\lambda}^k), \quad
  \hat{\hat{h}}_{3, Y}, \,
  \hat{\hat{h}}_{3, X_i Y}
  = \gamma^k O(\lambda^k), \\
  &\hat{\hat{h}}_{3, Y Y}
  = \gamma^k O(\gamma^k \hat{\gamma}^{-k}), \quad 
  \hat{\hat{h}}_{3, X_i Y Y}
  = \gamma^k (O(\gamma^k \hat{\gamma}^{-k})
  + O(\lambda^k \gamma^{2k} \hat{\gamma}^{-k})), \quad
  \hat{\hat{h}}_{3, Y Y Y}
  = \gamma^k O(1), \\
  &\hat{\hat{h}}_{3, \sigma'}, \,
  \hat{\hat{h}}_{3, X_i \sigma'}
    = \gamma^k O(\hat{\hat{\lambda}}^k), \quad
  \hat{\hat{h}}_{3, Y \sigma'}
  = \gamma^k O(k \lambda^k), \quad
  \hat{\hat{h}}_{3, \mu}
  = \gamma^k O(\gamma^k \hat{\gamma}^{-k}), \\
  &\hat{\hat{h}}_{3, X_i \mu}, \,
  \hat{\hat{h}}_{3, Y \mu}
  = \gamma^k (O(\gamma^k \hat{\gamma}^{-k})
  + O(\lambda^k \gamma^{2k} \hat{\gamma}^{-k}))
\end{aligned} \label{eq-03-peval_hhh3}
\end{align}
for any $Y$ with $|Y| \leq O(\lambda^k)$,
$\mu$ with $|\mu| \leq O(\lambda^k) + O(\gamma^{-k})$,
$i$, $j$, $l \in \{1, 2\}$,
and $\sigma' \in \{\omega, \rho\}$,
where $\hat{\hat{\lambda}}$ is a constant with \eqref{e31-hhlmd}.
\debugtext{
  Hint:
  Let
  \begin{align*}
    E := \gamma^k \left(
      O(\lambda^k)(1 + X_1 + X_2)Y
      + O(Y^3) + O(\hat{\lambda}^k)
    \right).
  \end{align*}
  Then,
  \begin{align*}
    &E_{X_j}, \,
    E_{X_j X_l}, \,
    E_{X_j X_l X_m}, \,
    E_{X_j X_l Y}, \,
    E_{X_j Y Y}
    = \gamma^k O(\hat{\lambda}^k), \quad
    E_{Y}, \,
    E_{X_j Y}, \,
    E_{Y Y}
    = \gamma^k O(\lambda^k), \\
    &E_{Y Y Y}
    = \gamma^k O(1), \quad
    E_{\sigma}, \,
    E_{X_j \sigma}
    = \gamma^k O(\hat{\hat{\lambda}}^k), \quad
    E_{Y \sigma}
    = \gamma^k O(k \lambda^k)
  \end{align*}
  for any $Y$ with $|Y| \leq O(\lambda^k)$,
  $i \in \{1, 2\}$,
  $\sigma \in \{\mu, \omega, \rho\}$.
}

\medskip

\textbf{(3) Shift.}\quad
Note that since $(f, \Gamma)$ holds \textbf{(QC)},
we may suppose $d \neq 0$ for any $\varepsilon \in R_\mathrm{prm}$
by replacing $\delta_\mathrm{prm}$ with a smaller one
according to Remark~\ref{rem-smaller}.
Consider the system of equations
\begin{align*}
  X_1 = H_k^{(1)}(\mu, \omega, \rho, X_1, X_2, Y), \quad
  X_2 = H_k^{(2)}(\mu, \omega, \rho, X_1, X_2, Y), \quad
  Y = H_k^{(3)}(\mu, \omega, \rho, X_1, X_2, Y)
\end{align*}
for $Y$ with $|Y| \leq O(\lambda^k)$
and $\mu$ with $|\mu| \leq O(\lambda^k) + O(\gamma^{-k})$, where
\begin{align*}
  &H_k^{(3)} := -(2 \gamma^k d)^{-1}
    \hat{\hat{h}}_{3, Y}(X_1, X_2, Y, \varepsilon), \\
  &H_k^{(1)} := \lambda^k \hat{\alpha}_1 (X_1 + x_1^+)
    + \lambda^k \hat{\beta}_1 (X_2 + x_2^+)
    + b H_k^{(3)} + \hat{\hat{h}}_1(X_1, X_2, H_k^{(3)}, \varepsilon), \\
  &H_k^{(2)} := \lambda^k \hat{\alpha}_2 (X_1 + x_1^+)
    + \lambda^k \hat{\beta}_2 (X_2 + x_2^+)
    + \hat{\hat{h}}_2(X_1, X_2, H_k^{(3)}, \varepsilon).
\end{align*}
By \eqref{eq-03-evalhhh12org},
\eqref{eq-03-evalhhh12},
and \eqref{eq-03-peval_hhh3},
we have
\begin{align*}
  &H_k^{(i)}
    = O(\lambda^k), \quad
  H_{k, \mu}^{(1)}, \,
  H_{k, \mu}^{(3)}
    = O(\gamma^k \hat{\gamma}^{-k})
    + O(\lambda^k \gamma^{2k} \hat{\gamma}^{-k}), \quad
  H_{k, \mu}^{(2)}, \,
  H_{k, \sigma'}^{(i)}
    = O(k \lambda^k), \\
  &H_{k, X_j}^{(i)}
    = O(\lambda^k), \quad
  H_{k, Y}^{(1)}, \,
  H_{k, Y}^{(3)}
    = O(\gamma^k \hat{\gamma}^{-k}), \quad
  H_{k, Y}^{(2)}
    = O(\hat{\lambda}^k)
\end{align*}
for any $i \in \{1, 2, 3\}$,
$j \in \{1, 2\}$,
and $\sigma' \in \{\omega, \rho\}$.
Thus, the Proposition~\ref{pc3-sys} solves the above system of equations
and we get the solutions $(X_1, X_2, Y)
= (X_{1, k}^*(\varepsilon), X_{2,k}^*(\varepsilon), Y_k^*(\varepsilon))
= (X_{1, k}^*, X_{2, k}^*, Y_k^*)$
that are $C^{r - 2}$ with respect to $\varepsilon$ and
\begin{align}
\begin{aligned}
  &X_{1, k}^*, \,
  X_{2, k}^*, \,
  Y_k^*
    = O(\lambda^k), \quad
  X_{1, k, \mu}^*, \,
  Y_{k, \mu}^*
    = O(\gamma^k \hat{\gamma}^{-k})
    + O(\lambda^k \gamma^{2k} \hat{\gamma}^{-k}), \\
  &X_{1, k, \sigma'}^*
    = O(\gamma^k \hat{\gamma}^{-k}), \quad
  X_{2, k, \mu}^*, \,
  X_{2, k, \sigma'}^*, \,
  Y_{k, \sigma'}^*
    = O(k \lambda^k)
\end{aligned} \label{eq-03-sol_ev}
\end{align}
for any $\sigma' \in \{\omega, \rho\}$.

We define the new coordinates
\begin{align}
  X_1^{new} := X_1 - X_{1,k}^*(\varepsilon), \quad
  X_2^{new} := X_2 - X_{2,k}^*(\varepsilon), \quad
  Y_2^{new} := Y - Y_k^*(\varepsilon).
  \label{e32-XYshift}
\end{align}
Then,
by dropping `new',
$T_k$ has the form
\begin{align}
\begin{aligned}
  \bar{X}_1 &= \lambda^k \hat{\alpha}_1 X_1
    + \lambda^k \hat{\beta}_1 X_2
    + b Y + \hat{h}_1(X_1, X_2, Y, \varepsilon), \\
  \bar{X}_2 &= \lambda^k \hat{\alpha}_2 X_1
    + \lambda^k \hat{\beta}_2 X_2
    + \hat{h}_2(X_1, X_2, Y, \varepsilon), \\
  \bar{Y} &= \hat{\mu}
    + \lambda^k \gamma^k \alpha^* X_1
    + \lambda^k \gamma^k \beta^* X_2
    + \gamma^k d Y^2
    + \hat{h}_3(X_1, X_2, Y, \varepsilon),
\end{aligned} \label{eq-03-ShiftForm}
\end{align}
where
\begin{align}
\begin{aligned}
  &\hat{h}_1
    = \hat{\hat{h}}_1(X_1 + X_{1, k}^*, X_2 + X_{2, k}^*, Y + Y_k^*, \varepsilon)
    - \hat{\hat{h}}_1(X_{1, k}^*, X_{2, k}^*, Y_k^*, \varepsilon), \\
  &\hat{h}_2
    = \hat{\hat{h}}_2(X_1 + X_{1, k}^*, X_2 + X_{2, k}^*, Y + Y_k^*, \varepsilon)
    - \hat{\hat{h}}_2(X_{1, k}^*, X_{2, k}^*, Y_k^*, \varepsilon), \\
  &\hat{h}_3
    = 2 \gamma^k d Y Y_k^*
    + \hat{\hat{h}}_3(X_1 + X_{1, k}^*, X_2 + X_{2, k}^*, Y + Y_k^*, \varepsilon)
    - \hat{\hat{h}}_3(X_{1, k}^*, X_{2, k}^*, Y_k^*, \varepsilon),
\end{aligned} \label{eq-03-hh123_def}
\end{align}
and
\begin{align*}
  \hat{\mu}
    &= \gamma^k \mu - y^- + \lambda^k \gamma^k (\alpha^* x_1^+ + \beta^* x_2^+) \\
    &\quad + \gamma^k \left(
      -\gamma^{-k} Y_k^* + \lambda^k (\alpha^* X_{1, k}^* + \beta^* X_{2, k}^*)
      + d (Y_k^*)^2
      + \gamma^{-k} \hat{\hat{h}}_3(X_{1, k}^*, X_{2, k}^*, Y_k^*, \varepsilon)
    \right).
\end{align*}
The $\Pi_k$ is given as
\begin{align*}
  \Pi_k = [-\delta_\mathrm{dom}, \delta_\mathrm{dom}]^3
    - (X_{1, k}^*, X_{2, k}^*, Y_k^*)
\end{align*}
in the new coordinates.
By \eqref{eq-03-hatC},  
\eqref{eq-03-peval_hhh3} implies \eqref{eq-03-muhat} with the desired estimate.
Moreover, \eqref{eq-03-evalhhh12org} and \eqref{eq-03-hhh3} yield
\begin{align}
\begin{aligned}
  \hat{h}_i &= O(\hat{\lambda}^k)(X_1 + X_2) + O(\lambda^k)Y + O(Y^2), \\
  \hat{h}_3 &= \gamma^k \left(
    O(\hat{\lambda}^k)(X_1 + X_2) + O(\lambda^k)(X_1 + X_2)Y
    + O(\lambda^k) Y^2 + O(Y^3)
\right)
\end{aligned} \label{eq-03-hh123}
\end{align}
for any $i \in \{1, 2\}$,
where $O(\lambda^k)$, $O(Y^2)$, $O(Y^3)$, and $O(\hat{\lambda}^k)$
are now different from the ones at the end of the step~(1).
\debugtext{
  Hint:
  We have
  \begin{align*}
    \hat{h}_i &= -\hat{\hat{h}}_i(X_{1, k}^*, X_{2, k}^*, Y_k^*, \varepsilon)
    + O(\lambda^k)(1 + X_1 + X_{1, k}^* + X_2 + X_{2, k}^*)(Y + Y_k^*)
    + O((Y + Y_k^*)^2) + O(\hat{\lambda}^k) \\
    &= O(\lambda^k) Y + O(Y^2) + O(\hat{\lambda}^k).
  \end{align*}
  On the other hand,  
  setting $Y = 0$ yields  
  \begin{align*}
    \hat{h}_i = O(\hat{\lambda}^k).  
  \end{align*}
}
By partial differentiation of \eqref{eq-03-hh123_def},
\eqref{eq-03-evalhhh12},
\eqref{eq-03-peval_hhh3},
and \eqref{eq-03-sol_ev} yield
\begin{align}
\begin{aligned}
  &\hat{h}_{i, X_j}, \,
  \hat{h}_{i, X_j X_l}, \,
  \hat{h}_{i, X_j X_l X_m}, \,
  \hat{h}_{i, X_j X_l Y}, \,
  \hat{h}_{i, X_j Y Y}
    = O(\hat{\lambda}^k), \quad
  \hat{h}_{i, Y}, \,
  \hat{h}_{i, X_j Y}
    = O(\lambda^k), \\
  &\hat{h}_{i, Y Y}, \,
  \hat{h}_{i, Y Y Y}
    = O(1), \quad
  \hat{h}_{i, \sigma}, \,
  \hat{h}_{i, X_j \sigma}
    = O(\hat{\lambda}^k), \\
  &\hat{h}_{i, Y \mu}
    = O(\gamma^k \hat{\gamma}^{-k})
    + O(\lambda^k \gamma^{2k} \hat{\gamma}^{-k}), \quad
  \hat{h}_{i, Y \sigma'}
    = O(k \lambda^k)
\end{aligned} \label{eq-03-hh12Higher}
\end{align}
and
\begin{align}
\begin{aligned}
  &\hat{h}_{3, X_i}, \,
  \hat{h}_{3, X_i X_j}, \,
  \hat{h}_{3, X_i X_j X_l}, \,
  \hat{h}_{3, X_i X_j Y}
    = \gamma^k O(\hat{\lambda}^k), \quad
  \hat{h}_{3, Y}, \,
  \hat{h}_{3, X_i Y}
    = \gamma^k O(\lambda^k), \\
  &\hat{h}_{3, Y Y}
    = \gamma^k O(\gamma^k \hat{\gamma}^{-k}), \quad 
  \hat{h}_{3, X_i Y Y}
    = \gamma^k (O(\gamma^k \hat{\gamma}^{-k})
    + O(\lambda^k \gamma^{2k} \hat{\gamma}^{-k})), \quad
  \hat{h}_{3, Y Y Y}
    = \gamma^k O(1), \\
  &\hat{h}_{3, \mu}
    = \gamma^k O(\gamma^k \hat{\gamma}^{-k}), \quad
  \hat{h}_{3, \sigma'}, \,
  \hat{h}_{3, X_i \sigma'}
    = \gamma^k O(\hat{\hat{\lambda}}^k), \\
  &\hat{h}_{3, X_i \mu}, \,
  \hat{h}_{3, Y \mu}
    = \gamma^k (O(\gamma^k \hat{\gamma}^{-k})
    + O(\lambda^k \gamma^{2k} \hat{\gamma}^{-k})), \quad
  \hat{h}_{3, Y \sigma'}
    = \gamma^k O(k \lambda^k)
\end{aligned} \label{eq-03-hh3Higher}
\end{align}
for any $Y$ with $|Y| \leq O(\lambda^k)$,
$\mu$ with $|\mu| \leq O(\lambda^k) + O(\gamma^{-k})$,
$i$, $j$, $l$, $m \in \{1, 2\}$,
$\sigma \in \{\mu, \omega, \rho\}$,
and $\sigma' \in \{\omega, \rho\}$.

\medskip

\textbf{(4) Normal form.}\quad
For any
$k \in \mathbb{Z}_{> \kappa(\delta_\mathrm{dom})}$
and $\varepsilon \in R_k^\mathrm{bd}$,
we introduce the new coordinates
\begin{align}
  Z := \alpha^*_k(\varepsilon) X_1 + \beta^*_k(\varepsilon) X_2, \quad
  W := X_2
  \label{e32-ZW}
\end{align}
on $\Pi_k$.
The $\Pi_k$ is given by
\begin{align*}
  \Pi_k &=
    \{(Z, W) \:|\: |W| \leq \delta_\mathrm{dom}, \,
      |Z - \beta^* W| \leq |\alpha^*| \delta_\mathrm{dom}\}
    \times [-\delta_\mathrm{dom}, \delta_\mathrm{dom}] \\
    &\quad - (\alpha_k^* X_{1, k}^* + \beta_k^* X_{2, k}^*,
    X_{2, k}^*, Y_k^*)
\end{align*}
in $(Z, Y, W)$ coordinates.
By the note after \eqref{e32-sinE}
and \eqref{e32-bsetc} of $\beta^*_k$,
there exist constants $C_1 = C_1(\mathbb{F}) > 0$
and $C_2 = C_2(\mathbb{F}) > 0$ such that
\begin{align*}
  |\alpha^*_k(\varepsilon)| \geq C_1, \quad
  |\beta^*_k(\varepsilon)| \leq C_2
\end{align*}
for any $\varepsilon \in R_k^\mathrm{bd}$.
Defining  
\begin{align}
  \delta_\mathrm{dom}' = \delta_\mathrm{dom}'(\delta_\mathrm{dom}; \mathbb{F})
    := \frac{1}{2} \delta_\mathrm{dom} \min \left\{  
    1, \, \frac{C_1}
      {1 + C_2}  
  \right\},  
  \label{e32-ddompD}
\end{align}  
we can verify that  
\begin{align*}
  \Pi_k' := [-\delta_\mathrm{dom}', \delta_\mathrm{dom}']^3  
\end{align*}  
in $(Z, Y, W)$ coordinates
is contained in $\Pi_k$
by replacing $\kappa(\delta_\mathrm{dom})$ with a larger one  
according to Remark~\ref{rem-smaller}.
The definition of $\delta_\mathrm{dom}'$ implies \eqref{e32-ddomp}.
We can rewrite \eqref{eq-03-ShiftForm} as \eqref{eq-03-NF},
where
\begin{align*}
  &\alpha_1 = \alpha_1^{(k)}(\varepsilon)
    := \hat{\alpha}_1
    + \hat{\alpha}_2
    \frac{\beta^*}{\alpha^*}, \quad
  \beta_1 = \beta_1^{(k)}(\varepsilon)
    := -\hat{\alpha}_1 \beta^*
    + \hat{\beta}_1 \alpha^*
    - \hat{\alpha}_2
    \frac{(\beta^*)^2}{\alpha^*}
    + \hat{\beta}_2 \beta^*, \\
  &\alpha_3 = \alpha_3^{(k)}(\varepsilon)
    := \frac{\hat{\alpha}_2}{\alpha^*}, \quad
  \beta_3 = \beta_3^{(k)}(\varepsilon)
    := -\hat{\alpha}_2
    \frac{\beta^*}{\alpha^*}
    + \hat{\beta}_2
\end{align*}
and
\begin{align*}
  h_1 = \alpha^* \hat{h}_1(M, \varepsilon)
    + \beta^* \hat{h}_2(M, \varepsilon), \quad
  h_2 = \hat{h}_3(M, \varepsilon), \quad
  h_3 = \hat{h}_2(M, \varepsilon), \quad
  M = (\frac{1}{\alpha^*}Z - \frac{\beta^*}{\alpha^*} W, W, Y).
\end{align*}
Note that for the quantities in \eqref{e32-alp} and \eqref{e32-bsetc},
we have
\begin{align*}
  &\alpha^*_k, \,
  \frac{1}{\alpha^*_k}, \,
  \beta^*_k, \,
  \hat{\alpha}^{(k)}_i, \,
  \hat{\beta}^{(k)}_i = O(1), \quad
  \alpha^*_{k, \tilde{\sigma}}, \,
  \partial_{\tilde{\sigma}} \left(\frac{1}{\alpha^*_k}\right), \,
  \beta^*_{k, \tilde{\sigma}}, \,
  \hat{\alpha}^{(k)}_{i, \tilde{\sigma}}, \,
  \hat{\beta}^{(k)}_{i, \tilde{\sigma}}
    = O(1), \\
  &\alpha^*_{k, \omega}, \,
  \partial_{\omega} \left(\frac{1}{\alpha^*_k}\right), \,
  \beta^*_{k, \omega}, \,
  \hat{\alpha}^{(k)}_{i, \omega}, \,
  \hat{\beta}^{(k)}_{i, \omega}
    = O(k)
\end{align*}
for any $\tilde{\sigma} \in \{\mu, \rho\}$,
and $i \in \{1, 2\}$.
The above formula implies \eqref{eq-03-ab_ev}.
The formula \eqref{eq-03-hh123} yields \eqref{eq-03-NF_higher1}.
By partial differentiation of the above equations,
\eqref{eq-03-hh12Higher} and \eqref{eq-03-hh3Higher} imply
\eqref{eq-03-NF_higher2} and \eqref{eq-03-NF_higher3}.
We complete the proof.
\end{proof}

\begin{remark} \label{rem-hhathat3p}
  We explain how to get \eqref{eq-03-evalGH}.
  Let
  \begin{align*}
    F_k(X_1, X_2, Y, \bar{Y}, \varepsilon)
      := (\text{left-hand side of \eqref{eq-03-SForm1.2}})
      - (\text{right-hand side of \eqref{eq-03-SForm1.2}}).
  \end{align*}
  Let $G_k(X_1, X_2, Y, \varepsilon)$ be the part
  from $\gamma^k \mu$ to $\gamma^k O(\hat{\lambda}^k)$
  on the right-hand side of \eqref{eq-03-SForm1.2}
  and $H_k(X_1, X_2, Y, \bar{Y}, \varepsilon)$ be
  the last term of the right-hand side of \eqref{eq-03-SForm1.2}.
  By the definitions of $G_k$,
  we have
  \begin{align*}
    &G_{k, X_i}
      = O(\lambda^k \gamma^k), \quad
    G_{k, Y}
      = O(\lambda^k \gamma^k) + O(\gamma^k Y), \quad
    G_{k, X_i X_j}
      = O(\hat{\lambda}^k \gamma^k), \quad
    G_{k, X_i Y}
      = O(\lambda^k \gamma^k), \\
    &G_{k, Y Y}
      = O(\gamma^k), \quad
    G_{k, X_i X_j X_l}, G_{k, X_i X_j Y}, G_{k, X_i Y Y}
      = O(\hat{\lambda}^k \gamma^k), \quad
    G_{k, Y Y Y}
      = O(\gamma^k)
  \end{align*}
  for any $i$, $j$, $l \in \{1, 2\}$ and
  \begin{align*}
    &G_{k, \mu}
      = O(k \gamma^k)(k^{-1} + \mu + Y^2), \quad
    G_{k, \sigma'}
      = O(k \gamma^k) (k^{-1} \gamma^{-k} + \mu + Y + \lambda^k), \\
    &G_{k, X_i \sigma}
      = O(k \lambda^k \gamma^k), \quad
    G_{k, Y \sigma}
      = O(k \lambda^k \gamma^k) + O(k \gamma^k Y)
  \end{align*}
  for any $\sigma' \in \{\omega, \rho\}$,
  any $\sigma \in \{\mu, \omega, \rho\}$,
  and any $i \in \{1, 2\}$.
  We assume that $Y$ and $\mu$ vary under the conditions  
  $|Y| \leq O(\lambda^k)$ and $|\mu| \leq O(\lambda^k) + O(\gamma^{-k})$.  
  Then, the $\mathcal{G}_k$ defined in Section~\ref{sC2-ev} are
  \begin{align*}
    &\mathcal{G}_k^{(X_i X_j)}, \,
    \mathcal{G}_k^{(X_i X_j X_l)}, \,
    \mathcal{G}_k^{(X_i X_j Y)}, \,
    \mathcal{G}_k^{(X_i Y Y)}
      = O(\hat{\lambda}^k \gamma^k), \quad
    \mathcal{G}_k^{(X_i)}, \,
    \mathcal{G}_k^{(Y)}, \,
    \mathcal{G}_k^{(X_i Y)}
      = O(\lambda^k \gamma^k), \\
    &\mathcal{G}_k^{(Y Y)}, \,
    \mathcal{G}_k^{(Y Y Y)}
      = O(\gamma^k), \quad
    \mathcal{G}_k^{(\mu)}
      = O(\gamma^k), \\
    &\mathcal{G}_k^{(\sigma')}
      = O(k \lambda^k \gamma^k) + O(k), \quad
    \mathcal{G}_k^{(X_i \sigma)}, \,
    \mathcal{G}_k^{(Y \sigma)}
      = O(k \lambda^k \gamma^k)
  \end{align*}
  for any $i$, $j$, $l \in \{1, 2\}$,
  $\sigma' \in \{\omega, \rho\}$,
  and $\sigma \in \{\mu, \omega, \rho\}$.
  By the definition of $H_k$,
  the $\mathcal{H}_k$ defined in Section~\ref{sC2-ev} are
  \begin{align*}
    &\mathcal{H}_k^{(X_i)}, \,
    \mathcal{H}_k^{(Y)}, \,
    \mathcal{H}_k^{(X_i X_j)}, \,
    \mathcal{H}_k^{(X_i Y)}, \,
    \mathcal{H}_k^{(Y Y)}, \,
    \mathcal{H}_k^{(X_i X_j X_l)}, \,
    \mathcal{H}_k^{(X_i X_j Y)}, \,
    \mathcal{H}_k^{(X_i Y Y)}, \,
    \mathcal{H}_k^{(Y Y Y)}
      = O(\gamma^k \hat{\gamma}^{-k}) \\
    &\mathcal{H}_k^{(\sigma)}, \,
    \mathcal{H}_k^{(X_i \sigma)}, \,
    \mathcal{H}_k^{(Y \sigma)}
    = O(k \gamma^k \hat{\gamma}^{-k})
  \end{align*}
  for any $i$, $j$, $l \in \{1, 2\}$,
  and $\sigma \in \{\mu, \omega, \rho\}$.
  By \eqref{eq-03-hatC},
  Proposition~\ref{pc2-ev} implies
  the desired estimate \eqref{eq-03-evalGH}.
\end{remark}

\subsection{Invariant cone fields}\label{s33-cone}
Recall the $\delta_\mathrm{dom}'$ and $\Pi_k'$
in Proposition~\ref{p32-NF},
the range of parameters $R_k^\mathrm{bd}$ in \eqref{e32-Rkbd},
and the tuple of the core objects $\mathbb{F}$ in \eqref{eq-03-mathbbF}.
We shall think of the domain of the first-return map $T_k$ as
$\Pi_k' = [-\delta_\mathrm{dom}', \delta_\mathrm{dom}']^3$
in the $(Z, Y, W)$-space.
We use $(z, y, w)$ to denote vectors in the tangent spaces.

\begin{proposition}[Existence of cone fields]\label{p33-CF}
  By replacing  
  $\hat{\delta}_\mathrm{dom} > 0$ with a smaller one and  
  $\kappa(\delta_\mathrm{dom})$ with a larger one  
  according to Remark~\ref{rem-smaller},
  there exists $K = K(\mathbb{F}) > 0$ such that
  the following statements hold
  for any $k \in \mathbb{Z}_{> \kappa(\delta_\mathrm{dom})}$,
  $\varepsilon \in R_k^\mathrm{bd}$,
  and $\delta_\mathrm{dom} \in (0, \hat{\delta}_\mathrm{dom})$:    
  \begin{enumerate}
    \item The cone field in $\Pi_k'$
      \begin{align}
        \mathcal{C}^{ss}(Z, Y, W)
          = \{ (z, y, w) \:|\: |z| + |y| < K \delta_\mathrm{dom} |w| \}
        \label{eq-03-coness}
      \end{align}
      is backward-invariant, in other words,
      if $\bar{M} \in \Pi_k'$ with $M = T^{-1}_k(\bar{M}) \in \Pi_k'$,
      then
      \begin{align*}
        D(T^{-1}_k)_{\bar{M}}(\mathcal{C}^{ss}(\bar{M}))
        \subset \mathcal{C}^{ss}(M).
      \end{align*}
    \item The cone field
      \begin{align}
        \mathcal{C}^{cu}(Z, Y, W)
          = \left\{ (z, y, w) \:\middle|\:
          |w| < K \left((|Y| + \lambda^k)|z| + |\gamma|^{-k} |y|\right) \right\}
        \label{eq-03-conecu}
      \end{align}
      is forward-invariant, in other words,
      if $M \in \Pi_k'$ with $\bar{M} = T_k(M) \in \Pi_k'$, then
      \begin{align*}
        D(T_k)_M(\mathcal{C}^{cu}(M)) \subset \mathcal{C}^{cu}(\bar{M}).
      \end{align*}
  \end{enumerate}
\end{proposition}

\begin{remark}
  As in Remark~\ref{rem-proof},
  the proof of the above lemma can be found in \cite{LLST2022}.
  For the sake of completeness, the proof is given below.
\end{remark}

\begin{proof}[Proof of Proposition~\ref{p33-CF}]
By the normal form of $T_k$ in \eqref{eq-03-NF},
if we put
$D(T_k)_M: (z, y, w) \mapsto (\bar{z}, \bar{y}, \bar{w})$,
then
\begin{align}
\begin{aligned}
  \bar{z} &= O(\lambda^k) z
    + \left( -E_k + O(\lambda^k) + O(Y) \right) y + O(\lambda^k) w, \\
  \bar{y} &= \gamma^k \left( \lambda^k + O(\lambda^k) \right) z
    + \gamma^k \left(O(Y) + O(\lambda^k) \right) y
    + \gamma^k \left( O(\lambda^k) Y + O(\hat{\lambda}^k) \right) w, \\
  \bar{w} &= O(\lambda^k) z + \left( O(\lambda^k) + O(Y) \right) y + O(\lambda^k) w
\end{aligned} \label{eq-03-DTkM}
\end{align}
for any $M \in \Pi_k'$,
where $E_k$ is the quantity defined in \eqref{eq-03-Ek}.

In the following, we often replace $\hat{\delta}_\mathrm{dom}$
and $\kappa(\delta_\mathrm{dom})$ with smaller and larger ones, respectively,
but we replace them according to the rules in Remark~\ref{rem-smaller}.

\medskip

\textbf{(1) The cone field $\mathcal{C}^{ss}$.}\quad
Choose any $K > 0$ and define $\mathcal{C}^{ss}$ by \eqref{eq-03-coness}.
Let $(\bar{z}, \bar{y}, \bar{w}) \in \mathcal{C}^{ss}(\bar{M})$
and $(z, y, w) = D(T_k^{-1})_{\bar{M}}(\bar{z}, \bar{y}, \bar{w})$,
where $\bar{M} \in \Pi_k'$ with $M = T_k^{-1}(\bar{M}) \in \Pi_k'$.

The equation of $\bar{w}$ in \eqref{eq-03-DTkM} implies
\begin{align}
  |\bar{w}| < C_1 \lambda^k (|z| + |w|) + C_1 \delta_\mathrm{dom}|y|
  \label{eq-03-bwineq}
\end{align}
for some $C_1 = C_1(\mathbb{F})$ independent of $K$
by replacing $\hat{\delta}_\mathrm{dom}$
and $\kappa(\delta_\mathrm{dom})$ with smaller and larger ones.
The equation of $\bar{z}$ in \eqref{eq-03-DTkM},
$|\bar{z}| < K \delta_\mathrm{dom} |\bar{w}|$,
\eqref{eq-03-bwineq},
and the fact that $|E_k|$ is bounded away from zero yield
\begin{align}
  |y| < C_2 \lambda^k (|z| + |w|)
  \label{eq-03-yineq}
\end{align}
for some $C_2 = C_2(\mathbb{F})$ independent of $K$
by replacing $\hat{\delta}_\mathrm{dom}$ with smaller one
$\hat{\delta}_\mathrm{dom}^{new}(K)$.
Note that the new one does depend on $K$.
In the following, we drop the `new'.
The equation of $\bar{y}$ in \eqref{eq-03-DTkM},
$|\bar{y}| < K \delta_\mathrm{dom} |\bar{w}|$,
\eqref{eq-03-bwineq},
and \eqref{eq-03-yineq} yield
\begin{align}
  |z| < C_3 \delta_\mathrm{dom} |w|
  \label{eq-03-zineq}
\end{align}
for some $C_3 = C_3(\mathbb{F})$ independent of $K$
by replacing $\hat{\delta}_\mathrm{dom}(K)$
and $\kappa(\delta_\mathrm{dom})$ with smaller and larger ones.
Note that the new $\kappa(\delta_\mathrm{dom})$ does depend on $K$:
$\kappa(\delta_\mathrm{dom}) = \kappa(\delta_\mathrm{dom}, K)$.
The equations \eqref{eq-03-yineq} and \eqref{eq-03-zineq} imply
\begin{align}
  |y| < C_4 \lambda^k |w|, \quad
  |z| + |y| < C_4 \delta_\mathrm{dom} |w|
  \label{eq-03-yzyineq}
\end{align}
for some $C_4 = C_4(\mathbb{F})$ independent of $K$
by replacing $\hat{\delta}_\mathrm{dom}(K)$
and $\kappa(\delta_\mathrm{dom}, K)$ with smaller and larger ones.
Taking $K$ greater than $C_4$ completes the proof of invariance of
$\mathcal{C}^{ss}$.

\medskip

\textbf{(2) The cone field $\mathcal{C}^{cu}$.}\quad
Forget $K$ in Step~(1) for the moment
and choose any $K > 0$ and define $\mathcal{C}^{cu}$ by \eqref{eq-03-conecu}.
Let $(z, y, w) \in \mathcal{C}^{cu}(M)$
and $(\bar{z}, \bar{y}, \bar{w}) = D(T_k)_M(z, y, w)$,
where $M \in \Pi_k'$ with $\bar{M} = T_k(M) \in \Pi_k'$.

Since $(z, y, w) \in \mathcal{C}^{cu}(M)$,
\begin{align}
  |w| < |z| + |y|
  \label{eq-03-wineq}
\end{align}
by replacing $\hat{\delta}_\mathrm{dom}(K)$
and $\kappa(\delta_\mathrm{dom}, K)$ with smaller and larger ones.
The equation of $\bar{w}$ in \eqref{eq-03-DTkM} and \eqref{eq-03-wineq} imply
\begin{align}
  |\bar{w}| < C_5 \lambda^k |z| + C_5 (\lambda^k + |Y|) |y|
  \label{eq-03-bwineq2}
\end{align}
for some $C_5 = C_5(\mathbb{F})$.
The equations of $\bar{z}$ and $\bar{y}$ in \eqref{eq-03-DTkM},
\eqref{eq-03-wineq},
and the fact that $|E_k|$ is bounded away from zero yield
\begin{align*}
  |y| < C_6 |\bar{z}| + C_6 \lambda^k |z|, \quad
  |z| < C_6 \lambda^{-k} |\gamma|^{-k} |\bar{y}|
    + C_6(\lambda^{-k}|Y| + 1) |y|
\end{align*}
for some $C_6 = C_6(\mathbb{F})$
by replacing $\hat{\delta}_\mathrm{dom}(K)$
and $\kappa(\delta_\mathrm{dom}, K)$ with smaller and larger ones.
Thus, we have
\begin{align*}
  |y| < C_7 |\bar{z}| + C_7 |\gamma|^{-k} |\bar{y}|, \quad
  |z| < C_7 \lambda^{-k} |\gamma|^{-k} |\bar{y}|
    + C_7 (\lambda^{-k} |Y| + 1) |\bar{z}|
\end{align*}
for some $C_7 = C_7(\mathbb{F})$
by replacing $\hat{\delta}_\mathrm{dom}(K)$
and $\kappa(\delta_\mathrm{dom}, K)$ with smaller and larger ones.
This and \eqref{eq-03-bwineq2} further yield
\begin{align*}
  |\bar{w}| < C_8(|Y| + \lambda^k) |\bar{z}| + C_8 |\gamma|^{-k} |\bar{y}|
\end{align*}
for some $C_8 = C_8(\mathbb{F})$
by replacing $\hat{\delta}_\mathrm{dom}(K)$
and $\kappa(\delta_\mathrm{dom}, K)$ with smaller and larger ones.
Taking $K$ greater than $C_8$ completes the proof of invariance of
$\mathcal{C}^{cu}$.
\end{proof}

\begin{remark} \label{rem-WContract}
  Substituting \eqref{eq-03-zineq}
  and the former inequality in \eqref{eq-03-yzyineq}
  into \eqref{eq-03-bwineq},
  we get
  \begin{align*}
    |\bar{w}| = O(\lambda^k) |w|.
  \end{align*}
  Thus, vectors in $\mathcal{C}^{ss}$ are uniformly contracted by $DT_k$.
  This will be used in the proof of Proposition~\ref{p34-Qk}.
\end{remark}

\subsection{Non-hyperbolic periodic points} \label{s34-Qk}
In this section,
we prove the existence of the non-hyperbolic fixed point $Q_k$,
or simply $Q$ of $T_k$.

\medskip

We define
\begin{align*}
  &I_k^\mathrm{ps}
    := \{ \omega \in (\omega^* + I_\mathrm{prm}) \:|\: 
    \sin(k \omega + \eta^*(0, \omega^*, 0)) < -2e_\mathrm{bd} \}
  \quad (\subset I_k^\mathrm{bd}), \\
  &R_k^\mathrm{ps} := I_\mathrm{prm} \times I_k^\mathrm{ps} \times I_\mathrm{prm}
  \quad (\subset R_k^\mathrm{bd}),
\end{align*}
where $e_\mathrm{bd}$ is the constant in \eqref{e32-Ikbd}.
The subscript `ps' indicates that $E_k$ in \eqref{eq-03-Ek} is positive.

\begin{proposition}[Existence of a non-hyperbolic fixed point]\label{p34-Qk}
  By replacing $\kappa(\delta_\mathrm{dom})$ with a larger one  
  according to Remark~\ref{rem-smaller}, we have the following statements:
  \begin{enumerate}
    \item For any $k \in \mathbb{Z}_{> \kappa(\delta_\mathrm{dom})}$,
      $t$ with $|t| \leq O(1)$,
      $\omega \in I_k^\mathrm{bd}$,
      and $\rho \in I_\mathrm{prm}$,
      by restricting $\mu = \mu_k(t, \omega, \rho)$,
      there exists a fixed point $Q = Q_k = Q_k(t, \omega, \rho) = (Z_Q, Y_Q, W_Q)$
      of $T_k$ such that
      \begin{align}
        &Z_Q, \,
        W_Q = O(\lambda^k), \quad
        Y_Q = \frac{E_k}{2d}\lambda^k t, \quad
        \mu_k(t, \omega, \rho) = O(\lambda^k) + O(\gamma^{-k}),
        \label{eq-03-ZQWQmuk} \\
        &\begin{aligned}
          &Z_{Q, \sigma}, \,
          W_{Q, \sigma}
            = O(k \lambda^k), \quad
          Y_{Q, t}
            = O(\lambda^k), \quad
          Y_{Q, \omega}, \,
          Y_{Q, \rho}
            = O(k \lambda^k), \\
          &\mu_{k, t}
            = O(\lambda^k), \quad
          \mu_{k, \omega}, \,
          \mu_{k, \rho}
            = O(k \lambda^k) + O(k \gamma^{-k}), \quad
        \end{aligned} \label{eq-03-ZQWQmuk_ev}
      \end{align}
      for any $\sigma \in \{t, \omega, \rho\}$.
      Moreover,
      $Q$, $Z_Q$, $Y_Q$, $W_Q$, and $\mu_k$ are $C^{r - 2}$
      with respect to $(t, \omega, \rho)$.
    \item For any $k \in \mathbb{Z}_{> \kappa(\delta_\mathrm{dom})} \cap 2\mathbb{Z}$,
      there exist $C^{r - 2}$ maps $t_k^-$,
      $t_k^+: I_k^\mathrm{ps} \to \mathbb{R}$ with
      \begin{align*}
        t_k^\pm(\omega) = O(1), \quad
        t_k^-(\omega) < t_k^+(\omega)
      \end{align*}
      such that for any
      $(t, \omega) \in \{(t, \omega) \:|\: t_k^-(\omega) < t < t_k^+(\omega), \,
      \omega \in I_k^\mathrm{ps} \}$,
      by restricting $\rho = \rho_k(t, \omega)$,
      the above $Q$ will be non-hyperbolic:
      the multipliers of $Q$ are $\nu_1$, $\nu_2$, and $\nu_3$ such that
      \begin{align}
        \nu_1 = \cos\psi + \mathrm{i} \sin\psi, \quad
        \nu_2 = \cos\psi - \mathrm{i} \sin\psi, \quad
        |\nu_3| < 1
        \label{e34-nu123}
      \end{align}
      for some $C^{r - 2}$ map $\psi = \psi(t, \omega) \in (0, \pi)$,
      where the map $(t_k^-(\omega), t_k^+(\omega))
      \ni t \mapsto \psi(t, \omega) \in (0, \pi)$
      is an orientation reversing $C^{r - 2}$ diffeomorphism
      for any fixed $\omega \in I_k^\mathrm{ps}$.
      Moreover, $\rho_k$ is $C^{r - 2}$ with respect to $(t, \omega)$ with
      \begin{gather}
        \rho_k(t, \omega) = O(k^{-1}),
          \label{eq-03-solrho} \\
        \rho_{k, t}(t, \omega)
          = O(k^{-1} \hat{\lambda}^k \gamma^k), \quad
        \rho_{k, \omega}(t, \omega) = O(1).
          \label{eq-03-solrho_ev}
      \end{gather}
    \item If the original $(f, \Gamma)$ holds the expanding condition \textbf{(EC)},
      then there exists
      $I_k^\mathrm{ex} \subset I_k^\mathrm{ps}$ such that
      \begin{align*}
        \rho_k(t, \omega) < 0
      \end{align*}
      for any $k \in \mathbb{Z}_{> \kappa(\delta_\mathrm{dom})} \cap 2\mathbb{Z}$,
      $\omega \in I_k^\mathrm{ex}$,
      and $t \in (t_k^-(\omega), t_k^+(\omega))$.
  \end{enumerate}
\end{proposition}

\begin{proof}[Proof of Proposition~\ref{p34-Qk}]
We divide the proof into three parts,
corresponding to the items in the lemma.

\medskip

\textbf{(1) First item.}\quad
We put
\begin{align*}
  Y_Q = Y_Q(t, \mu, \omega, \rho)
    := \frac{E_k}{2d}\lambda^k t
\end{align*}
for any $t$ with $|t| \leq O(1)$.
The \eqref{eq-03-NF} and \eqref{eq-03-muhat} imply that
the first-return map $T_k$ has a fixed point $(Z_Q, Y_Q, W_Q)$ if
\begin{align}
\begin{aligned}
  &Z_Q = H_k^{(1)}(t, \omega, \rho, Z_Q, W_Q, \mu), \\
  &W_Q = H_k^{(2)}(t, \omega, \rho, Z_Q, W_Q, \mu), \\
  &\mu = H_k^{(3)}(t, \omega, \rho, Z_Q, W_Q, \mu),
\end{aligned} \label{eq-03-Fsystem1}
\end{align}
where
\begin{align*}
  H_k^{(1)} &:= \lambda^k \alpha_1 Z_Q - \frac{E_k^2}{2d}\lambda^kt
    + \lambda^k \beta_1 W_Q
    + h_1(Z_Q, \frac{E_k}{2d}\lambda^kt, W_Q, \varepsilon), \\
  H_k^{(2)} &:= \lambda^k \alpha_3 Z_Q + \lambda^k \beta_3 W_Q
    + h_3(Z_Q, \frac{E_k}{2d}\lambda^kt, W_Q, \varepsilon), \\
  H_k^{(3)} &:= \frac{E_k}{2d}\lambda^k \gamma^{-k}t + \gamma^{-k} y^-
    - \lambda^k(\alpha^* x_1^+ + \beta^* x_2^+) + O(\hat{\lambda}^k) \\
    &\quad - \lambda^k Z_Q - \frac{E_k^2}{4d}\lambda^{2k}t^2
    - \gamma^{-k} h_2(Z_Q, \frac{E_k}{2d}\lambda^kt, W_Q, \varepsilon).
\end{align*}
\debugtext{
  Hint:
  \begin{align*}
    &\frac{E_k}{2d}\lambda^kt
    = \gamma^k \mu - y^- + \lambda^k \gamma^k (\alpha^* x_1^+ + \beta^* x_2^+)
      + \gamma^k O(\hat{\lambda}^k)
      + \lambda^k \gamma^k Z_Q + d \gamma^k \frac{E_k^2}{4d^2} \lambda^{2k} t^2 + h_2.
  \end{align*}
}
Note that the $O(\hat{\lambda}^k)$ in the above equation
is a function of $(\mu, \omega, \rho)$ and its first partial derivatives
with respect to $\mu$ and $(\omega, \rho)$ have estimates of
$O(\gamma^k \hat{\gamma}^{-k})$ and $O(\hat{\hat{\lambda}}^k)$, respectively,
see Proposition~\ref{p32-NF}.
By the definition of $E_k = E_k(\mu, \omega, \rho)$,
we have
\begin{align}
  E_k = O(1), \quad
  E_{k, \mu}, \,
  E_{k, \rho}
    = O(1), \quad
  E_{k, \omega}
    = O(k).
  \label{eq-03-Eko_evs}
\end{align}
Since we have \eqref{eq-03-NF_higher1}--\eqref{eq-03-NF_higher3},
\begin{align*}
  &H_k^{(1)}, H_k^{(2)} = O(\lambda^k), \quad
  H_k^{(3)} = O(\lambda^k) + O(\gamma^{-k}), \quad
  H_{k, X_Q}^{(1)}, \,
  H_{k, X_Q}^{(2)}, \,
  H_{k, Z_Q}^{(3)}
    = O(\lambda^k), \\
  &H_{k, W_Q}^{(3)}
    = O(\hat{\lambda}^k), \quad
  H_{k, \mu}^{(1)}, \,
  H_{k, \mu}^{(2)}
    = O(k \lambda^k), \quad
  H_{k, \mu}^{(3)}
    = O(\gamma^k \hat{\gamma}^{-k})
\end{align*}
for any $X_Q \in \{Z_Q, W_Q\}$.
Therefore, Proposition~\ref{pc3-sys} solves the system of equations
\eqref{eq-03-Fsystem1} with respect to $(Z_Q, W_Q, \mu)$
as functions of $t$, $\omega$, and $\rho$.
In such a way, we obtain the solutions in \eqref{eq-03-ZQWQmuk}.
Since $H_k^{(i)}$, $i \in \{1, 2, 3\}$, are at least $C^{r - 2}$
with respect to $(t, \omega, \rho, Z_Q, W_Q, \mu)$,  
the solutions are also $C^{r - 2}$
with respect to $(t, \omega, \rho)$.
In addition, \eqref{eq-03-NF_higher1}--\eqref{eq-03-NF_higher3} yield
\begin{align*}
  H_{k, t}^{(1)}
    = O(\lambda^k), \quad
  H_{k, t}^{(2)}, \,
  H_{k, t}^{(3)}
    = O(\hat{\lambda}^k), \quad
  H_{k, \sigma'}^{(1)}, \,
  H_{k, \sigma'}^{(2)}
    = O(k \lambda^k), \quad
  H_{k, \sigma'}^{(3)}
    = O(k \lambda^k) + O(k \gamma^{-k})
\end{align*}
for any $\sigma' \in \{\omega, \rho\}$.
Proposition~\ref{pc3-sys} further implies \eqref{eq-03-ZQWQmuk_ev}.
Here, to get the equations of $Y_{Q, \sigma}$, $\sigma \in \{t, \omega, \rho\}$
in \eqref{eq-03-ZQWQmuk_ev},
we used
\begin{align}
  \partial_t E_k(\mu_k, \omega, \rho)
    = O(\lambda^k), \quad
  \partial_\omega E_k(\mu_k, \omega, \rho)
    = O(k), \quad
  \partial_\rho E_k(\mu_k, \omega, \rho)
    = O(1).
  \label{eq-03-Ek_evs}
\end{align}
Note that in the above computation, $\mu = \mu_k(t, \omega, \rho)$  
is substituted and the chain rule is applied.  
That is, for instance, the computation of $\partial_\omega E_k(\mu_k, \omega, \rho)$
is given by  
\begin{align*}
  \partial_\omega E_k(\mu_k, \omega, \rho)
    = E_{k, \mu} \cdot \mu_{k, \omega}
    + E_{k, \omega}
    = O(k).  
\end{align*} 

\medskip

\textbf{(2) Second item.}\quad
Next, we prove the second item of the lemma.
By Proposition~\ref{p33-CF} and \eqref{eq-03-conecu},
in the tangent space at $Q \in \Pi_k'$,
there exists a forward-invariant subspace
$E^{cu} = E^{cu}(Q) \subset \mathcal{C}^{cu}(Q)$.
Any vector $v^{cu}\in E^{cu}$ has the form
\begin{align*}
  v^{cu} = (z, y, S(z, y)),
\end{align*}
where $S$ is a linear map such that $S(z, y) = S_1 z + S_2 y$ with
$S_1 = S_1(t, \omega, \rho)$ and $S_2 = S_2(t, \omega, \rho)$.

Consider $D(T_k)_Q|_{E^{cu}}$ as the linear transformation of $\mathbb{R}^2$
defined by
\begin{align*}
  D(T_k)_Q|_{E^{cu}}(z, y) = (\bar{z}, \bar{y}),
\end{align*}
where
\begin{align*}
  D(T_k)_Q(z, y, S(z, y)) = (\bar{z}, \bar{y}, S(\bar{z}, \bar{y})).
\end{align*}
Differentiating \eqref{eq-03-NF},
we get the formula for $D(T_k)_Q|_{E^{cu}}$:
\begin{align*}
  \bar{z} = A_k z + (-E_k + B_k) y, \quad
  \bar{y} = (\lambda^k \gamma^k + C_k) z + (E_k \lambda^k \gamma^k t + D_k) y,
\end{align*}
where
\begin{align*}
  &A_k = A_k(t, \omega, \rho)
    := \lambda^k \alpha_1 + h_{1, Z} + (\lambda^k \beta_1 + h_{1, W}) S_1, \\
  &B_k = B_k(t, \omega, \rho)
    := h_{1, Y} + (\lambda^k \beta_1 + h_{1, W}) S_2, \\
  &C_k = C_k(t, \omega, \rho)
    := h_{2, Z} + h_{2, W} S_1, \quad
  D_k = D_k(t, \omega, \rho)
    := h_{2, Y} + h_{2, W} S_2,
\end{align*}
and $h_{i, X} = h_{i, X}(Z_Q, Y_Q, W_Q, \mu_k, \omega, \rho)$
for any $i \in \{1, 2\}$ and $X \in \{Z, Y, W\}$.
Although $A_k$, $B_k$, $C_k$, and $D_k$ involve first partial derivatives,  
these are taken with respect to the spatial variables $(Z, Y, W)$,  
so they are $C^{r - 2}$ with respect to $(t, \omega, \rho)$,  
see Proposition~\ref{p32-NF}.  
Let $\nu_1$ and $\nu_2$ be the eigenvalues of $D(T_k)_Q|_{E^{cu}}$.
Then, we have
\begin{align*}
  &\nu_1 + \nu_2
    = E_k \lambda^k \gamma^k t + A_k + D_k, \quad
  \nu_1 \nu_2
    = E_k \lambda^k \gamma^k + R_k,
\end{align*}
where
\begin{align}
  &R_k = R_k(t, \omega, \rho)
    := A_k (E_k \lambda^k \gamma^k t + D_k)
    + E_k C_k - B_k (\lambda^k \gamma^k + C_k).
  \label{eq-03-Rk}
\end{align}
The $R_k$ is also $C^{r - 2}$ with respect to $(t, \omega, \rho)$.

The $\lambda$ is positive,  
but $\gamma$ may be negative due to the assumption,  
so we assume $k \in 2\mathbb{Z}$,  
and we have $\lambda^k \gamma^k = \mathrm{e}^{k \rho}$  
by the definition \eqref{eq-02-rhoD} of $\rho$.  
We further assume $\omega \in I_k^\mathrm{ps}$
and consider to make $\nu_1 \nu_2$ equal to 1:
\begin{align}
  \rho = H_k, \quad
  H_k = H_k(t, \omega, \rho)
    := -k^{-1} \log E_k + k^{-1}\log(1 - R_k).
  \label{eq-03-rhoeq}
\end{align}
In fact,  
\begin{align}
  S_{i} = O(\lambda^k) + O(\gamma^{-k}), \quad
  S_{i, t} = O(\lambda^k) + O(\gamma^{-k}), \quad
  S_{i, \sigma'} = O(1)
  \label{eq-03-SitSis}
\end{align}
hold for any $\sigma' \in \{\omega, \rho\}$ and $i \in \{1, 2\}$.  
This does not give the sharpest estimate,  
but for how to obtain this estimate, see Remark~\ref{r34-Sev}.  
The estimates \eqref{eq-03-NF_higher1}--\eqref{eq-03-NF_higher3},
as well as the estimates above for $S_1$ and $S_2$, imply
\begin{align*}
  A_k, \,
  B_k
    = O(\lambda^k), \quad
  C_k
    = O(\hat{\lambda}^k \gamma^k), \quad
  D_k
    = O(\lambda^k \gamma^k).
\end{align*}
Using the estimates \eqref{eq-03-NF_higher1}--\eqref{eq-03-NF_higher3},
\eqref{eq-03-ZQWQmuk_ev},
and \eqref{eq-03-SitSis},
we obtain
\begin{align}
\begin{aligned}
  &A_{k, t}, \,
  B_{k, t}
    = O(\hat{\lambda}^k), \quad
  A_{k, \sigma'}, \,
  B_{k, \sigma'}
    = O(k \lambda^k), \quad
  C_{k, t}, \,
  D_{k, t}
    = O(\hat{\lambda}^k \gamma^k), \\
  &C_{k, \sigma'}
    = O(\hat{\hat{\lambda}}^k \gamma^k), \quad
  D_{k, \sigma'}
    = O(k \lambda^k \gamma^k)
\end{aligned} \label{e34-ABCDev}
\end{align}
for any $\sigma' \in \{\omega, \rho\}$,
where we used \eqref{eq-03-hatC} to sort $O(\cdot)$ terms.
The above formulas and \eqref{eq-03-Ek_evs} yield
\begin{align*}
  &R_k = O(\hat{\lambda}^k \gamma^k), \quad
  R_{k, t} = O(\hat{\lambda}^k \gamma^k), \quad
  R_{k, \sigma'} = O(\hat{\hat{\lambda}}^k \gamma^k)
\end{align*}
for any $\sigma \in \{t, \omega, \rho\}$.
Thus, we obtain
\begin{align*}
  H_k = O(k^{-1}), \quad
  H_{k, t} = O(k^{-1} \hat{\lambda}^k \gamma^k), \quad
  H_{k, \omega} = O(1), \quad
  H_{k, \rho} = O(k^{-1})
\end{align*}
by using \eqref{eq-03-hatC} to sort $O(\cdot)$ terms.
Hence, Proposition~\ref{pc1-sgl} implies
the solution \eqref{eq-03-solrho} of the equation
with \eqref{eq-03-solrho_ev}.
Since $H_k$ is $C^{r - 2}$ with respect to $(t, \omega, \rho)$,
the solution is also $C^{r - 2}$ with respect to $(t, \omega)$.

The transformation $D(T_k)_Q|_{E^{cu}}$ has two eigenvalues
in the unit circle if and only if $\nu_1 \nu_2 = 1$ and $|\nu_1 + \nu_2| < 2$
(see e.g. \cite[Section~2.3.1]{L2016});
the boundary $\nu_1 + \nu_2 = 2$ corresponds to the multipliers equal to 1,
and $\nu_1 + \nu_2 = -2$ corresponds to they equal to $-1$.
Let
\begin{align*}
  \Sigma_k = \Sigma_k(t, \omega)
    := \nu_1 + \nu_2
\end{align*}
Then, using $\lambda^k \gamma^k = \mathrm{e}^{k \rho_k}$, we get
\begin{align}
  \Sigma_k = E_k \mathrm{e}^{k \rho_k} t + A_k + D_k
    = E_k \mathrm{e}^{k \rho_k} t + O(1)
  \label{eq-03-Tk}
\end{align}
and $\Sigma_k$ is $C^{r - 2}$ with respect to $(t, \omega)$.
Using \eqref{eq-03-ZQWQmuk_ev},
\eqref{eq-03-solrho},
\eqref{eq-03-Eko_evs},
\eqref{eq-03-Ek_evs},
and \eqref{e34-ABCDev}
we obtain
\begin{align*}
  &\partial_t E_k(\mu_k, \omega, \rho_k)
    = O(k^{-1} \hat{\lambda}^k \gamma^k), \quad
  \partial_t \mathrm{e}^{k \rho_k}
    = O(\hat{\lambda}^k \gamma^k), \\
  &\partial_t A_k(t, \omega, \rho_k)
    = O(\hat{\lambda}^k), \quad
  \partial_t D_k(t, \omega, \rho_k)
    = O(\hat{\lambda}^k \gamma^k).
\end{align*}
This implies there exists a constant $C = C(\mathbb{F}) > 0$ such that
\begin{align}
  \Sigma_{k, t}
    = E_k \mathrm{e}^{k \rho_k} + O(\hat{\lambda}^k \gamma^k) \geq C \quad (> 0)
  \label{eq-03-Tkt}
\end{align}
for any $k \in \mathbb{Z}_{> \kappa(\delta_\mathrm{dom})} \cap 2\mathbb{Z}$
by replacing $\kappa(\delta_\mathrm{dom})$ with a larger one  
according to Remark~\ref{rem-smaller}.
By using the intermediate value theorem,
\eqref{eq-03-Tk} and \eqref{eq-03-Tkt} yield,
for any $\omega \in I_k^\mathrm{ps}$,
there exist unique $t_k^- = t_k^-(\omega)$ and $t_k^+ = t_k^+(\omega)$
with $t_k^-(\omega) < t_k^+(\omega)$ such that
\begin{align*}
  \Sigma_k(t_k^-(\omega), \omega) = -2, \quad
  \Sigma_k(t_k^+(\omega), \omega) = 2, \quad
  t_k^\pm = O(1).
\end{align*}
By \eqref{eq-03-Tkt}, the implicit function theorem says $t_k^\pm$ are $C^{r - 2}$
with respect to $\omega$.
For any $(t, \omega) \in \{(t, \omega) \:|\: t_k^-(\omega) < t < t_k^+(\omega), \,
\omega \in I_k^\mathrm{ps} \}$,
we can write
\begin{align*}
  \nu_1 = \cos\psi + \mathrm{i} \sin\psi, \quad
  \nu_2 = \cos\psi - \mathrm{i} \sin\psi,
\end{align*}
where $\psi = \psi(t, \omega) := \arccos(\frac{\Sigma_k(t, \omega)}{2})
\in [0, \pi]$.
The $\psi$ is $C^{r - 2}$ with respect to $(t, \omega)$
and the map $(t_k^-, t_k^+) \ni t \mapsto \psi(t, \omega) \in (0, \pi)$
is an orientation reversing $C^{r - 2}$ diffeomorphism
for any $\omega \in I_k^\mathrm{ps}$
since we have \eqref{eq-03-Tkt} and the restriction
$(-1, 1) \ni x \mapsto \arccos(x) \in (0, \pi)$ is
an orientation reversing diffeomorphism.
The remaining eigenvalue of $D(T_k)_Q|_{E^{cu}}$ is
inside the unit circle by Remark~\ref{rem-WContract}.

\medskip

\textbf{(3) Third item.}\quad
Recall the coefficients in \eqref{eq-03-coepdep} and \eqref{eq-03-coedep}
and $c(\varepsilon)$ in \eqref{eq-03-defc}.
For the quantity $\mathcal{E}$ in \eqref{eq-02-mathcalE},
we have
\begin{align*}
  \mathcal{E} =
    \sqrt{(b_1'(\varepsilon^*))^2 + (b_2'(\varepsilon^*))^2}
    \sqrt{(c_1'(\varepsilon^*))^2 + (c_2'(\varepsilon^*))^2}
    = b(\varepsilon^*) \sqrt{(c_1(\varepsilon^*))^2 + (c_2(\varepsilon^*))^2}
    = b(\varepsilon^*) c(\varepsilon^*),
\end{align*}
where $\varepsilon^* = (0, \omega^*, 0)$.
Since $(f, \Gamma)$ holds \textbf{(EC)},
we have
\begin{align*}
  0 < \mathcal{E} - 1, \quad
  0 < \frac{\mathcal{E} - 1}{6\mathcal{E}}
    < \frac{\mathcal{E} - 1}{3\mathcal{E}} < \frac{1}{3}.
\end{align*}
We put
\begin{align*}
  \delta = \delta(\mathbb{F})
    := \frac{\mathcal{E} - 1}{3} > 0, \quad
  \delta' = \delta'(\mathbb{F})
    := \frac{\mathcal{E} - 1}{6\mathcal{E}} > 0.
\end{align*}
Thus, we may suppose
\begin{align}
  |b(\varepsilon) c(\varepsilon)
    - \mathcal{E}| < \delta, \quad
  \delta' < \frac{\mathcal{E} - 1}{3 b(\varepsilon) c(\varepsilon)}
  \label{eq-03-bc}
\end{align}
for any $\varepsilon \in R_\mathrm{prm}$
by replacing $\delta_\mathrm{prm} > 0$ with a smaller one
according to Remark~\ref{rem-smaller}.
By \eqref{e52-alp2},
and \eqref{eq-03-Ek}
\begin{align}
  E_k(\varepsilon) = -b(\varepsilon) c(\varepsilon)
    \sin(k\omega + \eta^*(\varepsilon)).
  \label{eq-03-Ekv}
\end{align}

On the other hand,
we define
\begin{align}
  &I_k^\mathrm{ex}
    := \{ \omega \in (\omega^* + I_\mathrm{prm}) \:|\: 
    \sin(k \omega + \eta^*(0, \omega^*, 0)) + 1 < \delta'/2 \}
  \quad (\subset I_k^\mathrm{ps}),
  \label{e34-Ikex} \\
  &R_k^\mathrm{ex} := I_\mathrm{prm} \times I_k^\mathrm{ex} \times I_\mathrm{prm}
    \quad (\subset R_k^\mathrm{ps}).
  \label{e34-Rkex}
\end{align}
By replacing $\delta_\mathrm{prm} > 0$ with a smaller one
according to Remark~\ref{rem-smaller}, we have
\begin{align}
  \sin(k \omega + \eta^*(\mu, \omega, \rho)) + 1 < \delta'
  \label{e34-sin}
\end{align}
for any $\varepsilon \in R_k^\mathrm{ex}$.
We take $\omega \in I_k^\mathrm{ex}$
and $t \in (t_k^-(\omega), t_k^+(\omega))$
and fix $\rho = \rho_k(t, \omega)$,
$\mu = \mu_k(t, \omega, \rho_k)$.
To simplify the notation,
let $\varepsilon_k := (\mu_k, \omega, \rho_k)$.
Using \eqref{eq-03-bc},
\eqref{eq-03-Ekv},
and \eqref{e34-sin}, we have
\begin{align*}
  |E_k(\varepsilon_k) - \mathcal{E}|
  &\leq |E_k(\varepsilon_k) - b(\varepsilon_k) c(\varepsilon_k)|
    + |b(\varepsilon_k) c(\varepsilon_k)
    - \mathcal{E}| \\
  &< |\sin(k \omega + \eta^*(\varepsilon_k)) + 1|
    |b(\varepsilon_k) c(\varepsilon_k)| + \delta \\
  &< \delta' |b(\varepsilon_k) c(\varepsilon_k)| + \delta \\
  &< 2\delta
\end{align*}
for any $k \in \mathbb{Z}_{> \kappa(\delta_\mathrm{dom})} \cap 2\mathbb{Z}$.
Hence,
\begin{align}
  E_k(\varepsilon_k) > \mathcal{E} - 2\delta = \frac{\mathcal{E} + 2}{3}.
  \label{eq-03-Ekj}
\end{align}
In contrast, recall the $R_k$ in \eqref{eq-03-Rk}.
Since $R_k = O(\hat{\lambda}^k \gamma^k)$,
\begin{align*}
  |R_k(t, \omega, \rho_k)| < \delta
\end{align*}
for any $k \in \mathbb{Z}_{> \kappa(\delta_\mathrm{dom})} \cap 2\mathbb{Z}$
by replacing $\kappa(\delta_\mathrm{dom})$ with a larger one  
according to Remark~\ref{rem-smaller}.
Thus,
\begin{align}
  1 - R_k(t, \omega, \rho_k) < 1 + \delta = \frac{\mathcal{E} + 2}{3}.
  \label{eq-03-Rkj}
\end{align}
The \eqref{eq-03-Ekj} and \eqref{eq-03-Rkj} imply
\begin{align*}
  1 - R_k(t, \omega, \rho_k) < E_k(\varepsilon_k).
\end{align*}
Hence, \eqref{eq-03-rhoeq} yields
\begin{align*}
  \rho_k(t, \omega) = H_k(t, \omega, \rho_k)
  = \frac{1}k \log
    \frac{1 - R_k(t, \omega, \rho_k)}{E_k(\varepsilon_k)} < 0.
\end{align*}
The desired statement is proved.  
\end{proof}

\begin{remark} \label{r34-Sev}
We explain how \eqref{eq-03-SitSis} is obtained.
First, the leftmost estimate in \eqref{eq-03-SitSis}  
follows from the definition of the cone field \eqref{eq-03-conecu}  
and the fact that $Y_Q = O(\lambda^k)$.  
From \eqref{eq-03-DTkM}, the $D(T_k)_Q: (z, y, w) \mapsto (\bar{z}, \bar{y}, \bar{w})$  
is expressed as  
\begin{align*}
  \bar{z} &= A_{11} z + A_{12} y + A_{13} w,
  & A_{11} &= O(\lambda^k),\quad
    A_{12} = -E_k + O(\lambda^k),\quad
    A_{13} = O(\lambda^k), \\
  \bar{y} &= A_{21} z + A_{22} y + A_{23} w,
  & A_{21} &= O(\lambda^k \gamma^k),\quad
    A_{22} = O(\lambda^k \gamma^k),\quad
    A_{23} = O(\hat{\lambda}^k \gamma^k), \\
  \bar{w} &= A_{31} z + A_{32} y + A_{33} w,
  & A_{31} &= O(\lambda^k),\quad
    A_{32} = O(\lambda^k),\quad
    A_{33} = O(\lambda^k).
\end{align*}
If a plane $w = S_1 z + S_2 y$ with \eqref{eq-03-SitSis} is mapped by $D(T_k)_Q$  
to a new plane $\bar{w} = \bar{S}_1 \bar{z} + \bar{S}_2 \bar{y}$,  
then they satisfy the following equation:  
\begin{align*}
  M_\mathrm{coe} \bar{S} = b_\mathrm{cst}, \quad
  M_\mathrm{coe} :=
  \begin{pmatrix}
    A_{11} + A_{13} S_1 & A_{21} + A_{23} S_1 \\
    A_{12} + A_{13} S_2 & A_{22} + A_{23} S_2
  \end{pmatrix}, \quad
  \bar{S} :=
  \begin{pmatrix}
    \bar{S}_1 \\
    \bar{S}_2
  \end{pmatrix}, \quad
  b_\mathrm{cst} :=
  \begin{pmatrix}
    A_{31} + A_{33} S_1 \\
    A_{32} + A_{33} S_2
  \end{pmatrix}.
\end{align*}
Differentiating both sides of this equation with respect to  
$\sigma \in \{t, \omega, \rho\}$,  
we obtain  
\begin{align*}
  \partial_\sigma \bar{S}
    = M_\mathrm{coe}^{-1}(\partial_\sigma b_\mathrm{cst}
    - (\partial_\sigma M_\mathrm{coe}) \bar{S}).
\end{align*}
From \eqref{eq-03-NF}, \eqref{eq-03-NF_higher2}, and \eqref{eq-03-NF_higher3},  
we obtain the rough estimates  
\begin{align*}
  \partial_t b_\mathrm{cst} - (\partial_t M_\mathrm{coe}) \bar{S}
    = (O(\lambda^k), O(\lambda^k))^\mathsf{T}, \quad
  \partial_{\sigma'} b_\mathrm{cst} - (\partial_{\sigma'} M_\mathrm{coe}) \bar{S}
    = (O(\hat{\lambda}^k \gamma^k), O(\hat{\lambda}^k \gamma^k))^\mathsf{T}
\end{align*}
for any $\sigma' \in \{\omega, \rho\}$.  
Since the inverse of $M_\mathrm{coe}$ is given by  
\begin{align*}
  M_\mathrm{coe} =
  \begin{pmatrix}
    O(1) & O(1) \\
    O(\lambda^{-k} \gamma^{-k}) & O(\gamma^{-k})
  \end{pmatrix},
\end{align*}
the new plane $\bar{w} = \bar{S}_1 \bar{z} + \bar{S}_2 \bar{y}$  
again satisfies the estimates in \eqref{eq-03-SitSis}.  
Since the plane $w = S_1 z + S_2 y$ in the proof  
is the limit of such an iteration,  
\eqref{eq-03-SitSis} holds.  
\end{remark}

\ifdefined\isMaster
\else
  \end{document}
\fi
\ifdefined\isMaster
\else
  \documentclass[11pt,a4paper]{article}

\usepackage[left=2cm,right=2cm,top=2cm,bottom=2cm]{geometry}
\usepackage{graphicx}
\usepackage{multirow}
\usepackage{amsmath,amssymb,amsfonts}
\usepackage{amsthm}
\usepackage{mathrsfs}
\usepackage[title]{appendix}
\usepackage{xcolor}
\usepackage{textcomp}
\usepackage{manyfoot}
\usepackage{booktabs}
\usepackage{algorithm}
\usepackage{algorithmicx}
\usepackage{algpseudocode}
\usepackage{listings}
\ifdefined\debug
  \usepackage[inline]{showlabels}
  \showlabels{ref}
  \showlabels{cite}
  \usepackage{fancyhdr}
\fi
\usepackage{bm}
\usepackage{hyperref}

\theoremstyle{plain}
\newtheorem{theorem}{Theorem}[section]
\newtheorem{proposition}[theorem]{Proposition}%
\newtheorem{lemma}[theorem]{Lemma}
\newtheorem{corollary}[theorem]{Corollary}
\newtheorem{claim}[theorem]{Claim}
\theoremstyle{definition}
\newtheorem{definition}[theorem]{Definition}
\newtheorem{remark}[theorem]{Remark}
\newtheorem{example}[theorem]{Example}
\newtheorem{question}[theorem]{Question}
\newtheorem*{notation}{Notation}
\theoremstyle{plain}
\newtheorem{maintheorem}{Theorem}
\renewcommand{\themaintheorem}{\Alph{maintheorem}}
\newcommand{\wiptext}[1]{\textcolor{blue}{#1}}
\ifdefined\debug
  \newcommand{\debugtext}[1]{\textcolor{red}{#1}}
  \pagestyle{fancy} 
  \fancyhead[L]{}
  \fancyhead[C]{}
  \fancyhead[R]{\textcolor{red}{[DEBUG] Debug mode is activated.}}
  \renewcommand{\headrulewidth}{0pt}
\else
  \newcommand{\debugtext}[1]{\ignorespaces}
\fi

\counterwithin{figure}{section}
\counterwithin{table}{section}

\raggedbottom

\numberwithin{equation}{section}

  \begin{document}
\fi

\section{Verifying the periodic point is a generic Hopf point}\label{s4-Hopf}
The periodic point $Q$,
or precisely $Q = Q_k = Q_k(t, \omega, \rho_k)$
in Proposition~\ref{p34-Qk} is non-hyperbolic
since it has a complex multiplier with norm one and it looks like that
$Q$ is a generic Hopf point.
To verify it,
we need to calculate the Lyapunov coefficient $\mathrm{LC}(Q)$
of $Q$ that determines whether it is attracting or repelling
on its two-dimensional local center manifold $W^c_\mathrm{loc}(Q)$
by seeing the higher order terms of $T_k|_{W^c_\mathrm{loc}(Q)}$.
In this section, we calculate it accurately.
In Section~\ref{s41-LC},
we first give the formula of the Lyapunov coefficient
in general settings.
In Section~\ref{s42-LCQk},
we give the precise formula of the Lyapunov coefficient of $Q$.
In Section~\ref{s43-rep},
we verify that
$Q$ is a repeller on its two-dimensional local center manifold
for appropriate value of $t$.

\subsection{Formula of the Lyapunov coefficient}\label{s41-LC}
In this section,
we give the formula for the Lyapunov coefficient for general systems.
The following discussion is based on
\cite[Section~7, 8]{RT1971},
\cite[Chapter~III]{I1979},
\cite[Section~6, 6A]{MM2012},
and \cite[Section~2.8]{D2018}.
See these references for details.

\medskip

Let $T: \mathbb{C} \ni z \mapsto \tilde{z} \in \mathbb{C}$
be a $C^r$, $r \geq 4$, map having the expansion
\begin{align}
  \tilde{z} = \nu z
    + \sum_{2 \leq p + q \leq 3} \tilde{z}^{(pq)} z^p \bar{z}^q
    + O(|z|^4), \quad
  \nu = \cos \psi + \mathrm{i} \sin \psi,
  \label{eq-04-z2zt}
\end{align}
where $O(|z|^4)$ is a term of fourth order or higher
and $\psi \in (0, \pi)$.
Here,
we always assume $p$ and $q$ are non-negative integers.
Recall $\Psi_\mathrm{reg}$ in \eqref{e12-Preg}.
For any $\psi \in \Psi_\mathrm{reg}$,
putting new coordinate $w \in \mathbb{C}$ as
\begin{align}
  w = z + \sum_{p + q = 2} \frac{\tilde{z}^{(pq)}}{\nu - \nu^p \bar{\nu}^q}
    z^p \bar{z}^q,
  \label{eq-04-z2w}
\end{align}
we have the new expression of $T: w \mapsto \tilde{w}$ as
\begin{align*}
  \tilde{w} = \nu w
    + \sum_{p + q = 3} \tilde{w}^{(pq)} w^p \bar{w}^q
    + O(|w|^4),
\end{align*}
where $O(|w|^4)$ is a term of fourth order or higher.
Note that the homogeneous quadratic terms completely disappear
in the above equation.
In fact,
the coefficient
\begin{align*}
  \alpha := \tilde{w}^{(21)}
\end{align*}
is the same as $\alpha$ in the normal form \eqref{eq-01-Cform}.
Recall that the Lyapunov coefficient is defined by
$\mathrm{LC}(0) = -\Re(\bar{\nu} \alpha)$.

\begin{proposition}[Formula of $\mathrm{LC}$]\label{p41-FmlaLC}
  The Lyapunov coefficient is given by
  \begin{align}
    \mathrm{LC}(0; w)
      = \Re&\left( -\tilde{z}^{(21)} \bar{\nu}
      + |\tilde{z}^{(02)}|^2 \frac{-4 + 2\bar{\nu}^3}{-2 + \nu^3 + \bar{\nu}^3}
      + |\tilde{z}^{(11)}|^2
        \frac{-2\bar{\nu} + \bar{\nu}^2}{(-1 + \bar{\nu})^2}
      + \tilde{z}^{(11)} \tilde{z}^{(20)}
        \frac{2 - 6\bar{\nu} + \bar{\nu}^2}{(-1 + \nu)^2}
    \right)
    \label{e41-LC}
  \end{align}
  for some coordinate $w$ that gives the normal form \eqref{eq-01-Cform}.
\end{proposition}

\begin{proof}[Proof of Proposition~\ref{p41-FmlaLC}]
By \eqref{eq-04-z2w},
we get
\begin{align*}
  \tilde{w} = \tilde{z}
    + \sum_{p + q = 2} \frac{\tilde{z}^{(pq)}}{\nu - \nu^p \bar{\nu}^q}
    \tilde{z}^p \bar{\tilde{z}}^q.
\end{align*}
Substituting \eqref{eq-04-z2zt} into the above equation,
we have
\begin{align}
\begin{aligned}
  \tilde{w} &= \nu z + \frac{\tilde{z}^{(20)}}{1 - \nu} z^2
    + \frac{\tilde{z}^{(11)}}{1 - \bar{\nu}} z \bar{z}
    + \frac{\tilde{z}^{(02)} \nu}{\nu - \bar{\nu}^2} \bar{z}^2
    + \left( \tilde{z}^{(30)}
      - \frac{2 (\tilde{z}^{(20)})^2}{-1 + \nu}
      + \frac{\tilde{z}^{(11)} \overline{\tilde{z}^{(02)}}}{1 - \bar{\nu}}
      \right) z^3 \\
    &\quad + \left( \tilde{z}^{(21)}
      + \tilde{z}^{(11)} \tilde{z}^{(20)} \frac{1 - 2\nu}{\nu^2 - \nu}
      + \frac{|\tilde{z}^{(11)}|^2}{1 - \bar{\nu}}
      + \frac{2 |\tilde{z}^{(02)}|^2}{\nu^2 - \bar{\nu}}
      \right) z^2 \bar{z} \\
    &\quad + \left( \tilde{z}^{(12)}
      + \frac{2 \tilde{z}^{(20)} \tilde{z}^{(02)}}{1 - \nu}
      + \frac{\tilde{z}^{(11)} \overline{\tilde{z}^{(20)}}}{1 - \bar{\nu}}
      + \frac{(\tilde{z}^{(11)})^2}{\nu^2 - \nu}
      + \frac{2 \tilde{z}^{(02)} \overline{\tilde{z}^{(11)}}}{\nu^2 - \bar{\nu}}
    \right) z \bar{z}^2 \\
    & \quad + \left( \tilde{z}^{(03)}
      + \frac{\tilde{z}^{(11)} \tilde{z}^{(02)}}{\nu^2 - \nu}
      + \frac{2 \tilde{z}^{(02)} \overline{\tilde{z}^{(20)}}}{\nu^2 - \bar{\nu}}
    \right) \bar{z}^3
    + O(|z|^4),
\end{aligned} \label{eq-04-z2wt}
\end{align}
where $O(|z|^4)$ is a term of fourth order or higher.
The \eqref{eq-04-z2w} implies the inverse transformation
\begin{align}
  z = w - \sum_{p + q = 2} \frac{\tilde{z}^{(pq)}}{\nu - \nu^p \bar{\nu}^q}
    w^p \bar{w}^q + O(|w|^3),
  \label{e41-w2z}
\end{align}
where $O(|w|^3)$ is a term of third order or higher.
Substituting the above equation into \eqref{eq-04-z2wt},
we obtain
\begin{align}
\begin{aligned}
  \alpha = \tilde{w}^{(21)}
  = \tilde{z}^{(21)}
  + |\tilde{z}^{(02)}|^2 \frac{4\nu - 2\bar{\nu}^2}{-2 + \nu^3 + \bar{\nu}^3}
  + |\tilde{z}^{(11)}|^2 \frac{2 - \bar{\nu}}{(-1 + \bar{\nu})^2}
  - \tilde{z}^{(11)} \tilde{z}^{(20)}
    \frac{-6 + 2\nu + \bar{\nu}}{(-1 + \nu)^2},
\end{aligned} \label{eq-04-FmlaAlp}
\end{align}
which yields the desired formula \eqref{e41-LC}.
\end{proof}

\begin{remark}
  The correct formula is \eqref{eq-04-FmlaAlp},
  although \cite[p.30]{I1979} gives a different formula
  for $\alpha$ than \eqref{eq-04-FmlaAlp}.
  The formula given there is the coefficient
  of $z^2 \bar{z}$ in \eqref{eq-04-z2wt}.
\end{remark}

\subsection{Lyapunov coefficient of the periodic point}
\label{s42-LCQk}
Recall that there exists a non-hyperbolic fixed point $Q$,
or more precisely $Q = Q_k = Q_k(t, \omega, \rho_k) = (Z_Q, Y_Q, W_Q)$,
of $T_k$ for any $(t, \omega) \in \{t_k^-(\omega) < t < t_k^+(\omega), \,
\omega \in I_k^\mathrm{ps}\}$
and for any $k \in \mathbb{Z}_{> \kappa(\delta_\mathrm{dom})} \cap 2\mathbb{Z}$
by Proposition~\ref{p34-Qk}.
In this section, we compute the Lyapunov coefficient of the periodic point $Q_k$.

\medskip

Recall the multipliers $\nu_1$, $\nu_2$,
and the argument $\psi = \psi(t, \omega)$ in \eqref{e34-nu123}.
Also, recall $\Psi_\mathrm{reg}$ in \eqref{e12-Preg}.
For any $(t, \omega) \in \{t_k^-(\omega) < t < t_k^+(\omega), \,
\omega \in I_k^\mathrm{ps}\}$
with $\psi(t, \omega) \in \Psi_\mathrm{reg}$, the Lyapunov coefficient
$\mathrm{LC}(Q_k)$ of $Q_k$ is defined, see Section~\ref{s11-plan}.
Our goal in this section is to prove the following proposition.  

\begin{proposition}[Lyapunov coeffcient of $Q_k$] \label{p42-LC}
  The Lyapunov coefficient of $Q_k$ is given by
  \begin{align}
    \mathrm{LC}(Q_k; w)
      = \mathcal{L}(\psi) + O(\lambda^k), \quad
    \mathcal{L}(\psi)
      := \frac{4 \cos \psi (1 + \cos \psi)}{(-1 + \cos \psi)(1 + 2\cos \psi)^2}
    \label{e42-LCQ}
  \end{align}
  for some coordinate $w$ that gives the normal form \eqref{eq-01-Cform}.
\end{proposition}

\begin{remark} \label{r42-psibd}
  We fix the constant $\psi_\mathrm{bd} = \pi/20$.  
  In fact, any number in the interval $\psi_\mathrm{bd} \in (0, \pi/2)$ would suffice,  
  but throughout this paper, we will always use the above value.  
  Under this setting, as long as $\psi \in (0, \pi/2 - \psi_\mathrm{bd}]$,  
  $\mathcal{L}(\psi)$ has a negative maximum.
  Therefore,  
  by replacing $\kappa(\delta_\mathrm{dom})$ with the larger one  
  according to Remark~\ref{rem-smaller},  
  there exists a constant $C = C(\psi_\mathrm{bd}) > 0$ such that  
  $\mathcal{L}(\psi), \, \mathrm{LC}(Q_k) \leq -C$ for any  
  $k \in \mathbb{Z}_{> \kappa(\delta_\mathrm{dom})} \cap 2\mathbb{Z}$
  and $(t, \omega) \in \{t_k^-(\omega) < t < t_k^+(\omega), \,  
  \omega \in I_k^\mathrm{ps}\}$ with  
  $\psi(t, \omega) \in (0, \pi/2 - \psi_\mathrm{bd}] \cap \Psi_\mathrm{reg}
  = (0, \pi/2 - \psi_\mathrm{bd}]$.
\end{remark}

The proof of the above lemma is carried out through three subsections.  
\begin{itemize}
  \item In Section~\ref{s42-ZY}, using the coordinates $(Z, Y)$  
    from Proposition~\ref{p32-NF}, we describe the restriction of the global map  
    to the local center manifold of $Q$ as $(Z, Y) \mapsto (\bar{Z}, \bar{Y})$.  
    We also estimate the coefficients in $\bar{Z}$ and $\bar{Y}$.  
  \item In Section~\ref{s42-uv}, we introduce new coordinates $(u, v)$  
    instead of $(Z, Y)$.  
    These coordinates define the complex variable $z = u + \mathrm{i} v$,  
    so that the restriction of the global map to the local center manifold of $Q$  
    is expressed as
    $z = u + \mathrm{i} v \mapsto \bar{u} + \mathrm{i} \bar{v} = \tilde{z}$  
    in the form of \eqref{eq-04-z2zt}.  
    We also estimate the coefficients in $\bar{u}$ and $\bar{v}$.  
  \item Finally, in Section~\ref{s42-prf}, we give the proof  
    of Proposition~\ref{p42-LC}.  
    Since we have already given the formula for the Lyapunov coefficient  
    in Proposition~\ref{p41-FmlaLC}, the proof is completed by applying it.  
\end{itemize}  

\subsubsection{Estimate for the original coordinates} \label{s42-ZY}
By Proposition~\ref{p34-Qk} and its proof,
the center manifold theorem \cite{K1967} and \cite[Section~5A]{HPS1970} says that
there is the two-dimensional local center manifold $W^c_\mathrm{loc}(Q)$
that is transverse to the $W$-direction at $Q$.
Let us move the origin to $Q$ by applying
\begin{align*}
  Z^{new} = Z - Z_Q, \quad Y^{new} = Y - Y_Q, \quad W^{new} = W - W_Q
\end{align*}
which allows us to rewrite \eqref{eq-03-NF} as
\begin{align}
  \begin{aligned}
    \bar{Z} &=
      \lambda^k \alpha_1 Z - E_k Y + \lambda^k \beta_1 W
      + h_1(M + Q, \varepsilon) - h_1(Q, \varepsilon), \\
    \bar{Y} &=
      \mathrm{e}^{k \rho_k} Z + E_k \mathrm{e}^{k \rho_k} t Y + d \gamma^k Y^2
      + h_2(M + Q, \varepsilon) - h_2(Q, \varepsilon), \\
    \bar{W} &=
      \lambda^k \alpha_3 Z + \lambda^k \beta_3 W
      + h_3(M + Q, \varepsilon) - h_3(Q, \varepsilon),
  \end{aligned}\label{e42-bZYW}
\end{align}
where $M = (Z, Y, W)$
and the label `new' was dropped.
Here, note that $\lambda^k \gamma^k = \mathrm{e}^{k \rho_k}$
and $Y_Q = (2d)^{-1} E_k \lambda^k t$,
see Proposition~\ref{p34-Qk} and its proof.
By using the above new coordinates, $W^c_\mathrm{loc}(Q)$ has the form
\begin{align*}
  W = w^c(Z, Y)
\end{align*}
for some at least $C^5$ map $w^c$ from a small open two-dimensional disk centered
$(0, 0)$ to a small open one-dimensional disk centered at 0 with $w^c(0, 0) = 0$.
The smoothness is at least $C^5$ because we are currently assuming $r \geq 5$.
By \eqref{e42-bZYW}, the system $T_k|_{W^c(Q)}$ is given by
\begin{align}
  \bar{Z} &=
    \lambda^k \alpha_1 Z -E_k Y + \lambda^k \beta_1 w^c(Z, Y)
    + h_1(M(Z, Y) + Q, \varepsilon) - h_1(Q, \varepsilon), \\
  \bar{Y} &=
    \mathrm{e}^{k \rho_k} Z + E_k \mathrm{e}^{k \rho_k} t Y + d \gamma^k Y^2
    + h_2(M(Z, Y) + Q, \varepsilon) - h_2(Q, \varepsilon),
  \label{e42-bZbYo}
\end{align}
where $M(Z, Y) = (Z, Y, w^c(Z, Y))$ and $\varepsilon = (\mu_k, \omega, \rho_k)$.
We write the Taylor expansion of this system at $(Z, Y) = 0$ as
\begin{align}
  \bar{Z} = \sum_{1 \leq p + q \leq 3} \bar{Z}^{(pq)} Z^p Y^q
    + O(\|(Z, Y)\|^4), \quad
  \bar{Y} = \sum_{1 \leq p + q \leq 3} \bar{Y}^{(pq)} Z^p Y^q
    + O(\|(Z, Y)\|^4),
  \label{e42-bZbY}
\end{align}
where $O(\cdot)$ are terms of fourth order or higher and $p$, $q \geq 0$.
The following holds for the coefficients $\bar{Z}^{(pq)}$ and $\bar{Y}^{(pq)}$  
in \eqref{e42-bZbY}.  

\begin{lemma}[Estimate of $\bar{Z}^{(pq)}$ and $\bar{Y}^{(pq)}$] 
\label{l42-evbZbY}
  We have
  \begin{align*}
    &\bar{Z}^{(pq)} = 
    \begin{cases}
      O(1) & \text{if} \quad (p, q) \in \{(0, 1), \, (0, 2), \, (0, 3)\}, \\
      O(\lambda^k) & \text{otherwise},
    \end{cases} \\
    &\bar{Y}^{(pq)} = 
    \begin{cases}
      O(\gamma^k) & \text{if} \quad (p, q) \in \{(0, 2), \, (1, 2), \, (0, 3)\}, \\
      O(1) & \text{otherwise}
    \end{cases}
  \end{align*}
  for any \( p, \, q \geq 0 \) with \( 1 \leq p + q \leq 3 \).
\end{lemma}

\begin{proof}
In rough estimates,
we have
\begin{align}
  w^c_X(0, 0), \,
  w^c_{X X'}(0, 0), \,
  w^c_{X X' X''}(0, 0)
    = O(1)
  \label{e42-wc}
\end{align}
for any $X$, $X'$, $X'' \in \{Z, Y\}$
since Proposition~\ref{p32-NF} yields that
the $C^2$ norm of $T_k$ in a small neighborhood of $Q_k$ is bounded
with respect to $k$.
By \eqref{eq-03-NF_higher2} and
\eqref{eq-03-NF_higher3}, we can estimate the partial derivatives of $h_1$ and $h_2$
at $(Z, Y, W) = Q$ as
\begin{align*}
  &h_{1, X}, \,
  h_{1, Y}, \,
  h_{1, X Y}, \,
  h_{1, X X'}, \,
  h_{1, X X' X''}, \,
  h_{1, X X' Y}, \,
  h_{1, X Y Y}, \,
    = O(\lambda^k), \quad
  h_{1, Y Y}, \,
  h_{1, Y Y Y}
    = O(1), \\
  &h_{2, X}, \,
  h_{2, Y}, \,
  h_{2, X X'}, \,
  h_{2, X Y}, \,
  h_{2, X X' X''}, \,
  h_{2, X X' Y}, \,
    = O(1), \quad
  h_{2, Y Y}, \,
  h_{2, X Y Y}, \,
  h_{2, Y Y Y}
    = O(\gamma^k)
\end{align*}
for any $X$, $X'$, $X'' \in \{Z, W\}$.
Thus, when $p + q = 1$,
using \eqref{e42-bZbYo} and \eqref{e42-wc},
we have
\begin{align}
\begin{aligned}
  &\bar{Z}^{(10)}
    = (\partial_Z \bar{Z})|_{(Z, Y) = 0}
    = \lambda^k \alpha_1 + \lambda^k \beta_1 w^c_Z + h_{1, Z} + h_{1, W} w^c_Z
    = O(\lambda^k), \\
  &\bar{Z}^{(01)}
    = (\partial_Y \bar{Z})|_{(Z, Y) = 0}
    = -E_k + \lambda^k \beta_1 w^c_Y + h_{1, Y} + h_{1, W} w^c_Y
    = O(1), \\
  &\bar{Y}^{(10)}
    = (\partial_Z \bar{Y})|_{(Z, Y) = 0}
    = \mathrm{e}^{k \rho_k} + h_{2, Z} + h_{2, W} w^c_Z
    = O(1), \\
  &\bar{Y}^{(01)}
    = (\partial_Y \bar{Y})|_{(Z, Y) = 0}
    = E_k \mathrm{e}^{k \rho_k} t + h_{2, Y} + h_{2, W} w^c_Y
    = O(1).
\end{aligned} \label{e42-bZbY1}
\end{align}
Next, when $p + q = 2$,
in a similar manner,
\begin{align}
\begin{aligned}
  &2 \bar{Z}^{(20)}
    = \lambda^k \beta_1 w^c_{Z Z} + h_{1, Z Z} + h_{1, Z W} w^c_Z
    + (h_{1, W Z} + h_{1, W W} w^c_Z) w^c_Z + h_{1, W} w^c_{Z Z}
    = O(\lambda^k), \\
  &2 \bar{Z}^{(11)}
    = \lambda^k \beta_1 w^c_{Z Y} + h_{1, Z Y} + h_{1, Z W} w^c_Y
    + (h_{1, W Y} + h_{1, W W} w^c_Y) w^c_Z + h_{1, W} w^c_{Z Y}
    = O(\lambda^k), \\
  &2 \bar{Z}^{(02)}
    = \lambda^k \beta_1 w^c_{Y Y} + h_{1, Y Y} + h_{1, Y W} w^c_Y
    + (h_{1, W Y} + h_{1, W W} w^c_Y) w^c_Y + h_{1, W} w^c_{Y Y}
    = O(1), \\
  &2 \bar{Y}^{(20)}
    = h_{2, Z Z} + h_{2, Z W} w^c_Z
    + (h_{2, W Z} + h_{2, W W} w^c_Z) w^c_Z + h_{2, W} w^c_{Z Z}
    = O(1), \\
  &2 \bar{Y}^{(11)}
    = h_{2, Z Y} + h_{2, Z W} w^c_Y
    + (h_{2, W Y} + h_{2, W W} w^c_Y) w^c_Z + h_{2, W} w^c_{Z Y}
    = O(1), \\
  &2 \bar{Y}^{(02)}
    = d \gamma^k + h_{2, Y Y} + h_{2, Y W} w^c_Y
    + (h_{2, W Y} + h_{2, W W} w^c_Y) w^c_Y + h_{2, W} w^c_{Y Y}
    = O(\gamma^k).
\end{aligned} \label{e42-bZbY2}
\end{align}
Finally, when $p + q = 3$, we can calculate $\bar{Z}^{(pq)}$ and $\bar{Y}^{(pq)}$  
in the same way, and one can find that  
\begin{itemize}
  \item the $\bar{Z}^{(30)}$, $\bar{Z}^{(21)}$, and $\bar{Z}^{(12)}$  
    do not include either $h_{1, Y Y}$ or $h_{1, Y Y Y}$,  
    so $\bar{Z}^{(30)}$, $\bar{Z}^{(21)}$, $\bar{Z}^{(12)} = O(\lambda^k)$;  
    the $\bar{Z}^{(03)}$ includes $h_{1, Y Y Y}$,  
    so $\bar{Z}^{(03)} = O(1)$;  
  \item the $\bar{Y}^{(30)}$, $\bar{Y}^{(21)}$ include none of  
    $h_{2, Y Y}$, $h_{2, Z Y Y}$, $h_{2, W Y Y}$, or $h_{2, Y Y Y}$,  
    so $\bar{Y}^{(30)}$, $\bar{Y}^{(21)} = O(1)$;  
    the $\bar{Y}^{(12)}$ includes $h_{2, Z Y Y}$ and $h_{2, W Y Y}$,  
    so $\bar{Y}^{(12)} = O(\gamma^k)$;  
    $\bar{Y}^{(03)}$ includes $h_{2, Y Y Y}$,  
    so $\bar{Y}^{(03)} = O(\gamma^k)$.  
\end{itemize}
Summarizing the above results, the desired statement is proved.  
\end{proof}

\subsubsection{Estimate for new coordinates} \label{s42-uv}
From \eqref{e42-bZbY1} and \eqref{e42-bZbY2},
using \eqref{eq-03-NF_higher2} and \eqref{e42-wc},
we have
\begin{align}
  \bar{Y}^{(10)} = \mathrm{e}^{k \rho_k} + O(\hat{\lambda}^k \gamma^k), \quad
  \bar{Y}^{(02)} = \gamma^k (d + O(\gamma^k \hat{\gamma}^{-k})).
  \label{e42-bY}
\end{align}
Since we have \eqref{eq-03-solrho} and $d \neq 0$,
\begin{align}
  \bar{Y}^{(10)} \neq 0, \quad
  \bar{Y}^{(02)} \neq 0
  \label{e42-Yneq}
\end{align}
for any $k > \kappa(\delta_\mathrm{dom})$
by replacing $\kappa(\delta_\mathrm{dom})$ with a larger one  
according to Remark~\ref{rem-smaller}.
Consider the new coordinates $(u, v)$ such that
\begin{align}
  (u, v)^\mathsf{T} := P_k^{-1} (Z, Y)^\mathsf{T},
  \label{e42-uv}
\end{align}
where
\begin{align}
  P_k = P_k(t, \omega)
    := -\frac{4 \sin^2 \psi}{\bar{Y}^{(10)} \bar{Y}^{(02)}}
    \begin{pmatrix}
      1 & \frac{\bar{Z}^{(10)} - \cos \psi}{\sin \psi} \\
      0 & \frac{\bar{Y}^{(10)}}{\sin \psi}
    \end{pmatrix}.
  \label{e42-Pk}
\end{align}
Note that $P_k$ is well-defined and regular by \eqref{e42-Yneq}
and $\psi \in (0, \pi)$.  
The matrix $P_k$ is chosen so that  
when $T_k|_{W^c_\mathrm{loc}(Q)}$ is written in the $(u, v)$ coordinates as  
$(u, v) \mapsto (\bar{u}, \bar{v})$,
its Taylor expansion at $(u, v) = 0$ is given by
\begin{align}
\begin{aligned}
  \bar{u} &= u \cos \psi - v \sin \psi  
    + \sum_{2 \leq p + q \leq 3} \bar{u}^{(pq)} u^p v^q
    + O(\|(u, v)\|^4), \\  
  \bar{v} &= u \sin \psi + v \cos \psi  
    + \sum_{2 \leq p + q \leq 3} \bar{v}^{(pq)} u^p v^q
    + O(\|(u, v)\|^4),
\end{aligned}\label{e42-bubv}
\end{align}
where $O(\cdot)$ are terms of fourth order or higher.

The estimate of $\bar{u}^{(pq)}$ and $\bar{v}^{(pq)}$ is given as follows.

\begin{lemma}[Estimate of $\bar{u}^{(pq)}$ and $\bar{v}^{(pq)}$]
\label{l42-RZYuv}
  For any $p$, $q \geq 0$ with $2 \leq p + q \leq 3$,
  \begin{align}
    &\bar{u}^{(pq)}, \,
    \bar{v}^{(pq)} = O(\lambda^k) \quad \text{if} \quad (p, q) \neq (0, 2),
      \label{e42-Ebubv} \\
    &\bar{u}^{(02)} = -4 \cos \psi + O(\lambda^k), \quad
    \bar{v}^{(02)} = -4 \sin \psi + O(\lambda^k).
      \label{e42-E2bubv}
  \end{align}
\end{lemma}

\begin{proof}[Proof of Lemma~\ref{l42-RZYuv}]
First, we will prove that
\begin{align}
\begin{aligned}
  &\bar{u}^{(pq)} = O(\lambda^{(p + q - 1) k})
    \sum_{\substack{p' + q' = p + q, \\ p' \geq p}}
      (\bar{Z}^{(p'q')} + \bar{Y}^{(p'q')}), \\
  &\bar{v}^{(pq)} = O(\lambda^{(p + q - 1) k})
    \sum_{\substack{p' + q' = p + q, \\ p' \geq p}} \bar{Y}^{(p'q')}
\end{aligned} \label{e42-bubvG}
\end{align}
for any $p$, $q \geq 0$ with $2 \leq p + q \leq 3$.
Note that $\bar{Y}^{(10)} = \mathrm{e}^{k \rho_k} + O(\hat{\lambda}^k \gamma^k)$.
From \eqref{e42-bY}, \eqref{e42-Pk} implies
\begin{align}
  P_k = O(\gamma^{-k})
    \begin{pmatrix}
      1 & 1 \\
      0 & 1
    \end{pmatrix}, \quad
  P_k^{-1} = O(\gamma^k)
  \begin{pmatrix}
    1 & 1 \\
    0 & 1
  \end{pmatrix}.
  \label{e42-Pkev}
\end{align}
Recall that $T_k|_{W^c_\mathrm{loc}(Q)}$ is written in the $(Z, Y)$ coordinates
as \eqref{e42-bZbY},  
and in the $(u, v)$ coordinates as \eqref{e42-bubv}
with the coordinate transformation given by \eqref{e42-uv}.  
Therefore,  
\begin{align*}
  (\bar{u}, \bar{v})^\mathsf{T} = P_k^{-1} (\bar{Z}, \bar{Y})^\mathsf{T},  
\end{align*}  
where $(\bar{Z}, \bar{Y})$ is obtained from \eqref{e42-bZbY}  
by substituting
$P_k (u, v)^\mathsf{T} = O(\gamma^{-k})(u + v, v)^\mathsf{T}$ into $(Z, Y)$.  
Then,
\begin{align}
  (\bar{u}, \bar{v})^\mathsf{T} = O(\gamma^k)
    \begin{pmatrix}
      1 & 1 \\
      0 & 1
    \end{pmatrix}
    \begin{pmatrix}
      \sum_{1 \leq p + q \leq 3} \bar{Z}^{(pq)}
        O(\gamma^{-(p + q)k}) (u + v)^p v^q
      + O(\|(u, v)\|^4) \\
      \sum_{1 \leq p + q \leq 3} \bar{Y}^{(pq)}
        O(\gamma^{-(p + q)k}) (u + v)^p v^q
      + O(\|(u, v)\|^4)
    \end{pmatrix},
  \label{e42-buvT}
\end{align}
where $O(\|(u, v)\|^4)$ are terms of fourth order or higher.
For any $p$, $q \geq 0$ with $2 \leq p + q \leq 3$,  
we obtain \eqref{e42-bubvG} by comparing the coefficients of $u^p v^q$,  
where note that $O(\gamma^{-k}) = O(\lambda^k)$  
since $\lambda^k \gamma^k = \mathrm{e}^{k \rho_k} = O(1)$.

Next, we will prove \eqref{e42-Ebubv} and \eqref{e42-E2bubv}.
The former follows from \eqref{e42-bubvG} and Lemma~\ref{l42-evbZbY}.  
Thus, it remains to compute $\bar{u}^{(02)}$ and $\bar{v}^{(02)}$ explicitly  
to verify the latter.  
One can find that the inverse matrix $P_k^{-1}$ of the matrix $P_k$ in \eqref{e42-Pk} is  
\begin{align*}
  P_k^{-1} = -\frac{\bar{Y}^{(02)}}{4 \sin^2 \psi}  
    \begin{pmatrix}
      \bar{Y}^{(10)} & -(\bar{Z}^{(10)} - \cos \psi) \\
      0 & \sin \psi
    \end{pmatrix}.
\end{align*}  
Let us denote $P_k$ and its inverse as  
\begin{align*}
  P_k =
  \begin{pmatrix}
    p_{11} & p_{12} \\
    0 & p_{22}
  \end{pmatrix}, \quad
  P_k^{-1} =
  \begin{pmatrix}
    \tilde{p}_{11} & \tilde{p}_{12} \\
    0 & \tilde{p}_{22}
  \end{pmatrix}.
\end{align*}  
Then \eqref{e42-buvT} can be written explicitly as  
\begin{align*}
  (\bar{u}, \bar{v})^\mathsf{T} =  
    \begin{pmatrix}
      \tilde{p}_{11} & \tilde{p}_{12} \\
      0 & \tilde{p}_{22}
    \end{pmatrix}
    \begin{pmatrix}
      \sum_{1 \leq p + q \leq 3} \bar{Z}^{(pq)}
        (p_{11} u + p_{12} v)^p (p_{22} v)^q
      + O(\|(u, v)\|^4) \\
      \sum_{1 \leq p + q \leq 3} \bar{Y}^{(pq)}
        (p_{11} u + p_{12} v)^p (p_{22} v)^q
      + O(\|(u, v)\|^4)
    \end{pmatrix}.
\end{align*}
By comparing the coefficients of $v^2$ and using Lemma~\ref{l42-evbZbY}, we obtain  
\begin{align*}
  \bar{u}^{(02)} &= \tilde{p}_{11} (
    \bar{Z}^{(20)} p_{12}^2
    + \bar{Z}^{(11)} p_{12} p_{22}
    + \bar{Z}^{(02)} p_{22}^2
  ) + \tilde{p}_{12} (
    \bar{Y}^{(20)} p_{12}^2
    + \bar{Y}^{(11)} p_{12} p_{22}
    + \bar{Y}^{(02)} p_{22}^2
  ) \\
  &= \bar{Y}^{(02)} \tilde{p}_{12} p_{22}^2 + O(\lambda^k) \\
  &= 4(\bar{Z}^{(10)} - \cos \psi) + O(\lambda^k) \\
  &= -4 \cos \psi + O(\lambda^k).
\end{align*}
Similarly, we obtain  
\begin{align*}
  \bar{v}^{(02)} &= \tilde{p}_{22} (
    \bar{Y}^{(20)} p_{12}^2
    + \bar{Y}^{(11)} p_{12} p_{22}
    + \bar{Y}^{(02)} p_{22}^2
  ) \\
  &= \bar{Y}^{(02)} \tilde{p}_{22} p_{22}^2 + O(\lambda^k) \\
  &= -4 \sin \psi + O(\lambda^k).
\end{align*}
This completes the verification of \eqref{e42-E2bubv}  
and the proof of the claim.  
\end{proof}

\subsubsection{Calculation of the Lyapunov coeffcient} \label{s42-prf}
Using the new coordinates $(u, v)$ defined in the previous section,  
we define the complex coordinate $z = u + \mathrm{i} v$.  
We naturally identify $\mathbb{R}^2$ with $\mathbb{C}$,  
and express $T_k|_{W^c_\mathrm{loc}(Q)}$ as $z \mapsto \tilde{z}$
using the complex coordinate $z$.  
This can be expanded as  
\begin{align}
  \tilde{z} = \nu_1 z
    + \sum_{2 \leq p + q \leq 3} \tilde{z}^{(pq)} z^p \bar{z}^q
    + O(|z|^4), \quad
  \nu_1 = \cos \psi + \mathrm{i} \sin \psi,
  \label{e42-tz}
\end{align}  
where $O(|z|^4)$ is a term of fourth order or higher.
Then the following holds.  

\begin{lemma}[Estimate of $\tilde{z}^{(pq)}$] \label{l42-evtz}
  We have
  \begin{align}
    \tilde{z}^{(20)}, \,
    \tilde{z}^{(02)}
      = \nu_1 + O(\lambda^k), \quad
    \tilde{z}^{(11)}
      = -2 \nu_1 + O(\lambda^k), \quad
    \tilde{z}^{(pq)}
      = O(\lambda^k)
    \label{e42-tzs}
  \end{align}
  for any $p$, $q \geq 0$ with $p + q = 3$.
\end{lemma}

\begin{proof}[Proof of Lemma~\ref{l42-evtz}]
  Recall the expression \eqref{e42-bubv} of $T_k|_{W^c_\mathrm{loc}(Q)}$
  in the $(u, v)$ coordinates.  
  The coordinate transformation between $z$ and $(u, v)$ is given by  
  \begin{align*}
    z = u + \mathrm{i} v, \quad
    \bar{z} = u - \mathrm{i} v, \quad
    u = \frac{1}{2}(z + \bar{z}), \quad
    v = \frac{1}{2\mathrm{i}}(z - \bar{z}),
  \end{align*}  
  so we have  
  \begin{align*}
    \tilde{z}
      &= \bar{u} + \mathrm{i} \bar{v} \\
      &= (u\cos\psi - v\sin\psi) + \mathrm{i}(u\sin\psi + v\cos\psi)
        + \sum_{2 \leq p + q \leq 3}
          (\bar{u}^{(pq)} + \mathrm{i}\bar{v}^{(pq)}) u^p v^q
        + O(\|(u, v)\|^4) \\
      &= (\cos\psi + \mathrm{i}\sin\psi)z
        + \sum_{2 \leq p + q \leq 3} \frac{(-\mathrm{i})^q}{2^{p + q}}
          (\bar{u}^{(pq)} + \mathrm{i}\bar{v}^{(pq)}) (z + \bar{z})^p (z - \bar{z})^q
        + O(\|(u, v)\|^4),
  \end{align*}  
  where $O(\|(u, v)\|^4)$ is a term of fourth order or higher.
  Calculating the cases $p + q = 2$, $3$ in the second term of the last equation yields
  \begin{align*}
    &\frac{1}{4} (\bar{u}^{(20)} + \mathrm{i}\bar{v}^{(20)})
      (z + \bar{z})^2
    + \frac{-\mathrm{i}}{4} (\bar{u}^{(11)} + \mathrm{i}\bar{v}^{(11)})
      (z + \bar{z}) (z - \bar{z})
    + \frac{(-\mathrm{i})^2}{4} (\bar{u}^{(02)} + \mathrm{i}\bar{v}^{(02)})
      (z - \bar{z})^2, \\
    &\frac{1}{8} (\bar{u}^{(30)} + \mathrm{i}\bar{v}^{(30)})
      (z + \bar{z})^3
    + \frac{-\mathrm{i}}{8} (\bar{u}^{(21)} + \mathrm{i}\bar{v}^{(21)})
      (z + \bar{z})^2 (z - \bar{z}) \\
    &\quad + \frac{(-\mathrm{i})^2}{8} (\bar{u}^{(12)} + \mathrm{i}\bar{v}^{(12)})
      (z + \bar{z}) (z - \bar{z})^2
    + \frac{(-\mathrm{i})^3}{8} (\bar{u}^{(03)} + \mathrm{i}\bar{v}^{(03)})
    (z - \bar{z})^3,
  \end{align*}
  respectively.
  By expanding the above equations and looking at the coefficient of $z^p \bar{z}^q$,
  Lemma~\ref{l42-RZYuv} yields the desired result \eqref{e42-tzs}.
  We complete the proof.
\end{proof}

Finally, let's prove the main consequence.

\begin{proof}[Proof of Proposition~\ref{p42-LC}]
  Using Lemma~\ref{l42-evtz}, we have
  \begin{align*}
    &\Re \left(
      -\tilde{z}^{(21)} \overline{\nu_1}
    \right)
      = O(\lambda^k), \quad
    \Re \left(
      |\tilde{z}^{(02)}|^2
      \frac{-4 + 2\overline{\nu_1}^3}{-2 + \nu_1^3 + \overline{\nu_1}^3}
    \right)
      = \frac{2 + 3\cos\psi - 4\cos^3\psi}{1 + 3\cos\psi - 4\cos^3\psi}
        + O(\lambda^k), \\
    &\Re \left(
      |\tilde{z}^{(11)}|^2
      \frac{-2\overline{\nu_1} + \overline{\nu_1}^2}{(-1 + \overline{\nu_1})^2}
    \right)
      = \frac{2(-2 + \cos\psi)}{-1 + \cos\psi}
      + O(\lambda^k), \\
    &\Re \left(
      \tilde{z}^{(11)} \tilde{z}^{(20)}
      \frac{2 - 6\overline{\nu_1} + \overline{\nu_1}^2}{(-1 + \nu_1)^2}
    \right)
      = - \frac{3(-2 + \cos\psi)}{-1 + \cos\psi}
      + O(\lambda^k).
  \end{align*}
  Adding the above quantities,
  Proposition~\ref{p41-FmlaLC} implies the desired result \eqref{e42-LCQ}.
  We complete the proof.
\end{proof}

\subsection{Parameters for weakly repelling behavior} \label{s43-rep}
In this section, we clarify the region of $(t, \omega)$ where $Q_k$ becomes
weakly repelling on the local center manifold.

\medskip

Recall that $Q_k$
has the multipliers $\nu_1$, $\nu_2$, and $\nu_3$ given by
\begin{align}
  \nu_1 = \cos\psi + \mathrm{i} \sin\psi, \quad
  \nu_2 = \cos\psi - \mathrm{i} \sin\psi, \quad
  |\nu_3| < 1
  \label{e51-nu123}
\end{align}
by restricting $\rho = \rho_k(t, \omega)$ and
$\mu = \mu_k(t, \omega, \rho_k)$
for any $k \in \mathbb{Z}_{> \kappa(\delta_\mathrm{dom})} \cap 2\mathbb{Z}$
and $(t, \omega) \in \{t_k^-(\omega) < t < t_k^+(\omega), \,
\omega \in I_k^\mathrm{ps}\}$.
Recall the constant $\psi_\mathrm{bd} = \pi/20$
in Remark~\ref{r42-psibd}.
Solving the equation
\begin{align*}
  \psi(t, \omega) = \frac{\pi}{2} - \psi_\mathrm{bd}
\end{align*}
by the implicit function theorem,
we obtain a solution
\begin{align*}
  t_k^{+0} : I_k^\mathrm{ps} \to \mathbb{R}.
\end{align*}
The solvability of this equation follows from~\eqref{eq-03-Tkt}
and the relation
$\psi(t, \omega) = \arccos(\frac{\Sigma_k(t, \omega)}{2})$;
see step~(2) in the proof of Lemma~\ref{eq-03-Tkt}.
In particular, the $t_k^{+0}$ is $C^{r - 2}$ and
\begin{align*}
  t_k^{+0}(\omega) = O(1), \quad
  t_k^-(\omega) < t_k^{+0}(\omega) < t_k^+(\omega)
\end{align*}
since we have
\begin{align*}
  \Sigma_k(t_k^-, \omega) = -2
  < \Sigma_k(t_k^{+0}, \omega) = 2\cos(\frac{\pi}{2} - \psi_\mathrm{bd})
  < \Sigma_k(t_k^+, \omega) = 2.
\end{align*}
We define the open set $\mathcal{R}_k^\mathrm{rep}$ by
\begin{align}
  \mathcal{R}_k^\mathrm{rep} &:= \{(t, \omega) \:|\:
    t_k^{+0}(\omega) < t < t_k^+(\omega), \,
    \omega \in I_k^\mathrm{ps}\}
    \label{e43-Rkrep}
\end{align}
for any $k \in \mathbb{Z}_{> \kappa(\delta_\mathrm{dom})} \cap 2\mathbb{Z}$.
Then for any $(t, \omega) \in \mathcal{R}_k^\mathrm{rep}$,  
we have $\psi(t, \omega) \in (0, \pi/2 - \psi_\mathrm{bd})
\cap \Psi_\mathrm{reg} = (0, \pi/2 - \psi_\mathrm{bd})$
since the map $(t_k^-(\omega), t_k^+(\omega)) \ni t \mapsto \psi(t, \omega) \in (0, \pi)$
is orientation reversing
for any fixed $\omega \in I_k^\mathrm{ps}$.
Remark~\ref{r42-psibd} implies that for any
$k \in \mathbb{Z}_{> \kappa(\delta_\mathrm{dom})} \cap 2\mathbb{Z}$
and $(t, \omega) \in \mathcal{R}_k^\mathrm{rep}$,
the point $Q_k$ becomes weakly repelling  
on the local center manifold
by replacing $\kappa(\delta_\mathrm{dom})$ with a larger one  
according to Remark~\ref{rem-smaller}.  

\ifdefined\isMaster
\else
  \end{document}
\fi
\ifdefined\isMaster
\else
  \documentclass[11pt,a4paper]{article}

\usepackage[left=2cm,right=2cm,top=2cm,bottom=2cm]{geometry}
\usepackage{graphicx}
\usepackage{multirow}
\usepackage{amsmath,amssymb,amsfonts}
\usepackage{amsthm}
\usepackage{mathrsfs}
\usepackage[title]{appendix}
\usepackage{xcolor}
\usepackage{textcomp}
\usepackage{manyfoot}
\usepackage{booktabs}
\usepackage{algorithm}
\usepackage{algorithmicx}
\usepackage{algpseudocode}
\usepackage{listings}
\ifdefined\debug
  \usepackage[inline]{showlabels}
  \showlabels{ref}
  \showlabels{cite}
  \usepackage{fancyhdr}
\fi
\usepackage{bm}
\usepackage{hyperref}

\theoremstyle{plain}
\newtheorem{theorem}{Theorem}[section]
\newtheorem{proposition}[theorem]{Proposition}%
\newtheorem{lemma}[theorem]{Lemma}
\newtheorem{corollary}[theorem]{Corollary}
\newtheorem{claim}[theorem]{Claim}
\theoremstyle{definition}
\newtheorem{definition}[theorem]{Definition}
\newtheorem{remark}[theorem]{Remark}
\newtheorem{example}[theorem]{Example}
\newtheorem{question}[theorem]{Question}
\newtheorem*{notation}{Notation}
\theoremstyle{plain}
\newtheorem{maintheorem}{Theorem}
\renewcommand{\themaintheorem}{\Alph{maintheorem}}
\newcommand{\wiptext}[1]{\textcolor{blue}{#1}}
\ifdefined\debug
  \newcommand{\debugtext}[1]{\textcolor{red}{#1}}
  \pagestyle{fancy} 
  \fancyhead[L]{}
  \fancyhead[C]{}
  \fancyhead[R]{\textcolor{red}{[DEBUG] Debug mode is activated.}}
  \renewcommand{\headrulewidth}{0pt}
\else
  \newcommand{\debugtext}[1]{\ignorespaces}
\fi

\counterwithin{figure}{section}
\counterwithin{table}{section}

\raggedbottom

\numberwithin{equation}{section}

  \begin{document}
\fi

\section{Creation of a Hopf-homoclinic cycle}
\label{s5-HHC}
We have completed most of the proof of Theorem~\ref{thm-third}
due to Proposition~\ref{p34-Qk}.
In this section,
we prove the remainder of that proof,
namely the existence of a homoclinic point to $Q_k = Q_k(t, \omega, \rho_k)$,
and we complete the proof of Theorem~\ref{thm-third},
where $Q_k$ is the non-hyperbolic periodic point in Proposition~\ref{p34-Qk}.
In this section, we always choose $(t, \omega)$ so that $Q_k$ becomes weakly repelling:
$(t, \omega) \in \mathcal{R}_k^\mathrm{rep}$,
where $\mathcal{R}_k^\mathrm{rep}$ is the set in \eqref{e43-Rkrep}.
In Section~\ref{s51-TInt},
we observe that the two-dimensional generalized unstable manifold
$\widetilde{W}^u(Q_k)$ defined by \eqref{e12-gmfds}
and the two-dimensional stable manifold $W^s(O(\mu_k, \omega, \rho_k))$ intersect
when $(f, \Gamma)$ holds accompanying condition \textbf{(AC)},
which is defined in Section~\ref{s23-resfam},
where $O(\mu, \omega, \rho)$ is a continuation of $O^*$.
In Section~\ref{s52-Hint},
we see that the two-dimensional generalized unstable manifold
$\widetilde{W}^u(Q_k)$ intersects
the one-dimensional generalized stable manifold
$\widetilde{W}^s(Q_k)$ by adjusting $\omega$ and
giving the proof of Theorem~\ref{thm-third}.

\subsection{Transverse intersection between the unstable and stable manifolds}
\label{s51-TInt}
Let us recall the accompanying condition \textbf{(AC)} defined in
Section~\ref{s23-resfam}.
The goal of this section is to prove the following proposition.
In the following, we write $\rho_k = \rho_k(t, \omega)$
and $\mu_k = \mu_k(t, \omega, \rho_k)$.

\begin{proposition}[$\widetilde{W}^u(Q_k) \cap W^s(O(\mu_k, \omega, \rho_k))
\neq \emptyset$]
\label{p51-Tint}
  Suppose that $(f, \Gamma)$ satisfies the accompanying condition \textbf{(AC)}.
  Then by replacing $\kappa(\delta_\mathrm{dom})$ with a larger one  
  according to Remark~\ref{rem-smaller},
  the two-dimensional generalized unstable manifold $\widetilde{W}^u(Q_k)$
  intersects the two-dimensional stable manifold
  $W^s(O(\mu_k, \omega, \rho_k))$
  transversely for any
  $(t, \omega) \in \mathcal{R}_k^\mathrm{rep}$.
\end{proposition}

\begin{remark}
  Although the situation is slightly different in the sense that
  $Q_k$ is non-hyperbolic, the idea of the proof is the same as
  in \cite{LLST2022}.
  For the sake of completeness, the proof is given below.
\end{remark}

The proof of the above proposition will be given in the following subsections.
\begin{itemize}
  \item In Section~\ref{s51-det}, we prove that the restriction
    of the first-return map $T_k$ to the local center manifold of $Q_k$
    is area expanding on an annular region excluding $Q_k$.
  \item In Section~\ref{s51-GUM}, using the above area expansion,
    we show that the generalized unstable manifold
    $\widetilde{W}^u(Q_k)$ becomes sufficiently
    large in the $Y$-direction.
  \item Finally, in Section~\ref{s51-Tint}, we prove
    Proposition~\ref{p51-Tint}.
\end{itemize}

\subsubsection{Area expanding property} \label{s51-det}
Recall the coordinates $(Z, Y, W)$ defined in Proposition~\ref{p32-NF}.
As explained in Section~\ref{s42-ZY}, the two-dimensional local center manifold
of $Q_k = (Z_Q, Y_Q, W_Q)$ exists and can be written in the form
\begin{align*}
  W = w^c_0(Z, Y),
\end{align*}
where $w^c_0$ is at least $C^5$,
$w^c_0$ is a map from a small open two-dimensional disk $D^c_0$ centered at $(Z_Q, Y_Q)$
to a small open one-dimensional disk centered at $W_Q$,
and satisfies $w^c_0(Z_Q, Y_Q) = W_Q$.
Let $S^c_0 := \{(Z, Y, w^c_0(Z, Y)) \:|\: (Z, Y) \in D^c_0\}$.

We begin by extending $S^c_0$ as follows.
Since $D^c_0$ can be taken sufficiently small,
we may assume that $S^c_0 \subset \Pi_k'$ initially,
where $\Pi_k' = [-\delta_\mathrm{dom}', \delta_\mathrm{dom}']^3$
is the domain of the $(Z, Y, W)$ coordinates;
see Proposition~\ref{p32-NF}.
The surface $S^c_0$ is tangent to the center-unstable cone field
$\mathcal{C}^{cu}$, meaning that at every point
$M \in S^c_0$, the tangent space satisfies
$T_M S^c_0 \subset \mathcal{C}^{cu}(M)$;
see \eqref{eq-03-conecu} for the definition of the center-unstable cone field.
By Proposition~\ref{p33-CF},
if $S^c_i$ is tangent to $\mathcal{C}^{cu}$,
then so is $S^c_{i+1}$, where $S^c_{i+1}$ is defined as
the connected component of $T_k(S^c_i) \cap \Pi_k'$
that contains $Q_k$, for each $i \in \{0, 1, 2, \cdots\}$.
In this way, we define $S^c_i$ inductively for all
$i \in \{0, 1, 2, \cdots\}$ and obtain the surface
in $\Pi_k'$ tangent to $\mathcal{C}^{cu}$:
\begin{align*}
  S^c := \bigcup_{i = 0}^\infty S^c_i.
\end{align*}
By replacing $\hat{\delta}_\mathrm{dom}$
and $\kappa(\delta_\mathrm{dom})$ with smaller and larger ones
according to Remark~\ref{rem-smaller},  
the size of the center-unstable cone in \eqref{eq-03-conecu}
can be made arbitrarily small.
Therefore, $S^c$ can be parameterized as  
\begin{align*}
  W = w^c(Z, Y),
\end{align*}  
where $w^c$ is at least $C^5$ and satisfies  
$w^c(Z_Q, Y_Q) = W_Q$.  
The domain of $w^c$ is the projection of $S^c$
onto the $(Z, Y)$-plane, which is a two-dimensional disk
containing $(Z_Q, Y_Q)$; we denote it by $D^c$.

\medskip

The sets $T_k^{-n}(S^c_n)$, $n \in \{0, 1, \cdots\}$,  
form a nested family of sets as  
\begin{align*}
  \cdots \subset T_k^{-2}(S^c_2)  
  \subset T_k^{-1}(S^c_1)  
  \subset S^c_0  
\end{align*}  
and any point in $T_k^{-n}(S^c_n)$ remains within $\Pi_k'$  
under $n$ iterations of $T_k$.  
Therefore, letting $\Omega^c_n$,
$n \in \{0, 1, \cdots\}$,
denote the projection of  
$T_k^{-n}(S^c_n)$ onto the $(Z, Y)$-plane,  
we can write for any $(Z_0, Y_0) \in \Omega^c_n$,  
\begin{align*}
  (Z_i, Y_i, W_i) = T_k^i(Z_0, Y_0, w^c(Z_0, Y_0)),  
  \quad i \in \{0, 1, \cdots, n\}.  
\end{align*} 
Our goal is to demonstrate the following fact.

\begin{lemma}[Area expanding property]\label{lem-AEP}
  For any $\Delta > 0$,
  there exists $n = n(\Delta) \in \{1, 2, \cdots\}$ such that
  \begin{align*}
    \left| \det \frac{\partial (Z_n, Y_n)}{\partial (Z_0, Y_0)} \right|
    \geq 2
  \end{align*}
  for any $(Z_0, Y_0) \in \Omega^c_n \cap \{|Z_0 - Z_Q|$, $|Y_0 - Y_Q| \geq \Delta\}$,
  $(t, \omega) \in \mathcal{R}_k^\mathrm{rep}$.
\end{lemma}

\begin{remark}
  If $\Omega^c_n$ becomes strictly decreasing with respect to inclusion,  
  it may happen that  
  $\Omega^c_n \cap \{|Z_0 - Z_Q|, \, |Y_0 - Y_Q| \geq \Delta\} = \emptyset$.  
  However, the fact that $\Omega^c_n$ becomes
  strictly decreasing with respect to inclusion
  implies that $S^c_i$ grows sufficiently to reach the boundary of $\Pi_k'$,  
  in which case it is not necessary to explicitly state  
  the area expanding property;  
  see Lemma~\ref{l51-eumfd} and its proof.    
\end{remark}

\begin{proof}[Proof of Lemma~\ref{lem-AEP}]
We divide the proof into steps.

\medskip

\textbf{(1) Reduction to the argument in the $(u, v)$ coordinates.} \quad
As in Section~\ref{s42-ZY},  
we prepare new coordinates $(Z^{new}, Y^{new},  
W^{new})$ centered at $Q_k$.  
Since these coordinates are defined via translation,  
it suffices to verify the area expanding property
in these new coordinates.  
Hereafter, we drop the `new'.
As in Section~\ref{s42-ZY},  
if we write
$T_k|_{S^c \cap T_k^{-1}(S^c)}: (Z, Y) \mapsto (\bar{Z}, \bar{Y})$,
then the $(\bar{Z}, \bar{Y})$ is given by \eqref{e42-bZbYo}.

Recall the coordinates $(u, v)$ prepared in  
Section~\ref{s42-uv}, defined as \eqref{e42-uv}. 
Let
\begin{align*}
  (u_i, v_i)^\mathsf{T} = P_k^{-1}(Z_i, Y_i)^\mathsf{T}
\end{align*} 
for any $i \in \{0, 1, \cdots, n\}$,
where $P_k$ is the matrix defined in \eqref{e42-Pk}.
Since the coordinate transformation is defined  
by a linear map via the matrix $P_k$, we have  
\begin{align*}
  \det \frac{\partial (Z_n, Z_n)}{\partial (Z_0, Y_0)}
  = \det \frac{\partial (u_n, v_n)}{\partial (u_0, v_0)}.
\end{align*}
Therefore, it suffices to show the area expanding property  
with respect to $(u, v)$.

\medskip

\textbf{(2) Reduction to the argument in the $(x, y)$ coordinates.} \quad
Recall the complex coordinate $z = u + \mathrm{i} v$ defined in  
Section~\ref{s42-prf}.  
Although this complex coordinate $z$ does not give the normal form  
\eqref{eq-01-Cform}, a new complex coordinate
\begin{align}
  w = z + \sum_{p + q = 2} \frac{\tilde{z}^{(pq)}}{\nu_1 - \nu_1^p \nu_2^q}
  z^p \bar{z}^q
  \label{e51-weq}
\end{align}
via \eqref{eq-04-z2w} does provide the normal form,
where the $\tilde{z}^{(pq)}$ in the above equation
are estimated as in Lemma~\ref{l42-evtz}.
We define
$z_i$,
$w_i$,
and $(x_i, y_i) \in \mathbb{R}^2$ as
\begin{align}
  z_i = u_i + \mathrm{i} v_i, \quad
  w_i = z_i + \sum_{p + q = 2}
    \frac{\tilde{z}^{(pq)}}{\nu_1 - \nu_1^p \nu_2^q}
    z_i^p \bar{z_i}^q
  = x_i + \mathrm{i} y_i
  \label{e51-xytzwbxy}
\end{align}
for any $i \in \{0, 1, \cdots, n\}$.
Then, by the chain rule, we can compute
\begin{align}
  \det \frac{\partial (u_n, v_n)}{\partial (u_0, v_0)}
  = \det \frac{\partial (u_n, v_n)}
    {\partial (x_n, x_n)} \cdot
    \det \frac{\partial (x_n, y_n)}
    {\partial (x_0, y_0)} \cdot
    \det \frac{\partial (x_0, y_0)}{\partial (u_0, v_0)}
  \label{e51-deddd}
\end{align}
The key question is whether absolute value of this becomes greater than $1$.

\medskip

\textbf{(3) Computation of the area expansion ratio  
in the normal form.} \quad
We begin by computing the middle term  
on the right-hand side of \eqref{e51-deddd}:
\begin{align}
  \det \frac{\partial (x_n, y_n)}{\partial (x_0, y_0)}
  = \prod_{i = 0}^{n - 1} \det \frac{\partial (x_{i + 1}, y_{i + 1})}
    {\partial (x_i, y_i)}.
  \label{e51-dxynxy0}
\end{align}
In general, for a complex function  
$g: \mathbb{C} \ni w \mapsto \tilde{w} \in \mathbb{C}$,  
its Jacobian determinant is given by  
$|g_w|^2 - |g_{\bar{w}}|^2$.
Therefore, if $g$ is given in the normal form \eqref{eq-01-Cform},  
its Jacobian determinant can be computed as  
\begin{align*}
  |g_w|^2 - |g_{\bar{w}}|^2 
  &= (\nu + 2\alpha w \bar{w})(\bar{\nu} + 2\bar{\alpha} \bar{w} w) 
    + O(|w|^3) \\
  &= 1 + 2(\bar{\nu} \alpha + \nu \bar{\alpha}) w \bar{w}
    + O(|w|^3) \\
  &= 1 + 4 \Re (\bar{\nu} \alpha) |w|^2 + O(|w|^3) \\
  &= 1 - 4 \mathrm{LC} |w|^2 + O(|w|^3),
\end{align*}
where $\mathrm{LC} = -\Re(\bar{\nu} \alpha)$  
is the Lyapunov coefficient at the origin.
Thus, we obtain  
\begin{align*}
  \frac{\partial (x_{i + 1}, y_{i + 1})}{\partial (x_i, y_i)}  
  &= 1 - 4 \mathcal{L}(\psi) |w_i|^2
  + O(\lambda^k) + O(|w_i|^3) \\
  &= 1 - 4 \mathcal{L}(\psi) |z_0|^2
  + O(\lambda^k) + O(|z_0|^3)
\end{align*}
for any $i \in \{0, 1, \cdots, n - 1\}$,
where $\mathcal{L}(\psi)$ is the function defined in \eqref{e42-LCQ}.
Here, we used $w_i = z_i + O(|z_i|^2)$ from \eqref{e51-weq}
and $|z_i| = |z_0| + O(|z_0|^2)$ from \eqref{e42-tz}.
Hence, from \eqref{e51-dxynxy0}, we conclude
\begin{align*}
  \det \frac{\partial (x_n, y_n)}{\partial (x_0, y_0)}
  = 1 - 4 n \mathcal{L}(\psi) |z_0|^2 + O(\lambda^k) + O(|z_0|^3).
\end{align*}

\medskip

\textbf{(4) Computation of the area expansion ratio  
in the coordinate transformation.} \quad
Next, we compute the leftmost and rightmost terms  
on the right-hand side of \eqref{e51-deddd}.  
From \eqref{e51-xytzwbxy},
\eqref{e51-nu123},
and Lemma~\ref{l42-evtz}, we have
\begin{align*}
  \det \frac{\partial(x_0, y_0)}{\partial(u_0, v_0)} 
  &= |\partial_{z_0} w_0|^2 - |\partial_{\bar{z}_0} w_0|^2 \\
  &= \left| 1 + \frac{2z_0}{1 - \nu_1} - \frac{2\bar{z}_0}{1 - \nu_2} \right|^2
    - \left| - \frac{2z_0}{1 - \nu_2} + \frac{2\bar{z}_0}{1 - \nu_2^3} \right|^2
    + O(\lambda^k) \\
  &= 1 + O(|z_0|^2) + O(\lambda^k).
\end{align*}
\debugtext{
Actually, we have
\begin{align}
  \begin{aligned}
    \det \frac{\partial(x_0, y_0)}{\partial(u_0, v_0)} 
    &= 1 + \mathcal{L}(\psi) (u_0, v_0)
      A_L
       \begin{pmatrix}
         u_0 \\ v_0
       \end{pmatrix} + O(\lambda^k),
  \end{aligned}
  \end{align}
where
\begin{align*}
  &A_L = A_L(\psi) :=
    \begin{pmatrix}
      L_{20} & \frac{L_{11}}{2} \\
      \frac{L_{11}}{2} & L_{02}
    \end{pmatrix}, \quad
  L_{20} = L_{20}(\psi)
    := -4(1 + \cos\psi), \\
  &L_{11} = L_{11}(\psi)
    := -\frac{4 \sin\psi}{\cos\psi (1 + 2 \cos\psi)}, \quad
  L_{02} = L_{02}(\psi)
    := 4 \cos\psi.
\end{align*}
}
Since the inverse of \eqref{e51-weq} is given by \eqref{e41-w2z},  
we have
\begin{align*}
  z_n = w_n - \sum_{p + q = 2}
    \frac{\tilde{z}^{(pq)}}{\nu_1 - \nu_1^p \nu_2^q}
    w_n^p \bar{w}_n^q
    + O(|w_n|^3).
\end{align*}
Note that since we have \eqref{e51-xytzwbxy} and \eqref{e42-bubv},
\begin{align*}
  |w_n| = O(|z_0|).
\end{align*}
Using the above two results, we have 
\begin{align*}
  \det \frac{\partial(u_n, v_n)}{\partial(x_n, y_n)} 
  &= |\partial_{w_n} z_n|^2
     - |\partial_{\bar{w_n}} z_n|^2 \\
  &= 1 + O(|w_n|^2) + O(\lambda^k) \\
  &= 1 + O(|z_0|^2) + O(\lambda^k).
\end{align*}

\medskip

\textbf{(5) Possession of the area expanding property.} \quad  
Combining the results of Step~(3) and (4),
there exists a constant $C_1 = C_1(\mathbb{F}) > 0$ such that
\begin{align}  
  \left| \det \frac{\partial (u_n, v_n)}{\partial (u_0, v_0)} \right|
  &\geq 1 + 4 n C |z_0|^2 - C_1(|z_0|^2 + \lambda^k),
\end{align}  
where $C$ is the constant in Remark~\ref{r42-psibd}.
Here, from the assumption and \eqref{e42-Pkev}, we have  
\begin{align*}
  C_2 \Delta |\gamma|^k \leq |z_0|^2 \leq C_3 |\gamma|^k  
\end{align*}  
for some constants $C_2 = C_2(\mathbb{F}) > 0$  
and $C_3 = C_3(\mathbb{F}) > 0$.  
Therefore,  
\begin{align*}  
  \left| \det \frac{\partial (u_n, v_n)}{\partial (u_0, v_0)} \right|  
  &\geq 1 + 4 n C C_2 \Delta |\gamma|^k - C_1(C_3 |\gamma|^k + \lambda^k).  
\end{align*}  
In order for this to be at least $2$, it suffices that  
\begin{align*}  
  n \geq \frac{1 + C_1 (C_3 + 1)}{4 C C_2 \Delta}.  
\end{align*}  
This completes the proof.  
\end{proof}

\subsubsection{Analysis of the size of the two-dimensional
generalized unstable manifold}
\label{s51-GUM}

Recall that the domain where the $(Z, Y, W)$ coordinates are defined
in Proposition~\ref{p32-NF} is
$\Pi_k' = [-\delta_\mathrm{dom}', \delta_\mathrm{dom}']^3$.
By replacing $\hat{\delta}_\mathrm{dom}$  
and $\kappa(\delta_\mathrm{dom})$ with smaller and larger ones  
according to Remark~\ref{rem-smaller},
we have the following.

\begin{lemma}[$\widetilde{W}^u(Q_k)$ is sufficiently large
in the $Y$-direction]
\label{l51-eumfd}
The two-dimensional generalized unstable manifold
$\widetilde{W}^u(Q_k)$ intersects $\{Y = (\delta_\mathrm{dom}')^2\}$ or
$\{Y = -(\delta_\mathrm{dom}')^2\}$
for any $(t, \omega) \in \mathcal{R}_k^\mathrm{rep}$.
\end{lemma}

\begin{proof}[Proof of Lemma~\ref{l51-eumfd}]
Let us recall the manifolds $S^c_i$, $i \in \{0, 1, \cdots\}$,  
which serve as extensions of the local center manifold of $Q_k$,  
defined at the beginning of Section~\ref{s51-det}.    
Let us recall $\Omega^c_i$, $i \in \{0, 1, \cdots\}$,  
defined before Lemma~\ref{lem-AEP}.  
Logically, the following two cases may occur:  
\begin{itemize}
  \item $\Omega^c_{i + 1} \subsetneq \Omega^c_i$ for some $i \in \{0, 1, \cdots\}$, or  
  \item $\Omega^c_{i + 1} = \Omega^c_i$ for any $i \in \{0, 1, \cdots\}$.  
\end{itemize}  

\medskip

\textbf{(1) In the former case.}\quad
There exists $I \geq 0$ such that  
$S^c_{I + 1}$ intersects $\partial^Y \pi_k' \cup \partial^Z \pi_k'$,  
and for all $i \in \{0, 1, \cdots, I\}$,  
$S^c_i$ does not intersect $\partial^Y \pi_k' \cup \partial^Z \pi_k'$,
where  
\begin{align}
\begin{aligned}
  &\pi_k'
    := \{(Z, Y, W) \in \Pi_k' \:|\: |Z| \leq \delta_\mathrm{dom}', \,  
      |Y| \leq (\delta_\mathrm{dom}')^2\}, \\
  &\partial^Y \pi_k'
    := \{(Z, Y, W) \in \pi_k' \:|\: |Y| = (\delta_\mathrm{dom}')^2\}, \\
  &\partial^Z \pi_k' 
    := \{(Z, Y, W) \in \pi_k' \:|\: |Z| = \delta_\mathrm{dom}'\}.  
\end{aligned} \label{e51-pikp}
\end{align}  
Here, $D_0$ is taken sufficiently small so that  
$S^c_0$ does not intersect  
$\partial^Y \pi_k' \cup \partial^Z \pi_k'$,  
by replacing $\kappa(\delta_\mathrm{dom})$ with a larger one  
according to Remark~\ref{rem-smaller} if necessary.  
In fact, by replacing $\hat{\delta}_\mathrm{dom}$  
and $\kappa(\delta_\mathrm{dom})$ with smaller and larger ones  
according to the same remark,  
$S^c_{I + 1}$ must intersect $\partial^Y \pi_k'$.  
Indeed, by the normal form \eqref{eq-03-NF},
as long as $|Z|$, $|Y| \leq O((\delta_\mathrm{dom}')^2)$,  
we have  
\begin{align*}
  |\bar{Z}| \leq O((\delta_\mathrm{dom}')^2).  
\end{align*}  
From the construction of $S^c_i$ and the fact that $Q_k$ is weakly repelling  
on the local center manifold,  
we have that $S^c_i$ is contained in the generalized unstable manifold  
$\widetilde{W}^u(Q_k)$.  
Therefore, the desired statement holds.  

\medskip

\textbf{(2) In the latter case.}\quad
We denote the projection of $S^c_i$ onto the $(Z, Y)$-plane by $D^c_i$,  
for each $i \in \{0, 1, \cdots\}$.  
We take $\Delta > 0$ sufficiently small so that  
$D^c_0 \cap D(\Delta)$,  
with $D(\Delta) := \{|Z - Z_Q|, \, |Y - Y_Q| \geq \Delta\}$,  
is non-empty and homeomorphic to an annulus. 
By Lemma~\ref{lem-AEP}, there exists $n = n(\Delta) > 0$.  
By the assumption, $\Omega^c_{i n} = \Omega^c_0 = D^c_0$
for any $i \in \{0, 1, \cdots\}$,  
and hence we note that $\Omega^c_{i n} \cap D(\Delta) \neq \emptyset$
for any $i \in \{0, 1, \cdots\}$.
In particular, the following inequality holds:  
\begin{align*}
  \mathrm{Area}(D^c_{i n} \cap D(\Delta)) \geq  
  2^i \, \mathrm{Area}(D^c_0 \cap D(\Delta))  
\end{align*}  
for any $i \in \{0, 1, \cdots\}$,
where $\mathrm{Area}(X)$ denotes the Euclidean area of a region $X$  
in the $(Z, Y)$-plane, viewed as $\mathbb{R}^2$.  
Therefore, there exists $I \geq 0$ such that  
$S^c_{I + 1}$ intersects $\partial^Y \pi_k' \cup \partial^Z \pi_k'$,  
and for all $i \in \{0, 1, \cdots, I\}$,  
$S^c_i$ does not intersect $\partial^Y \pi_k' \cup \partial^Z \pi_k'$.
By an argument similar to that in Step~(1),  
$S^c_{I + 1}$ must intersect $\partial^Y \pi_k'$.  
This completes the proof.  
\end{proof}

\subsubsection{Observing a transverse intersection} \label{s51-Tint}

\begin{proof}[Proof of Proposition~\ref{p51-Tint}]
We divide the proof into steps.

\medskip

\textbf{(1) Obtaining a segment of the stable manifold.} \quad
Let us recall that the points $M^-$ and $M^+$  
at $\varepsilon = \varepsilon^*$
were denoted by $M^-_0$ and $M^+_0$, respectively;  
see Section~\ref{s21-PU}.  
Also recall the $\varepsilon$-dependent coordinates  
$(x_1, x_2, y)$ defined in Section~\ref{s31-FRM}.  
By the accompanying condition \textbf{(AC)},  
$W^s(O^*)$ intersects $W^u_\mathrm{loc}(O^*)$  
at the point $M_0^{-+} = (0, 0, y^{-+}_0)$  
in $(x_1, x_2, y)$ coordinates at $\varepsilon = \varepsilon^*$ with  
\begin{align*}
  0 < y^{-+}_0 - y^-_0 < (\delta_\mathrm{dom}')^2,  
\end{align*}  
where $y^-_0$ denotes the $y$-coordinate of $M_0^-$, that is,  
$M_0^- = (0, 0, y^-_0)$.
We denote the small neighborhood of $M_0^{-+}$ in $W^s(O^*)$ as $W_0^{s+}$.
Also, $W^s(O^*)$ intersects $W^u_\mathrm{loc}(O^*)$  
at the point $M_0^{--} = (0, 0, y^{--}_0)$  
in $(x_1, x_2, y)$ coordinates at $\varepsilon = \varepsilon^*$ with  
\begin{align*}
  -(\delta_\mathrm{dom}')^2 < y^{--}_0 - y^-_0 < 0.
\end{align*}
We denote the small neighborhood of $M_0^{--}$ in $W^s(O^*)$ as $W_0^{s-}$.

We denote by $O(\varepsilon)$ the continuation with respect to $\varepsilon$  
of the hyperbolic periodic point $O^*$
with $O(\varepsilon^*) = O^*$.
Since $W_0^{s+}$ intersects $W^u_\mathrm{loc}(O^*)$ transversely,  
by replacing $\delta_\mathrm{prm} > 0$ with a smaller one  
according to Remark~\ref{rem-smaller},  
we can consider the continuations with respect to the parameter $\varepsilon$  
of $M_0^{-+}$, $y_0^{-+}$, and $W_0^{s+}$ as  
$M^{-+}(\varepsilon)$,  
$y^{-+}(\varepsilon)$,  
and $W^{s+}(\varepsilon)$.  
Here, $W^{s+}$ is a subset of $W^s(O)$ that intersects $W^u_\mathrm{loc}(O)$  
transversely at the point $M^{-+} = (0, 0, y^{-+})$  
in $(x_1, x_2, y)$ coordinates, and satisfies  
\begin{align*}
  M^{-+}(\varepsilon^*) = M_0^{-+}, \quad  
  y^{-+}(\varepsilon^*) = y^{-+}_0, \quad  
  W^{s+}(\varepsilon^*) = W_0^{s+}.  
\end{align*}  
In a similar manner, we can consider the continuations  
with respect to the parameter $\varepsilon$  
of $M_0^{--}$, $y_0^{--}$, and $W_0^{s-}$ as  
$M^{--}(\varepsilon)$,  
$y^{--}(\varepsilon)$,  
and $W^{s-}(\varepsilon)$.  
By replacing $\delta_\mathrm{prm} > 0$ with a smaller one
according to the same remark,
we may suppose
\begin{align*}
  -(\delta_\mathrm{dom}')^2 < y^{--} - y^- < 0
    < y^{-+} - y^- < (\delta_\mathrm{dom}')^2.
\end{align*}

\medskip

\textbf{(2) Pullback to the $(Z, Y, W)$ space.} \quad
We write the coordinates near $M^-$ as $(\tilde{x}_1, \tilde{x}_2, \tilde{y})$.
The manifolds $W^{s\sigma}$, $\sigma \in \{+, -\}$,  
can be expressed as the graphs of functions of the form  
\begin{align}
  \tilde{y} = y^{-\sigma} + w^{s\sigma}(\tilde{x}_1, \tilde{x}_2, \varepsilon)  
  \label{e51-wss}
\end{align}  
where $w^{s\sigma}$ and its first and second partial derivatives  
with respect to $(\tilde{x}_1, \tilde{x}_2)$ are $C^{r - 2}$
and $w^{s\sigma}(0, 0, \varepsilon) \equiv 0$.
We write $T_0^k: (x_1, x_2, y)  
\mapsto (\tilde{x}_1, \tilde{x}_2, \tilde{y})$.  
By substituting the expressions for $\tilde{x}_1$, $\tilde{x}_2$  
from \eqref{eq-03-T0kprm} into $\tilde{x}_1$, $\tilde{x}_2$  
in \eqref{e51-wss}, we obtain an equation of  
$(x_1, x_2, \varepsilon, \tilde{y})$.
By Proposition~\ref{pc1-sgl}, this equation can be solved for $\tilde{y}$ as  
\begin{align}
  \tilde{y} = y^{-\sigma} + \hat{w}^{s\sigma}(x_1, x_2, \varepsilon),
  \label{e51-wssR}
\end{align}  
and using Proposition~\ref{pc2-ev},  
we can estimate the partial derivatives as  
\begin{align}
  \hat{w}^{s\sigma}_{x_1}, \,  
  \hat{w}^{s\sigma}_{x_2} = O(\lambda^k).
  \label{e51-hwssEV}
\end{align}  
Indeed, defining  
\begin{align*}
  H_k(x_1, x_2, \varepsilon, \tilde{y})  
    := w^{s\sigma}(\tilde{x}_1, \tilde{x}_2, \varepsilon),  
\end{align*}  
we have  
\begin{align*}
  H_k = O(\lambda^k), \quad  
  H_{k, \tilde{y}} = O(\lambda^k), \quad  
  H_{k, x_1}, \,  
  H_{k, x_2} = O(\lambda^k),  
\end{align*}  
which yields the desired results.

From \eqref{eq-03-Shilnikov}, \eqref{e32-XYshift}, and \eqref{e32-ZW}, we have  
\begin{align*}
  \begin{pmatrix}
    x_1 \\
    x_2
  \end{pmatrix}
  &=
  \begin{pmatrix}
    F_{1, k}(Z, W, \varepsilon) \\
    F_{2, k}(Z, W, \varepsilon)
  \end{pmatrix}, \quad
  \begin{pmatrix}
    F_{1, k} \\
    F_{2, k}
  \end{pmatrix}
  :=
  \begin{pmatrix}
    \frac{1}{\alpha^*} & -\frac{\beta^*}{\alpha^*} \\
    0 & 1
  \end{pmatrix}
  \begin{pmatrix}
    Z \\
    W
  \end{pmatrix}
  +
  \begin{pmatrix}
    x_1^+ + X_{1, k}^* \\
    x_2^+ + X_{2, k}^*
  \end{pmatrix}.
\end{align*}  
Substituting this into $(x_1, x_2)$ in \eqref{e51-wssR}  
and using $Y = \tilde{y} - y^-$, we obtain  
\begin{align}
  Y = y^{-\sigma} - y^- + \tilde{w}^{s\sigma}(Z, W, \varepsilon), \quad
  \tilde{w}^{s\sigma} := \hat{w}^{s\sigma}(F_{1, k}, F_{2, k}, \varepsilon).  
  \label{e51-twss}
\end{align}  
By \eqref{e51-hwssEV},
noting that $\omega \in I_k^\mathrm{bd}$,
we have
\begin{align}
  \tilde{w}^{s\sigma}_Z, \,
  \tilde{w}^{s\sigma}_W = O(\lambda^k).
  \label{e51-twssEV}
\end{align}

\medskip

\textbf{(3) Conclusion.} \quad
By replacing $\kappa(\delta_\mathrm{dom})$ with a larger one  
according to Remark~\ref{rem-smaller},  
\eqref{e51-twss} is defined on
$(Z, W) \in [-\delta_\mathrm{dom}', \delta_\mathrm{dom}']^2$,  
and its graph represents a part of $T_0^{-k}(W^{s\sigma}) \subset W^s(O)$  
for any $k \in \mathbb{Z}_{> \kappa(\delta_\mathrm{dom})}$.  
From \eqref{e32-ddompD},  
the choice of $y_0^{-\sigma}$ depends only on $\delta_\mathrm{dom}$,  
so by Proposition~\ref{p34-Qk}, we have  
\begin{align*}
  y^{--} - y^- < Y_Q < y^{-+} - y^-  
\end{align*}  
by replacing $\delta_\mathrm{prm}$
and $\kappa(\delta_\mathrm{dom})$ with smaller and larger ones
according to the same remark.
Therefore, from \eqref{e51-twssEV},  
by further replacing $\kappa(\delta_\mathrm{dom})$,  
the parts of $T_0^{-k}(W^{s+})$ and $T_0^{-k}(W^{s-})$  
represented by \eqref{e51-twss}  
become surfaces nearly parallel to the $(Z, W)$-plane  
contained in $\pi_k' \cap \{Y > Y_Q\}$ and $\pi_k' \cap \{Y < Y_Q\}$, respectively,
where the $\pi_k'$ is defined by \eqref{e51-pikp}.
Thus, Lemma~\ref{l51-eumfd} yields that
$W^s(O)$ and $\widetilde{W}^u(Q_k)$  
have a transverse intersection
for any $k \in \mathbb{Z}_{> \kappa(\delta_\mathrm{dom})} \cap 2\mathbb{Z}$
and $(t, \omega) \in \mathcal{R}_k^\mathrm{rep}$.
This completes the proof.  
\end{proof}

\subsection{Homoclinic intersection  
between center and stable manifolds}\label{s52-Hint}  
In this section, we find a homoclinic point of $Q_k$.  
That is, we prove the following.  

\begin{proposition}[Existence of Hopf-homoclinic cycle] \label{p52-HHC}
  Assume that $(f, \Gamma)$ holds the accompanying condition \textbf{(AC)}.
  Then there exists an infinite subset $\mathcal{K} \subset
  \mathbb{Z}_{> \kappa(\delta_\mathrm{dom})} \cap 2\mathbb{Z}$
  associated with sequences $\{(t_k, \omega_k)\}$ in $\mathcal{R}_k^\mathrm{rep}$
  such that
  \begin{itemize}
    \item $\omega_k$ converges to $\omega^*$ as $k \to \infty$ and
    \item $\widetilde{W}^s(Q_k) \cap \widetilde{W}^u(Q_k) \neq \emptyset$
      at $\mu = \mu_k(t_k, \omega_k, \rho_k)$,
      $\omega = \omega_k$, and $\rho = \rho_k(t_k, \omega_k)$
      for any $k \in \mathcal{K}$.
  \end{itemize}
  Moreover, if $(f, \Gamma)$ holds the expanding condition \textbf{(EC)},
  the above $\{(t_k, \omega_k)\}$ can be chosen so that
  $(t_k, \omega_k) \in R_k^\mathrm{ex}$ for any $k \in \mathcal{K}$,
  where $R_k^\mathrm{ex}$ is the set in \eqref{e34-Rkex}.
\end{proposition}

\begin{remark}
  This proof is based on the argument presented in \cite{LLST2022}.  
  For more detailed results and rigorous arguments,  
  the reader is referred to the cited work.  
\end{remark}

\begin{proof}
We divide the proof into several parts.

\medskip

\textbf{(1) Equation of a segment of $\widetilde{W}^u(Q_k)$.}\quad
As shown in Proposition~\ref{p51-Tint},  
since $(f, \Gamma)$ satisfies \textbf{(AC)},  
$\widetilde{W}^u(Q_k)$ intersects $W^s(O)$ transversely  
at some point $M^{tv}_k(t, \omega)$.
Let $\hat{W}^{u*}_k(t, \omega)$ be a one-dimensional small open disk  
in $\widetilde{W}^u(Q_k)$ that contains $M^{tv}_k(t, \omega)$.  
Since $M^{tv}_k(t, \omega) \in W^s(O)$,  
there exists a large $I_k > 0$ such that  
$f^{I_k \mathrm{per}(O)}(M^{tv}_k(t, \omega)) \in W^s_\mathrm{loc}(O)$
and $f^{i \mathrm{per}(O)}(M^{tv}_k(t, \omega)) \notin W^s_\mathrm{loc}(O)$
for any $i \in \{0, 1, \cdots, I_k - 1\}$,  
where $\mathrm{per}(O)$ denotes the period of $O$.  
We write  
\begin{align*}
  f^{I_k \mathrm{per}(O)}(M^{tv}_k(t, \omega))
    = (x_{1, k}^*(t, \omega), x_{2, k}^*(t, \omega), 0)
\end{align*}
in the $(x_1, x_2, y)$ coordinates defined in Section~\ref{s31-FRM}.  
Let $W^{u*}_k(t, \omega)$ be a small neighborhood of  
$f^{I_k \mathrm{per}(O)}(M^{tv}_k(t, \omega))$
in $f^{I_k \mathrm{per}(O)}(\hat{W}^{u*}_k(t, \omega))$.  
Then $W^{u*}_k(t, \omega)$ is described by the following equation:  
\begin{align}
  (x_1, x_2) = (x_{1, k}^*(t, \omega), x_{2, k}^*(t, \omega)) + (O(y), O(y)),
  \label{e52-x12}
\end{align}
where the above $O(y)$ are at least $k$-dependent $C^1$ functions of $(y, t, \omega)$.

Substituting \eqref{e52-x12} into the $y$ equation in \eqref{eq-03-T0kprm},
we get the equation
\begin{align*}
  &y = H_j(\tilde{y}, t, \omega, y), \quad
  H_j := \gamma^{-j} \tilde{y}
    + \hat{\gamma}^{-j} q_j^{(3)}
      (x_{1, k}^*(t, \omega) + O(y), x_{2, k}^*(t, \omega) + O(y),
      \tilde{y}, \varepsilon),
\end{align*}
of $(\tilde{y}, t, \omega, y)$.
Since
\begin{align*}
  H_j = O(\gamma^{-j}), \quad
  H_{j, y} = O(\hat{\gamma}^{-j}), \quad
  H_{j, \tilde{y}} = O(\gamma^{-j})
\end{align*}
Proposition~\ref{pc1-sgl} and Proposition~\ref{sC2-ev} give the solution
\begin{align*}
  y = O(\gamma^{-j}),
\end{align*}
where $O(\gamma^{-j})$ is a $(k, j)$-dependent at least $C^1$ function
of $(\tilde{y}, t, \omega)$
and its first partial derivative with respect to $\tilde{y}$
is also $O(\gamma^{-j})$.
Substituting it into \eqref{e52-x12}
and \eqref{e52-x12} into $\tilde{x}_1$ and $\tilde{x}_2$ equations
in \eqref{eq-03-T0kprm},
the image $S_{k, j}(t, \omega) := T_0^j(W^{u*}_k(t, \omega))$ is given by  
\begin{align}
\begin{aligned}
  \tilde{x}_1 
  &= \lambda^j (x_{1, k}^*(t, \omega) \cos (j\omega)
    - x_{2, k}^*(t, \omega) \sin (j\omega)) 
  + O(\hat{\lambda}^j), \\
  \tilde{x}_2 
  &= \lambda^j (x_{1, k}^*(t, \omega) \sin (j\omega)
    + x_{2, k}^*(t, \omega) \cos (j\omega)) 
  + O(\hat{\lambda}^j),
\end{aligned} \label{e52-tx12}
\end{align}
where the above $O(\hat{\lambda}^j)$ are $(k, j)$-dependent at least $C^1$ functions
of $(\tilde{y}, t, \omega)$
with $\tilde{y} - y^- \in [-\delta_\mathrm{dom}, \delta_\mathrm{dom}]$
and their first partial derivatives with respect to $\tilde{y}$
are also $O(\hat{\lambda}^j)$.

\medskip

\textbf{(2) Pullback to the $(Z, Y, W)$ space.} \quad
By substituting \eqref{e52-tx12} into \eqref{eq-03-T1prm},
the image $T_1(S_j)$ satisfies
\begin{align}
\begin{aligned}
  x_1 - x_1^+ 
  &= \lambda^j A_{1, k}(t, \omega) + b(\tilde{y} - y^-)
     + O((\tilde{y} - y^-)^2) + O(\hat{\lambda}^j), \\
  x_2 - x_2^+ 
  &= \lambda^j A_{2, k}(t, \omega) + O((\tilde{y} - y^-)^2) + O(\hat{\lambda}^j), \\
  y 
  &= \mu_k + \lambda^j A_{3, k}(t, \omega) + d(\tilde{y} - y^-)^2 
     + O((\tilde{y} - y^-)^3) + O(\hat{\lambda}^j),
\end{aligned} \label{e52-x12y}
\end{align}
where 
\begin{align}
\begin{aligned}
  A_{k, i}(t, \omega) 
  &= a_{i1}(x_{1, k}^* \cos (j\omega) - x_{2, k}^* \sin (j\omega)) 
     + a_{i2}(x_{1, k}^* \sin (j\omega) + x_{2, k}^* \cos (j\omega)), 
     \quad i = 1, 2, \\
  A_{3, k}(t, \omega)
  &= c_1(x_{1, k}^* \cos (j\omega) - x_{2, k}^* \sin (j\omega)) 
     + c_2(x_{1, k}^* \sin (j\omega) + x_{2, k}^* \cos (j\omega)),
\end{aligned} \label{e52-A123}
\end{align}
the $O(\hat{\lambda}^j)$ terms have the same property as \eqref{e52-tx12},
$O((\tilde{y} - y^-)^i)$, $i \in \{2, 3\}$, are $(k, j)$-dependent
at least $C^1$ functions of $(\tilde{y}, t, \omega)$,
and their first partial derivatives with respect to $\tilde{y}$
are $O((\tilde{y} - y^-)^{i - 1})$.

Apply \eqref{eq-03-Shilnikov} and express \eqref{e52-x12y}  
using the Shilnikov coordinates on $\Pi_k$.
Substituting \eqref{e52-x12y} into the $y$ equation in \eqref{eq-03-T0kprm},
we obtain
\begin{align*}
  X_1 
  &= \lambda^j A_{1, k} + b s + O(s^2) + O(\hat{\lambda}^j), \\
  X_2 
  &= \lambda^j A_{2, k} + O(s^2) + O(\hat{\lambda}^j), \\
  Y 
  &= \gamma^k \mu_k - y^-
  + \lambda^j \gamma^k A_{3, k} + \gamma^k d s^2 + \gamma^k O(s^3)
  + O(\hat{\lambda}^j \gamma^k),
\end{align*}
where we put $s := \tilde{y} - y^-$ and
the $O(\cdot)$ terms have the same property as \eqref{e52-x12y}.

Applying the remaining coordinate transformation
\eqref{e32-XYshift} and \eqref{e32-ZW},
we have
\begin{align}
\begin{aligned}
  Z 
  &= \lambda^j (\alpha^* A_{1, k} + \beta^* A_{2, k})
     - E_k s + O(s^2) + O(\hat{\lambda}^j), \\
  Y 
  &= \gamma^k \mu_k - y^-
  + \lambda^j \gamma^k A_{3, k} + \gamma^k d s^2 + \gamma^k O(s^3)
  + \gamma^k O(\hat{\lambda}^j), \\
  W 
  &= O(s^2) + O(\lambda^j),
\end{aligned} \label{e52-ZYW}
\end{align}
where the $O(\lambda^j)$ is a $(k, j)$-dependent
at least $C^1$ function of $(s, t, \omega)$
and its first partial derivative with respect to $s$ is $O(\hat{\lambda}^k)$,
the other $O(\cdot)$ terms have the same property as \eqref{e52-x12y},
and the quantity $E_k$ is defined by \eqref{eq-03-Ek}.

\medskip

\textbf{(3) Non-transverse intersection.} \quad
By the same argument as in \cite[Lemma~4.3]{LLST2022} and its proof,
the stable manifold of $Q_k$ is given by  
\begin{align}
  Z = Z_Q + Z^s(W, t, \omega), \quad
  Y = Y_Q + Y^s(W, t, \omega),
  \label{e52-ZY}
\end{align}
where $Z^s$ and $Y^s$ are $C^r$ with respect to $W$,  
and $C^{r - 2}$ with respect to the parameters, satisfying  
\begin{align*}
  Z^s_W = O(\lambda^{-k} \hat{\lambda}^k), \quad
  Z^s, \,
  Y^s, \,
  Y^s_W
    = O(\lambda^k).
\end{align*}
We solve the system of equations \eqref{e52-ZYW} and \eqref{e52-ZY}  
to find a homoclinic point of $Q_k$.  
Substituting the $W$ equation in \eqref{e52-ZYW} into \eqref{e52-ZY},  
we obtain  
\begin{align}
  Z = Z_Q + Z^s(O(s^2) + O(\lambda^j), t, \omega), \quad
  Y = Y_Q + Y^s(O(s^2) + O(\lambda^j), t, \omega).
  \label{e52-ZY2}
\end{align}
Substituting the $Z$ equation \eqref{e52-ZYW} into the $Z$ equation \eqref{e52-ZY2},  
we obtain the following equation in $(s, t, \omega)$:  
\begin{align*}
  s = H_j(t, \omega, s), \quad
  H_j &:= E_k^{-1} (
    \lambda^j (\alpha^* A_{1, k} + \beta^* A_{2, k})
    - Z_Q \\
    &\quad - Z^s(O(s^2) + O(\lambda^j), t, \omega)
    + O(s^2) + O(\hat{\lambda}^j)
  ).
\end{align*}
By Proposition~\ref{p34-Qk}, for $|s| \leq O(\lambda^j)$, we have  
\begin{align*}
  H_j = O(\lambda^j), \quad
  H_{j, s} = O(\lambda^j)
\end{align*}
and thus Proposition~\ref{pc1-sgl} gives the solution  
\begin{align*}
  s = O(\lambda^j),
\end{align*}
where $O(\lambda^j)$ is a $(k, j)$-dependent
at least $C^1$ function of $(t, \omega)$.  
Note that $\mu_k - \gamma^{-k} y^- = O(\lambda^k)$
since $\mu_k$ is the solution of the system \eqref{eq-03-Fsystem1}.
Substituting this into the $Y$ equation in \eqref{e52-ZYW} and \eqref{e52-ZY2},  
and comparing the $Y$ values,
we obtain the following equation in $(s, t, \omega)$:  
\begin{align}
  A_{3, k} = O(\lambda^{-j} \lambda^k),
  \label{e52-LEq}
\end{align}
where $O(\lambda^{-j} \lambda^k)$ is a $(k, j)$-dependent
at least $C^1$ function of $(t, \omega)$. 
By the definition of $A_{3, k}$ in \eqref{e52-A123},
we can rewrite
\begin{align}
  A_{3, k} = A_{3, k}^*\sin (j \omega + \varphi^*), \quad
  A_{3, k}^*= A_{3, k}^*(t, \omega)
    := \sqrt{(c_1^2 + c_2^2)((x_{1, k}^*)^2 + (x_{2, k}^*)^2)},
  \label{e52-A3R}
\end{align}
where $\varphi^* = \varphi^*(t, \omega)$ is the angle determined by
\begin{align*}
  \varphi^*(t, \omega) = \arctan_2(
    c_2 x_{1, k}^* - c_1 x_{2, k}^*, c_1 x_{1, k}^* + c_2 x_{2, k}^*
  ),
\end{align*}
where $\arctan_2$ is a function defined by \eqref{e31-at2}.
By the definition of $(x_{1, k}^*, x_{2, k}^*)$ in the Step~(1)
and the note after \eqref{e32-sinE},
there exists a constant $C = C(\mathbb{F}) > 0$ such that $A_{3, k}^* \geq C$.

Fix $t$ by  
\begin{align*}
  t = t_k(\omega) := \frac{t_k^{+0}(\omega) + t_k^+(\omega)}{2}
\end{align*}
and consider varying only $\omega$.  
Referring to \eqref{e32-Ikbd}, define  
\begin{align}
  \Phi_\mathrm{bd} := \{ \varphi \in \mathbb{R} \:|\:
    |\sin(\varphi + \eta^*(0, \omega^*, 0))| > 2e_\mathrm{bd} \}.
  \label{e52-Pbd}
\end{align}
Take a constant $N = N(\mathbb{F}) \in \mathbb{Z}_{> 0}$ such that  
\begin{align*}
  \{a + i((2\pi)/N) \:|\: i \in \mathbb{Z}\} \cap \Phi_\mathrm{bd} \neq \emptyset
\end{align*}
holds for any $a \in \mathbb{R}$.  
For each $j$, define the value of $k$ by  
\begin{align*}
  k = k_j := 2(N + 1)\left\lfloor j/N \right\rfloor,
\end{align*}
where $\lfloor \cdot \rfloor$ denotes the floor function.  
Now define $\omega_j^*$ by  
\begin{align}
  \omega_j^* := j^{-1}(n_j \pi - \varphi^*(t_k(\omega^*), \omega^*)), \quad
  n_j := \left\lfloor (j\omega^*)/\pi \right\rfloor + i_j,
  \label{e52-omgj}
\end{align}
where $i_j$ is an integer with $0 < i_j \leq N$.  
In fact, we can choose $i_j$ so that $k_j \omega_j^* \in \Phi_\mathrm{bd}$.  
In fact, when $i_j$ increases by $1$, the increment of $k_j \omega_j^*$  
is at most $(2\pi)/N \mod 2\pi$.  
Therefore, such an $i_j$ can be chosen so that $k_j \omega_j^* \in \Phi_\mathrm{bd}$,  
and hence $\omega_j^* \in I_{k_j}^\mathrm{bd}$.  

Now, introduce a new parameter $\Delta\omega$ near $0$  
such that $\omega_j^* + \Delta\omega \in I_{k_j}^\mathrm{bd}$.  
Then, by the definition of $I_{k_j}^\mathrm{bd}$,  
we have $|\Delta\omega| \leq O(k_j^{-1})$.  
Next, equation \eqref{e52-LEq} can be rewritten using \eqref{e52-A3R} as  
\begin{align*}
  \sin(j \omega + \varphi^*(t_{k_j}(\omega), \omega))
  = O(\lambda^{-j} \lambda^{k_j}), \quad
  \omega = \omega_j^* + \Delta\omega.
\end{align*}
Note that $\omega_j^*$ converges to $\omega^*$ as $j \to \infty$ by \eqref{e52-omgj}.
Thus, the above equation becomes  
\begin{align}
  O(\Delta\omega) + O(\lambda^{-j} \lambda^{k_j}) = 0
  \label{e52-lsteq}
\end{align}
Now, as we vary $\Delta\omega$ from its minimum to maximum allowed value,  
if $j$ is sufficiently large,  
the left-hand side of \eqref{e52-lsteq} changes sign.  
By the intermediate value theorem,  
there exists a solution $\omega = \omega_{k_j} = \omega_j^* + \Delta\omega_{k_j}^*$
to \eqref{e52-lsteq}.  
Letting $t_{k_j} = t_{k_j}(\omega_{k_j})$,  
we obtain $(t_{k_j}, \omega_{k_j}) \in \mathcal{R}_k^\mathrm{rep}$,  
and $\omega_{k_j}$ converges to $\omega^*$.
This completes the proof of the first part of the proposition. 

\medskip

\textbf{(4) For the case of (EC).} \quad
Next, we consider the case where $(f, \Gamma)$ satisfies \textbf{(EC)}.  
The proof proceeds in exactly the same way as above.  
Instead of using $\Phi_\mathrm{bd}$ in \eqref{e52-Pbd},  
we define, referring to \eqref{e34-Ikex},  
\begin{align*}
  \Phi_\mathrm{ex} := \{ \varphi \in \mathbb{R} \:|\:
    \sin(\varphi + \eta^*(0, \omega^*, 0)) + 1 < \delta'/2 \}.
\end{align*}
Then, we reselect $N$ accordingly.  
This completes the proof.
\end{proof}

Now we are ready to prove the third theorem.

\begin{proof}[Proof of Theorem~\ref{thm-third}]
Since $(f, \Gamma)$ holds \textbf{(AC)},  
there exist $\mathcal{K}$ and $(t_k, \omega_k) \in \mathcal{R}_k^\mathrm{rep}$,  
$k \in \mathcal{K}$,  
in Proposition~\ref{s52-Hint}.  
Let  
\begin{align*}
  \varepsilon_k := (\mu_k(t_k, \omega_k, \rho_k), \omega_k, \rho_k(t_k, \omega_k))
\end{align*}
for any $k \in \mathcal{K}$,  
where $\mu_k$ and $\rho_k$ are defined in Proposition~\ref{p34-Qk}.  
By the definition of $\mathcal{R}_k^\mathrm{rep}$ and Proposition~\ref{s52-Hint},  
$Q_k$ is a generic Hopf point with a negative Lyapunov coefficient  
and has a Hopf-homoclinic cycle.  
In addition,  
Propositions~\ref{s52-Hint} and \ref{p34-Qk} imply that  
$\varepsilon_k$ converges to $(0, \omega^*, 0)$ as $k \to \infty$.  
This completes the proof of the first case.  

The case where $(f, \Gamma)$ satisfies \textbf{(EC)} can be proved in a similar way,  
because Proposition~\ref{p34-Qk} yields $\rho_k(t_k, \omega_k) < 0$  
for any $k \in \mathcal{K}$.  
This completes the proof.
\end{proof}

\ifdefined\isMaster
\else
  \end{document}
\fi
\ifdefined\isMaster
\else
  \documentclass[11pt,a4paper]{article}

\usepackage[left=2cm,right=2cm,top=2cm,bottom=2cm]{geometry}
\usepackage{graphicx}
\usepackage{multirow}
\usepackage{amsmath,amssymb,amsfonts}
\usepackage{amsthm}
\usepackage{mathrsfs}
\usepackage[title]{appendix}
\usepackage{xcolor}
\usepackage{textcomp}
\usepackage{manyfoot}
\usepackage{booktabs}
\usepackage{algorithm}
\usepackage{algorithmicx}
\usepackage{algpseudocode}
\usepackage{listings}
\ifdefined\debug
  \usepackage[inline]{showlabels}
  \showlabels{ref}
  \showlabels{cite}
  \usepackage{fancyhdr}
\fi
\usepackage{bm}
\usepackage{hyperref}

\theoremstyle{plain}
\newtheorem{theorem}{Theorem}[section]
\newtheorem{proposition}[theorem]{Proposition}%
\newtheorem{lemma}[theorem]{Lemma}
\newtheorem{corollary}[theorem]{Corollary}
\newtheorem{claim}[theorem]{Claim}
\theoremstyle{definition}
\newtheorem{definition}[theorem]{Definition}
\newtheorem{remark}[theorem]{Remark}
\newtheorem{example}[theorem]{Example}
\newtheorem{question}[theorem]{Question}
\newtheorem*{notation}{Notation}
\theoremstyle{plain}
\newtheorem{maintheorem}{Theorem}
\renewcommand{\themaintheorem}{\Alph{maintheorem}}
\newcommand{\wiptext}[1]{\textcolor{blue}{#1}}
\ifdefined\debug
  \newcommand{\debugtext}[1]{\textcolor{red}{#1}}
  \pagestyle{fancy} 
  \fancyhead[L]{}
  \fancyhead[C]{}
  \fancyhead[R]{\textcolor{red}{[DEBUG] Debug mode is activated.}}
  \renewcommand{\headrulewidth}{0pt}
\else
  \newcommand{\debugtext}[1]{\ignorespaces}
\fi

\counterwithin{figure}{section}
\counterwithin{table}{section}

\raggedbottom

\numberwithin{equation}{section}

  \begin{document}
\fi

\begin{appendices}

\section{Toy model on 3-sphere satisfying expanding condition} \label{sA-toy}
In this appendix, we construct a concrete $C^r$, $r \geq 1$,
diffeomorphism $f$ on 3-sphere $S^3$ satisfying
the assumptions of Theorem~\ref{thm-main} and the expanding
condition \textbf{(EC)}.
Hence, Question~\ref{q-CQ} is resolved due to the existence of such a system.

\medskip

\begin{figure}[h]
  \centering
  \includegraphics{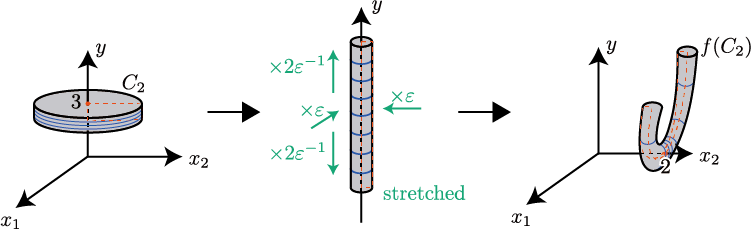}
  \caption{
    The transformation of $C_2$ under $f$. After being linearly 
    stretched, it is further modified by nonlinear transformations such as rotation 
    and bending, resulting in the configuration shown in the rightmost diagram.
  } \label{fig-toymodel}
\end{figure}

We define
\begin{align*}
  C := \{(x_1, x_2, y) \:|\: x_1^2 + x_2^2 \leq 3^2 ,\, 0 \leq y \leq 3\}, \quad
  C_1 := C \cap \{0 \leq y \leq 1\}, \quad
  C_2 := C \cap \{2 \leq y \leq 3\}.
\end{align*}
The $f|_{C_1}: (x_1, x_2, y) \mapsto (\hat{x}_1, \hat{x}_2, \hat{y})$
and $f|_{C_2}: (x_1, x_2, y) \mapsto (\bar{x}_1, \bar{x}_2, \bar{y})$
are assumed to be given as follows:
\begin{align*}
  \hat{x}_1 = \frac{x_1}{3} \cos\frac{\pi}{6} - \frac{x_2}{3} \sin\frac{\pi}{6}, \quad
  \hat{x}_2 = \frac{x_1}{3} \sin\frac{\pi}{6} + \frac{x_2}{3} \cos\frac{\pi}{6}, \quad
  \hat{y} = 3y,
\end{align*}
and
\begin{align*}
  \bar{x}_1 = 2\varepsilon^{-1}(y - 2.5) - 4\varepsilon^{-2} (y - 2.5)^2, \quad
  \bar{x}_2 = -\varepsilon x_2 + 2, \quad
  \bar{y} = \varepsilon x_1 + 4\varepsilon^{-2} (y - 2.5)^2,
\end{align*}
where $\varepsilon > 0$ is a small number.
$f|_{C_1}$ is a simple linear map, and the image of $C_2$ under $f$ 
is deformed as shown in Figure~\ref{fig-toymodel}.
The projection of $f(C_1)$ under $\mathrm{pr}_2(x_1, x_2, y) = x_2$ has the 
image $\{-1 \leq x_2 \leq 1\}$, while the image of $f(C_2)$ is 
$\{-3\varepsilon + 2 \leq x_2 \leq 3\varepsilon + 2\}$.
Hence, for $\varepsilon$ with $0 < \varepsilon < 1/3$,
$f(C_1)$ and $f(C_2)$ are disjoint.
Fix $\varepsilon \in (0, 1/3)$, and choose a sufficiently large open ball 
$B \subset \mathbb{R}^3$ centered at the origin such that $f(C_1), f(C_2) 
\subset B$.  
Extend the domain of $f$ to $C$ so that $f(C) \subset B$
and $f$ remains injective.  
Further extend the domain of $f$ to $B$ so that $f(B) \subset B$
and $f$ remains injective.  
Finally, by adding the point $\infty$ to $\mathbb{R}^3$ and identifying 
it with the 3-sphere $S^3$, extend $f$ to a $C^r$ diffeomorphism $f: S^3 \to S^3$ such
that $f$ has the source $\infty$ with $S^3 \setminus f(B) \subset W^u(\infty)$.

Note that the origin $O^*$ is a hyperbolic fixed point of $f$
and its multipliers are $\frac{1}{3}(\frac{\sqrt{3}}{2} \pm \mathrm{i} \frac{1}{2})$
and $3$.
In particular,
$\left|\frac{1}{3}(\frac{\sqrt{3}}{2} \pm \mathrm{i} \frac{1}{2})\right| \cdot 3 = 1$.
Also, note that the segment $\ell^u := C_2 \cap \{x_1 = 0, \, x_2 = 0\}$
is contained $W^u(O^*)$ by the definition of $f|_{C_1}$.
The image $f(\ell^u) \subset W^u(O^*)$ is given by
\begin{align*}
  \{(2\varepsilon^{-1}t - 4\varepsilon^{-2}t^2, \, 2, \, 4\varepsilon^{-2}t^2) \:|\: 
    -0.5 \leq t \leq 0.5\},
\end{align*}
which has the tangency $M^+_0 = (0, 2, 0)$ with $\{x_1^2 + x_2^2 \leq 3, \, y = 0\}
\subset W^s(O^*)$.
Thus, $f$ satisfies the assumption of Theorem~\ref{thm-main}.

Finally, let us verify that $(f, \Gamma)$ satisfies
the expanding condition \textbf{(EC)},
where $\Gamma$ is an orbit of $M^+_0$.
Let $U_0$ be a small neighborhood of
$C_1 \cup \{x_1 = 0, \, x_2 = 0, \, 0 \leq y \leq 3\}$,
let $\mathbb{U}_0^*$ be a pair of $U_0$ and the coordinates $(x_1, x_2, y)$,
and let $M^-_0 := f^{-1}(M^+_0) = (0, 0, 2.5)$.
Recall the quantity
$\mathcal{E}(f, \Gamma, \mathbb{U}_0^*, M^-_0, M^+_0)$
defined in \eqref{eq-02-mathcalE}.
In our settings,
\begin{align*}
  \mathcal{E}(f, \Gamma, \mathbb{U}_0^*, M^-_0, M^+_0)
    = \sqrt{(2\varepsilon^{-1})^2 + 0} \sqrt{\varepsilon^2 + 0}
    = 2 > 1.
\end{align*}
Hence, $f$ satisfies the \textbf{(EC)}.
From the above, Question~\ref{q-CQ} has been resolved affirmatively.

\ifdefined\isMaster
\else
  \end{appendices}
  \end{document}
\fi
\ifdefined\isMaster
\else
  \documentclass[11pt,a4paper]{article}

\usepackage[left=2cm,right=2cm,top=2cm,bottom=2cm]{geometry}
\usepackage{graphicx}
\usepackage{multirow}
\usepackage{amsmath,amssymb,amsfonts}
\usepackage{amsthm}
\usepackage{mathrsfs}
\usepackage[title]{appendix}
\usepackage{xcolor}
\usepackage{textcomp}
\usepackage{manyfoot}
\usepackage{booktabs}
\usepackage{algorithm}
\usepackage{algorithmicx}
\usepackage{algpseudocode}
\usepackage{listings}
\ifdefined\debug
  \usepackage[inline]{showlabels}
  \showlabels{ref}
  \showlabels{cite}
  \usepackage{fancyhdr}
\fi
\usepackage{bm}
\usepackage{hyperref}

\theoremstyle{plain}
\newtheorem{theorem}{Theorem}[section]
\newtheorem{proposition}[theorem]{Proposition}%
\newtheorem{lemma}[theorem]{Lemma}
\newtheorem{corollary}[theorem]{Corollary}
\newtheorem{claim}[theorem]{Claim}
\theoremstyle{definition}
\newtheorem{definition}[theorem]{Definition}
\newtheorem{remark}[theorem]{Remark}
\newtheorem{example}[theorem]{Example}
\newtheorem{question}[theorem]{Question}
\newtheorem*{notation}{Notation}
\theoremstyle{plain}
\newtheorem{maintheorem}{Theorem}
\renewcommand{\themaintheorem}{\Alph{maintheorem}}
\newcommand{\wiptext}[1]{\textcolor{blue}{#1}}
\ifdefined\debug
  \newcommand{\debugtext}[1]{\textcolor{red}{#1}}
  \pagestyle{fancy} 
  \fancyhead[L]{}
  \fancyhead[C]{}
  \fancyhead[R]{\textcolor{red}{[DEBUG] Debug mode is activated.}}
  \renewcommand{\headrulewidth}{0pt}
\else
  \newcommand{\debugtext}[1]{\ignorespaces}
\fi

\counterwithin{figure}{section}
\counterwithin{table}{section}

\raggedbottom

\numberwithin{equation}{section}

  \begin{document}
  \begin{appendices}
\fi

\section{Proof of Proposition~\ref{p23-IndEC}}
\begin{proof}[Proof of Proposition~\ref{p23-IndEC}]
We had verified that the validity of the expanding condition \textbf{(EC)}
does not depend on the choice of $U_0$, $M_0^-$, and $M_0^+$.
It remains to show that it is also independent of the choice of coordinates.

We take $C^r$ coordinates $(\bm{u}, v)$, $\bm{u} = (u_1, u_2)$ on $U_0$ such that
\begin{align}
  W^s_\mathrm{loc}(O^*) = \{ v = 0 \}, \quad
  W^u_\mathrm{loc}(O^*) = \{ \bm{u} = 0\},
  \label{eq-A2-Wsloc_Wuloc}
\end{align}
and $T_0^*: (\bm{u}, v)
\mapsto (\hat{\bm{u}}, \hat{v})$ has the form
\begin{align}
  \hat{\bm{u}} = \lambda^* \bm{u} R(\omega^*)^\mathsf{T}
    + q_{12}^*(\bm{u}, v), \quad
  \hat{v} = \gamma^* v + q_3^*(\bm{u}, v)
  \label{eq-A2-T00_uv}
\end{align}
where $R(\theta)$ denotes the rotation matrix of angle $\theta$
define in \eqref{eq-03-rotMat}.
Here, $q_{12}^*$ and $q_3^*$ are $C^r$ maps with
\begin{align*}
  q_{12}^*(0, 0) = 0, \quad
  q_3^*(0, 0) = 0, \quad
  \frac{\partial q_{12}^*}{\partial (\bm{u}, v)}(0, 0) = 0, \quad
  \frac{\partial q_3^*}{\partial (\bm{u}, v)}(0, 0) = 0, \quad
  q_{12}^*(0, v) \equiv 0, \quad
  q_3^*(\bm{u}, 0) \equiv 0.
\end{align*}

Let $\bm{s} = (s_1, s_2)$.
Recall that the global map $T_1^*: (\bm{s}, t) \mapsto (\bar{\bm{s}}, \bar{t})$
was given by \eqref{eq-02-ExpT10}.
We put
\begin{align*}
  \bm{A}^* := (a_{ij}^*)_{i, j \in \{1, 2\}}, \quad
  \bm{b}^* := (b_1^*, b_2^*), \quad
  \bm{c}^* := (c_1^*, c_2^*), \quad
  \bm{s}^+ := (s_1^+, s_2^+)
\end{align*}  
and we rewrite the global map as
\begin{align*}
  \bar{\bm{s}} - \bm{s}^+ &= \bm{s} (\bm{A}^*)^\mathsf{T} + \bm{b}^*(t - t^-)  
    + R_{12}(\bm{s}, t), \\
  \bar{t} &= \langle \bm{c}^*, \bm{s} \rangle  
  + R_3(\bm{s}, t),
\end{align*}
where
$\langle \cdot, \cdot \rangle$ is the Euclidean inner product,
and $R_{12}(\bm{s}, t)$ and $R_3(\bm{s}, t)$
are terms of second order or higher of the Taylor expansion,
in other words,
they hold
\begin{align}
\begin{aligned}
  &R_{12}(0, t^-) = 0, \quad
  R_{12, \bm{s}}(0, t^-) = 0, \quad
  R_{12, t}(0, t^-) = 0, \\
  &R_3(0, t^-) = 0, \quad
  R_{3, \bm{s}}(0, t^-) = 0, \quad
  R_{3, t}(0, t^-) = 0.
\end{aligned}\label{eq-A2-R123_prop}
\end{align}
Here, we used a similar notation in \eqref{e31-Fsig} for partial derivatives
to simplify the notation;
for instance, $R_{12, t} = \frac{\partial R_{12}}{\partial t}$
and $R_{3, \bm{s}} = \frac{\partial R_3}{\partial \bm{s}}$
are a $2 \times 1$ matrix and a $1 \times 2$ matrix, respectively.
Using the coordinates $(\bm{u}, v)$,
let us express the global map
$T_1^*: (\bm{u}, v) \mapsto (\bar{\bm{u}}, \bar{v})$:
\begin{align*}
  \bar{\bm{u}} - \bm{u}^+ &= \bm{u} (\bm{A}^*_{new})^\mathsf{T}
    + \bm{b}^*_{new}(v - v^-)  
    + R_{12}^{new}(\bm{u}, v), \\
  \bar{v} &= \langle \bm{c}^*_{new}, \bm{u} \rangle  
  + R_3^{new}(\bm{u}, v),
\end{align*}
where $M_0^- = (0, v^-)$ and $M_0^+ = (\bm{u}^+, 0)$ in the $(\bm{u}, v)$ coordinates,
and $R_{12}^{new}(\bm{u}, v)$ and $R_3^{new}(\bm{u}, v)$
are terms of second order or higher of the Taylor expansion.
We would like to check
$\| \bm{b}^*_{new} \| \| \bm{c}^*_{new} \| = \| \bm{b}^* \| \| \bm{c}^* \|$.

We denote the coordinate transformation and its inverse by  
\begin{align*}
  &\bm{u} = \tau_{12}(\bm{s}, t), \quad
  v = \tau_3(\bm{s}, t), \\
  &\bm{s} = \sigma_{12}(\bm{u}, v), \quad
  t = \sigma_3(\bm{u}, v).
\end{align*}
Then, we have
\begin{align*}
  &\bar{\bm{u}} = \tau_{12}(
    \bm{s}^+ + \sigma_{12} (\bm{A}^*)^\mathsf{T} + \bm{b}^*(\sigma_3 - t^-)
    + R_{12}(\sigma_{12}, \sigma_3),
    \langle \bm{c}^*, \sigma_{12} \rangle  
    + R_3(\sigma_{12}, \sigma_3)
  ), \\
  &\bar{v} = \tau_3(
    \bm{s}^+ + \sigma_{12} (\bm{A}^*)^\mathsf{T} + \bm{b}^*(\sigma_3 - t^-)
    + R_{12}(\sigma_{12}, \sigma_3),
    \langle \bm{c}^*, \sigma_{12} \rangle  
    + R_3(\sigma_{12}, \sigma_3)
  ),
\end{align*}
where $\sigma_{12} = \sigma_{12}(\bm{u}, v)$ and $\sigma_3 = \sigma_3(\bm{u}, v)$.
Thus, by using \eqref{eq-A2-R123_prop}, we get
\begin{align*}
  &(\bm{b}^*_{new})^\mathsf{T} = \frac{\partial \bar{\bm{u}}}{\partial v}(0, v^-)
    = \tau_{12, \bm{s}} \cdot \left(
      \bm{A}^* \sigma_{12, v} + (\bm{b}^*)^\mathsf{T} \sigma_{3, v}
    \right)
    + \tau_{12, t} \cdot \langle \bm{c}^*, (\sigma_{12, v})^\mathsf{T} \rangle, \\
  &\bm{c}^*_{new} = \frac{\partial \bar{v}}{\partial \bm{u}}(0, v^-)
    = \tau_{3, \bm{s}} \cdot \left(
      \bm{A}^* \sigma_{12, \bm{u}} + (\bm{b}^*)^\mathsf{T} \sigma_{3, \bm{u}}
    \right)
    + \tau_{3, t} \cdot \bm{c}^* \sigma_{12, \bm{u}},
\end{align*}
where
\begin{align*}
  &\tau_{12, \bm{s}} = \tau_{12, \bm{s}}(\bm{s}^+, 0), \quad
  \tau_{12, t} = \tau_{12, t}(\bm{s}^+, 0), \quad
  \tau_{3, \bm{s}} = \tau_{3, \bm{s}}(\bm{s}^+, 0), \quad
  \tau_{3, t} = \tau_{3, t}(\bm{s}^+, 0), \\
  &\sigma_{12, \bm{u}} = \sigma_{12, \bm{u}}(0, v^-), \quad
  \sigma_{12, v} = \sigma_{12, v}(0, v^-), \quad
  \sigma_{3, \bm{u}} = \sigma_{3, \bm{u}}(0, v^-), \quad
  \sigma_{3, v} = \sigma_{3, v}(0, v^-).
\end{align*}
In fact, the following hold (proof will be given later):  
\begin{align}
  &\tau_{3, \bm{s}}(\bm{s}^+, 0) = 0, \quad
  \sigma_{12, v}(0, v^-) = 0, \label{eq-A2-factD} \\
  &\sigma_{12, \bm{u}}(0, v^-) = (\tau_{12, \bm{s}}(\bm{s}^+, 0))^{-1}, \quad
  \sigma_{3, v}(0, v^-) = (\tau_{3, t}(\bm{s}^+, 0))^{-1}, \label{eq-A2-factI} \\
  &(\tau_{12, \bm{s}}(\bm{s}^+, 0))^\mathsf{T}
    = d (\tau_{12, \bm{s}}(\bm{s}^+, 0))^{-1}, \quad
  d := \mathrm{det} \, \tau_{12, \bm{s}}(0, 0) \quad (> 0). \label{eq-A2-factO}
\end{align} 
Note that
$\tau_{12, \bm{s}}$ is a $2 \times 2$ matrix and $\tau_{3, t}$ is a real number.
Using the above fact,
we obtain
\begin{align*}
  \bm{b}_{new}^* = \tau_{3, t}^{-1} \bm{b}^* d \tau_{12, \bm{s}}^{-1}, \quad
  (\bm{b}_{new}^*)^\mathsf{T}
    = \tau_{12, \bm{s}} (\bm{b}^*)^\mathsf{T} \tau_{3, t}^{-1}, \quad
  \bm{c}_{new}^* = \tau_{3, t} \bm{c}^* \tau_{12, \bm{s}}^{-1}, \quad
  (\bm{c}_{new}^*)^\mathsf{T}
    = d^{-1} \tau_{12, \bm{s}} (\bm{c}^*)^\mathsf{T} \tau_{3, t}.
\end{align*}
Thus, we get the desired result:
\begin{align*}
  \| \bm{b}^*_{new} \|^2 \| \bm{c}^*_{new} \|^2
  &= \bm{b}^*_{new} (\bm{b}^*_{new})^\mathsf{T}
    \bm{c}^*_{new} (\bm{c}^*_{new})^\mathsf{T} \\
  &= \tau_{3, t}^{-1} \bm{b}^* d \tau_{12, \bm{s}}^{-1}
     \tau_{12, \bm{s}} (\bm{b}^*)^\mathsf{T} \tau_{3, t}^{-1}
     \tau_{3, t} \bm{c}^* \tau_{12, \bm{s}}^{-1}
     d^{-1} \tau_{12, \bm{s}} (\bm{c}^*)^\mathsf{T} \tau_{3, t}
  = \| \bm{b}^* \|^2 \| \bm{c}^* \|^2.
\end{align*}

It remains to prove \eqref{eq-A2-factD} -- \eqref{eq-A2-factO}.
The \eqref{eq-A2-factD} follows from
\eqref{eq-02-CCWsWu} and \eqref{eq-A2-Wsloc_Wuloc}.  
Indeed, it follows from these that
$\tau_{12}(0, t) \equiv 0$,
$\tau_3(\bm{s}, 0) \equiv 0$, 
$\sigma_{12}(0, v) \equiv 0$, 
$\sigma_3(\bm{u}, 0) \equiv 0$, and hence,
\begin{align*}
  \tau_{12, t}(0, t) \equiv 0, \quad
  \tau_{3, \bm{s}}(\bm{s}, 0) \equiv 0, \quad
  \sigma_{12, v}(0, v) \equiv 0, \quad
  \sigma_{3, \bm{u}}(\bm{u}, 0) \equiv 0.
\end{align*}  
Next, let us verify \eqref{eq-A2-factI} and \eqref{eq-A2-factO}.  
First, note that $\tau_{12, t}(0, t^-) = 0$ and
$\sigma_{12, v}(0, v^-) = 0$ imply
\begin{align*}
  \sigma_{12, \bm{u}}(0, v^-) = (\tau_{12, \bm{s}}(0, t^-))^{-1}, \quad
  \sigma_{3, v}(0, v^-) = (\tau_{3, t}(0, t^-))^{-1},
\end{align*}
respectively.
Thus, we need to verify
\begin{align}
  \tau_{12, \bm{s}}(\bm{s}^+, 0) = \tau_{12, \bm{s}}(0, t^-), \quad
  \tau_{3, t}(\bm{s}^+, 0) = \tau_{3, t}(0, t^-),
  \label{eq-A2-goal}
\end{align}
and \eqref{eq-A2-factO}.

Note that the $2 \times 2$ matrix $A := \tau_{12, \bm{s}}(0, 0)$ commutes with  
the rotation matrix $R(\omega^*)$.  
Indeed, since the differential at the origin of the composition of  
\eqref{eq-02-ExpT10} and the coordinate transformation  
$\tau: (\bm{s}, t) \mapsto (\tau_{12}(\bm{s}, t), \tau_3(\bm{s}, t))$  
coincides with the differential at the origin of the composition of  
$\tau$ and \eqref{eq-A2-T00_uv}, we obtain  
\begin{align*}
  D(\tau)_{(0, 0)}  
  \begin{pmatrix}
    R(\omega^*) & 0 \\
    0 & \gamma^*
  \end{pmatrix}
  =
  \begin{pmatrix}
    R(\omega^*) & 0 \\
    0 & \gamma^*
  \end{pmatrix}
  D(\tau)_{(0, 0)},
\end{align*}
which implies $A R(\omega^*) = R(\omega^*) A$.
Thus, since $\omega^* \in (0, \pi)$,
we can write $A = aI + bJ$ for some $a$, $b \in \mathbb{R}$,
where $I = \begin{pmatrix}
  1 & 0 \\
  0 & 1
\end{pmatrix}$
and $J = \begin{pmatrix}
  0 & -1 \\
  1 & 0
\end{pmatrix}$.
Hence, $A^\mathsf{T} = (a^2 + b^2) A^{-1}$
and $A$ commutes with any rotation matrix $R$.
Indeed,
\begin{align*}
  A^\mathsf{T} A = (aI - bJ)(aI + bJ) = (a^2 + b^2)I, \quad
  A R = (aI - bJ) R = R (aI - bJ) = R A.
\end{align*}

Let $(\bm{s}_n, 0) := (T_0^*)^n(\bm{s}^+, 0)$ for any $n \in \mathbb{Z}_{> 0}$
with $(T_0^*)^n(M_0^+) \in U_0$.
It is well-defined for sufficiently large $n \in \mathbb{Z}_{> 0}$.
Since the differential at $(\bm{s}^+, 0)$ of the composition of  
$(T_0^*)^n$ and $\tau$
coincides with the differential at $(\bm{s}^+, 0)$ of the composition of  
$\tau$ and $(T_0^*)^n$, we get
\begin{align*}
  D(\tau)_{(\bm{s}^+, 0)} = 
  \begin{pmatrix}
    (R(\omega^*))^{-n} & 0 \\
    0 & (\gamma^*)^{-n}
  \end{pmatrix}
  D(\tau)_{(\bm{s}_n, 0)}
  \begin{pmatrix}
    (R(\omega^*))^n & 0 \\
    0 & (\gamma^*)^n
  \end{pmatrix}.
\end{align*}
By the compactness of the space of all rotation matrices,  
there exists a subsequence $\{ n_i \}_{i \in \mathbb{Z}_{> 0}}$ such that  
$(R(\omega^*))^{n_i}$ converges to some rotation matrix $R$.  
Taking $n = n_i$ in the above equation and letting $i \to \infty$,  
we obtain the following since $\tau$ is at least $C^1$:  
\begin{align}
  \tau_{12, \bm{s}}(\bm{s}^+, 0) = R^{-1} A R = A, \quad
  \tau_{3, t}(\bm{s}^+, 0) = \tau_{3, t}(0, 0).
  \label{eq-A2-base_spz}
\end{align}  
The \eqref{eq-A2-factO} have been proven.

By repeating a similar argument for the sequence  
$(0, t_n) := (T_0^*)^{-n}(0, t^-)$ ($n \in \mathbb{Z}_{> 0}$),  
we obtain  
\begin{align*}
  \tau_{12, \bm{s}}(0, t^-) = A, \quad
  \tau_{3, t}(0, t^-) = \tau_{3, t}(0, 0).
\end{align*}  
Combining this with \eqref{eq-A2-base_spz}, \eqref{eq-A2-goal} is proven.  
We complete the proof.
\end{proof}

\ifdefined\isMaster
\else
  \end{appendices}
  \end{document}
\fi
\ifdefined\isMaster
\else
  \documentclass[11pt,a4paper]{article}

\usepackage[left=2cm,right=2cm,top=2cm,bottom=2cm]{geometry}
\usepackage{graphicx}
\usepackage{multirow}
\usepackage{amsmath,amssymb,amsfonts}
\usepackage{amsthm}
\usepackage{mathrsfs}
\usepackage[title]{appendix}
\usepackage{xcolor}
\usepackage{textcomp}
\usepackage{manyfoot}
\usepackage{booktabs}
\usepackage{algorithm}
\usepackage{algorithmicx}
\usepackage{algpseudocode}
\usepackage{listings}
\ifdefined\debug
  \usepackage[inline]{showlabels}
  \showlabels{ref}
  \showlabels{cite}
  \usepackage{fancyhdr}
\fi
\usepackage{bm}
\usepackage{hyperref}

\theoremstyle{plain}
\newtheorem{theorem}{Theorem}[section]
\newtheorem{proposition}[theorem]{Proposition}%
\newtheorem{lemma}[theorem]{Lemma}
\newtheorem{corollary}[theorem]{Corollary}
\newtheorem{claim}[theorem]{Claim}
\theoremstyle{definition}
\newtheorem{definition}[theorem]{Definition}
\newtheorem{remark}[theorem]{Remark}
\newtheorem{example}[theorem]{Example}
\newtheorem{question}[theorem]{Question}
\newtheorem*{notation}{Notation}
\theoremstyle{plain}
\newtheorem{maintheorem}{Theorem}
\renewcommand{\themaintheorem}{\Alph{maintheorem}}
\newcommand{\wiptext}[1]{\textcolor{blue}{#1}}
\ifdefined\debug
  \newcommand{\debugtext}[1]{\textcolor{red}{#1}}
  \pagestyle{fancy} 
  \fancyhead[L]{}
  \fancyhead[C]{}
  \fancyhead[R]{\textcolor{red}{[DEBUG] Debug mode is activated.}}
  \renewcommand{\headrulewidth}{0pt}
\else
  \newcommand{\debugtext}[1]{\ignorespaces}
\fi

\counterwithin{figure}{section}
\counterwithin{table}{section}

\raggedbottom

\numberwithin{equation}{section}

  \begin{document}
  \begin{appendices}
\fi

\section{System of equations}
We often encounter situations where we need to solve a system of equations  
and estimate the partial derivatives of its solution.  
In this appendix,  
we first explain the method for solving a single  
equation in Section~\ref{sC1-sgl}, see Proposition~\ref{pc1-sgl}.  
Next, in Section~\ref{sC2-ev}, we describe how to estimate the partial  
derivatives of the solution, see Proposition~\ref{pc2-ev}.  
Finally, in Section~\ref{sC3-sys}, we discuss the application of these  
methods to solve a system of equations and estimate the partial derivatives,
see Proposition~\ref{pc3-sys}.  

There is no relationship between the symbols that appear in this appendix
and those that appear in the other sections.

\subsection{Single equation}
\label{sC1-sgl}
In this section,
we explain how to solve a single equation.

\medskip

Let $\{ G_k: U \to \mathbb{R} \}_{k = 1}^\infty$
be a sequence of $C^r$,
$r \in \mathbb{Z}_{> 0} \cup \{\infty, \omega\}$,
functions from an open set $U \subset \mathbb{R}^n$,
$n \in \mathbb{Z}_{> 0}$,
to $\mathbb{R}$.
Let $\{ H_k: U \times \mathbb{R} \to \mathbb{R} \}_{k = 1}^\infty$
be a sequence of $C^r$ functions.
For the above core objects, we set
\begin{align*}
  \mathbb{G}_1 := (\{ G_k: U \to \mathbb{R} \}_{k = 1}^\infty,
  \{ H_k: U \times \mathbb{R} \to \mathbb{R} \}_{k = 1}^\infty).
\end{align*}
Let $U \times \mathbb{R}$ has the coordinates
$(\bm{x}, y)$, where $\bm{x} = (x_1, x_2, \cdots, x_n)$.
In the following,
we use the notation in \eqref{e31-Fsig} for partial derivatives.

\begin{proposition}[Solution method for a single equation] \label{pc1-sgl}
  Assume
  \begin{gather*}
    \mathcal{H}_k = \mathcal{H}_k(\mathbb{G}_1)
      := \sup_{\bm{x} \in U, \, y \in \mathbb{R}} |H_k(\bm{x}, y)|
    \to 0 \quad \text{as} \quad k \to \infty, \\
    \sup_{\bm{x} \in U, \, y \in \mathbb{R}} |H_{k, y}(\bm{x}, y)|
    \to 0 \quad \text{as} \quad k \to \infty.
  \end{gather*}
  Then, there exists $\kappa = \kappa(\mathbb{G}_1) > 0$ such that
  the equation of $(\bm{x}, y) \in U \times \mathbb{R}$
  \begin{align*}
    y = G_k(\bm{x}) + H_k(\bm{x}, y)
  \end{align*}
  has the solution
  \begin{align*}
    y = G_k(\bm{x}) + I_k(\bm{x})
  \end{align*}
  for any $k > \kappa$,
  where $I_k: U \to \mathbb{R}$
  are $C^r$ functions such that
  there exists a constant $C = C(\mathbb{G}_1) > 0$ satisfying
  \begin{align*}
    |I_k(\bm{x})| \leq C \mathcal{H}_k
  \end{align*}
  for any $\bm{x} \in U$ and $k > \kappa$.
\end{proposition}

\begin{proof}[Proof of Proposition~\ref{pc1-sgl}]
We define $F_k(\bm{x}, y) := y - G_k(\bm{x}) - H_k(\bm{x}, y)$
for any $(\bm{x}, y) \in U \times \mathbb{R}$.
By the definition of $F_k$, we have
\begin{align}
  F_k(\bm{x}, G_k(\bm{x}) + \Delta y)
    = \Delta y - H_k(\bm{x}, G_k(\bm{x}) + \Delta y)
  \label{eq-A3-Fkexp}
\end{align}
for any $\bm{x} \in U$, $\Delta y \in \mathbb{R}$,
and $k \in \mathbb{Z}_{> 0}$.
Differentiating the above equation with respect to $\Delta y$, we obtain  
\begin{align*}
  \partial_{\Delta y} F_k(\bm{x}, G_k(\bm{x}) + \Delta y)  
    = 1 - H_{k, y}(\bm{x}, G_k(\bm{x}) + \Delta y).  
\end{align*}  
Since the assumptions of the lemma hold,
there exists $\kappa = \kappa(\mathbb{G}_1) > 0$ such that
\begin{align}
  |H_k(\bm{x}, G_k(\bm{x}) + \Delta y)|, \quad
  |H_{k, y}(\bm{x}, G_k(\bm{x}) + \Delta y)| \leq \frac{1}{2}
  \label{eq-A3-Hk_Hky_ev}
\end{align}
for any $\bm{x} \in U$, $\Delta y \in \mathbb{R}$, and $k > \kappa$.
When $\Delta y$ moves from $-1$ to 1,
the sign of \eqref{eq-A3-Fkexp} must change
from negative to positive.
Thus, by using the intermediate value theorem,
there exists unique $I_k(\bm{x}) \in \mathbb{R}$ such that
\begin{align}
  F_k(\bm{x}, G_k(\bm{x}) + I_k(\bm{x})) = 0
  \label{eq-A3-Ikprop}
\end{align}
for each $\bm{x} \in U$ and $k > \kappa$.
Since $F_{k, y} = 1 - H_{k, y} \neq 0$,
the implicit function theorem yields
$I_k: U \to \mathbb{R}$ are, in fact, $C^r$ functions.

By the mean value theorem,
there exists $\theta_k = \theta_k(\bm{x}) \in (0, 1)$ such that
\begin{align*}
  0 = F_k(\bm{x}, G_k(\bm{x}) + I_k(\bm{x}))
  &= F_k(\bm{x}, G_k(\bm{x}))
    + F_{k, y}(\bm{x}, G_k(\bm{x}) + \theta_k I_k(\bm{x}))
    I_k(\bm{x}) \\
  &= -H_k(\bm{x}, G_k(\bm{x})) + (
      1 - H_{k, y}(\bm{x}, G_k(\bm{x}) + \theta_k I_k(\bm{x}))
  ) I_k(\bm{x}).
\end{align*}
Thus, by using \eqref{eq-A3-Hk_Hky_ev},
\begin{align*}
  |I_k(\bm{x})|
  = \left|\frac{H_k(\bm{x}, G_k(\bm{x}))}
  {1 - H_{k, y}(\bm{x}, G_k(\bm{x}) + \theta_k I_k(\bm{x}))} \right|
  \leq 2 \mathcal{H}_k.
\end{align*}
We complete the proof.
\end{proof}

\subsection{Estimate of partial derivatives}
\label{sC2-ev}
In this section,
we assume $r \geq 3$ and
give estimates of partial derivatives of $I_k(\bm{x})$ up to order three,
where $I_k(\bm{x})$ is the function in Proposition~\ref{pc1-sgl}.

\medskip

For any finite $l \in \mathbb{Z}_{> 0}$ with $l \leq r$
and $\sigma_1$, $\sigma_2$, $\cdots$, $\sigma_l \in \{x_1, x_2, \cdots, x_n\}$,
we define
\begin{align*}
  &\mathcal{G}_k^{(\sigma_1 \sigma_2 \cdots \sigma_l)}
    = \mathcal{G}_k^{(\sigma_1 \sigma_2 \cdots \sigma_l)}(\mathbb{G}_1)
    := \sup_{\bm{x} \in U}
    |G_{k, \sigma_1 \sigma_2 \cdots \sigma_l}(\bm{x})|, \\
  &\mathcal{H}_k^{(\sigma_1 \sigma_2 \cdots \sigma_l)}
    = \mathcal{H}_k^{(\sigma_1 \sigma_2 \cdots \sigma_l)}(\mathbb{G}_1)
    := \max_{ \sigma_1', \sigma_2', \cdots, \sigma_l'}
      \sup_{\bm{x} \in U, \, y \in \mathbb{R}}
      |H_{k, \sigma_1' \sigma_2' \cdots \sigma_l'}(\bm{x}, y)|,
\end{align*}
where the variable $\sigma_i'$ is either equal to $\sigma_i$ or $y$
for each $i \in \{1, 2, \dots, l\}$.
We further define
\begin{align*}
  \hat{\mathcal{H}}_k^{(\sigma_1 \sigma_2 \cdots \sigma_l)}
    = \hat{\mathcal{H}}_k^{(\sigma_1 \sigma_2 \cdots \sigma_l)}(\mathbb{G}_1)
    := \max_{\tau \sqsubset (\sigma_1 \sigma_2 \cdots \sigma_l)}
      \mathcal{H}_k^{(\tau)},
\end{align*}
where $\tau \sqsubset (\sigma_1 \sigma_2 \cdots \sigma_l)$ means that
$\tau$ is a nonempty subsequence of the sequence $\sigma_1 \sigma_2 \cdots \sigma_l$
that preserves the original order.
That is, there exist indices $1 \le i_1 < i_2 < \cdots < i_{l'} \le l$
such that $\tau = \sigma_{i_1} \sigma_{i_2} \dots \sigma_{i_{l'}}$.

\begin{proposition}[Estimate for a single equation] \label{pc2-ev}
  For any $\sigma_1$, $\sigma_2$, $\sigma_3 \in \{x_1, x_2, \cdots, x_n\}$,
  we have the following three statements:
  \begin{enumerate}
    \item If $G_k(\bm{x})$ is constant, then
      \begin{align}
        |I_{k, \sigma_1}| \leq C \sup_{\bm{x} \in U, \, y \in \mathbb{R}}
          |H_{k, \sigma_1}(\bm{x}, y)|
        \label{eq-A3-Ik1}
      \end{align}
      for some constant $C = C(\mathbb{G}_1) > 0$.
      Otherwise, we have  
      \begin{align}
        |I_{k, \sigma_1}| 
          \leq C (1 + \mathcal{G}_k^{(\sigma_1)}) \mathcal{H}_k^{(\sigma_1)}
          \label{eq-A3-calFIk1}
      \end{align}
      for some constant $C = C(\mathbb{G}_1) > 0$.
    \item Assume $r \geq 2$ and $\hat{\mathcal{H}}_k^{(\sigma_1 \sigma_2)} \to 0$
      as $k \to \infty$.
      Then, the second partial derivatives of $I_k(\bm{x})$
      are estimated as
      \begin{align}
        &|I_{k, \sigma_1 \sigma_2}|
        \leq C (1 + \mathcal{G}_k^{(\sigma_1)}
          + \mathcal{G}_k^{(\sigma_2)}
          + \mathcal{G}_k^{(\sigma_1 \sigma_2)}
          + \mathcal{G}_k^{(\sigma_1)} \mathcal{G}_k^{(\sigma_2)})
            \hat{\mathcal{H}}_k^{(\sigma_1 \sigma_2)},
        \label{eq-A3-calFIk2}
      \end{align}
      for some constant $C = C(\mathbb{G}_1) > 0$.
    \item Assume $r \geq 3$ and
      $\hat{\mathcal{H}}_k^{(\sigma_1 \sigma_2 \sigma_3)} \to 0$
      as $k \to \infty$.
      Then, the third partial derivatives of $I_k(\bm{x})$
      are estimated as
      \begin{align}
        &\begin{aligned}
          |I_{k, \sigma_1 \sigma_2 \sigma_3}|
          &\leq C (1 + \mathcal{G}_k^{(\sigma_1)}
          + \mathcal{G}_k^{(\sigma_2)}
          + \mathcal{G}_k^{(\sigma_3)} 
          + \mathcal{G}_k^{(\sigma_1 \sigma_2)}
          + \mathcal{G}_k^{(\sigma_2 \sigma_3)} 
          + \mathcal{G}_k^{(\sigma_1 \sigma_3)}
          + \mathcal{G}_k^{(\sigma_1 \sigma_2 \sigma_3)} \\
          &\quad + \mathcal{G}_k^{(\sigma_1)} \mathcal{G}_k^{(\sigma_2)} 
          + \mathcal{G}_k^{(\sigma_2)} \mathcal{G}_k^{(\sigma_3)} 
          + \mathcal{G}_k^{(\sigma_1)} \mathcal{G}_k^{(\sigma_3)} \\
          &\quad + \mathcal{G}_k^{(\sigma_1 \sigma_2)} \mathcal{G}_k^{(\sigma_3)} 
          + \mathcal{G}_k^{(\sigma_1)} \mathcal{G}_k^{(\sigma_2 \sigma_3)} 
          + \mathcal{G}_k^{(\sigma_1 \sigma_3)} \mathcal{G}_k^{(\sigma_2)}
          + \mathcal{G}_k^{(\sigma_1)} \mathcal{G}_k^{(\sigma_2)}
            \mathcal{G}_k^{(\sigma_3)})
            \hat{\mathcal{H}}_k^{(\sigma_1 \sigma_2 \sigma_3)}
        \end{aligned} \label{eq-A3-calFIk3}
      \end{align}
      for some constant $C = C(\mathbb{G}_1) > 0$.
  \end{enumerate}
\end{proposition}

\begin{proof}[Proof of Proposition~\ref{pc2-ev}]
We prove the three assertions of the above Proposition in parallel.  
By differentiating \eqref{eq-A3-Ikprop} with respect to 
$\sigma_1 \in \{x_1, x_2, \cdots, x_n\}$, we obtain
\begin{align}
  I_{k, \sigma_1} = H_{k, \sigma_1} \tilde{H}
    + G_{k, \sigma_1} H_{k, y} \tilde{H},
  \label{eq-A3-Iks1}
\end{align}
where $\tilde{H} := (1 - H_{k, y})^{-1}$.
Thus, we get the desired formulas \eqref{eq-A3-Ik1} and \eqref{eq-A3-calFIk1}.
Note that $G_{k, \sigma} + I_{k, \sigma}
= (G_{k, \sigma} + H_{k, \sigma})\tilde{H}$
for any $\sigma \in \{x_1, x_2, \cdots, x_n\}$.
By the chain rule, we have  
\begin{align*}
  &\partial_{\sigma} H_{k, \sigma_1 \sigma_2 \cdots \sigma_l}  
  = H_{k, \sigma_1 \sigma_2 \cdots \sigma_l \sigma}  
    + (G_{k, \sigma} + H_{k, \sigma})
      H_{k, \sigma_1 \sigma_2 \cdots \sigma_l y} \tilde{H} \\
  &\partial_{\sigma} \tilde{H}  
  = H_{k, y \sigma} \tilde{H}^2
    + (G_{k, \sigma} + H_{k, \sigma}) H_{k, y y} \tilde{H}^3
\end{align*}
for any $\sigma \in \{x_1, x_2, \cdots, x_n\}$,
$\sigma_1 \sigma_2 \cdots \sigma_l \in \{x_1, x_2, \cdots, x_n\}^l$,
and $l < r$.
By differentiating \eqref{eq-A3-Iks1} with respect to 
$\sigma_2 \in \{x_1, x_2, \cdots, x_n\}$, the above formulas imply
\begin{align}
\begin{aligned}
  &I_{k, \sigma_1 \sigma_2}
  = H_{k, \sigma_1 \sigma_2} \tilde{H}
  + G_{k, \sigma_2} H_{k, \sigma_1 y} \tilde{H}^2
  + H_{k, \sigma_1 y} H_{k, \sigma_2} \tilde{H}^2
  + G_{k, \sigma_1 \sigma_2} H_{k, y} \tilde{H}
  + G_{k, \sigma_1} H_{k, y \sigma_2} \tilde{H} \\
  &\quad + H_{k, \sigma_1} H_{k, y \sigma_2} \tilde{H}^2
  + G_{k, \sigma_1} H_{k, y} H_{k, y \sigma_2} \tilde{H}^2
  + G_{k, \sigma_1} G_{k, \sigma_2} H_{k, yy} \tilde{H}^2
  + G_{k, \sigma_1} H_{k, \sigma_2} H_{k, yy} \tilde{H}^2 \\
  &\quad + G_{k, \sigma_2} H_{k, \sigma_1} H_{k, yy} \tilde{H}^3
  + H_{k, \sigma_1} H_{k, \sigma_2} H_{k, yy} \tilde{H}^3
  + G_{k, \sigma_1} G_{k, \sigma_2} H_{k, y} H_{k, yy} \tilde{H}^3
  + G_{k, \sigma_1} H_{k, \sigma_2} H_{k, y} H_{k, yy} \tilde{H}^3.
\end{aligned} \label{ec2-Ik12}
\end{align}
Furthermore, when we take partial derivatives of each term in \eqref{eq-A3-Iks1},  
the coefficients that appear with respect to the partial derivatives of $H_k$  
and the variables $\tilde{H}$ are summarized in Table~\ref{tc2-I12}.  
\begin{table}[h]
  \centering
  \caption{Coefficient terms appearing after differentiating terms
    in \eqref{eq-A3-Iks1}}
  \label{tc2-I12}
  \begin{tabular}{l|l}
    \midrule
    \textbf{Term in \eqref{eq-A3-Iks1}} & \textbf{Coefficients after differentiation} \\
    \midrule
    $H_{k, \sigma_1} \tilde{H}$ & $1$, $G_{k, \sigma_2}$ \\
    $G_{k, \sigma_1} H_{k, y} \tilde{H}$ &
      $G_{k, \sigma_1 \sigma_2}$,
      $G_{k, \sigma_1}$,
      $G_{k, \sigma_1} G_{k, \sigma_2}$ \\
    \midrule
  \end{tabular}
\end{table}
Since all the absolute values of the partial derivatives of the $H_k$
in \eqref{ec2-Ik12}
are bounded by $\hat{\mathcal{H}}_k^{(\sigma_1 \sigma_2)}$  
and $\hat{\mathcal{H}}_k^{(\sigma_1 \sigma_2)}$ is infinitesimal,  
we obtain the desired formula \eqref{eq-A3-calFIk2}.  
Analogously, differentiating the above relation \eqref{ec2-Ik12}  
with respect to $\sigma_3 \in \{x_1, x_2, \cdots, x_n\}$,  
we obtain coefficients as summarized in Table~\ref{tc2-I123}.  
\begin{table}[h]
  \centering
  \caption{Coefficient terms appearing after differentiating terms
    in \eqref{ec2-Ik12}}
  \label{tc2-I123}
  \begin{tabular}{l|l}
    \midrule
    \textbf{Term in \eqref{ec2-Ik12}} &
      \textbf{Coefficients after differentiation} \\
    \midrule
    $H_{k, \sigma_1 \sigma_2} \tilde{H}$ & $1$, $G_{k, \sigma_3}$ \\
    $G_{k, \sigma_2} H_{k, \sigma_1 y} \tilde{H}^2$ &
      $G_{k, \sigma_2 \sigma_3}$, $G_{k, \sigma_2}$,
      $G_{k, \sigma_2} G_{k, \sigma_3}$ \\
    $H_{k, \sigma_1 y} H_{k, \sigma_2} \tilde{H}^2$ & $1$, $G_{k, \sigma_3}$ \\
    $G_{k, \sigma_1 \sigma_2} H_{k, y} \tilde{H}$ &
      $G_{k, \sigma_1 \sigma_2 \sigma_3}$,
      $G_{k, \sigma_1 \sigma_2}$,
      $G_{k, \sigma_1 \sigma_2}$, $G_{k, \sigma_3}$ \\
    $G_{k, \sigma_1} H_{k, y \sigma_2} \tilde{H}$ &
      $G_{k, \sigma_1 \sigma_3}$,
      $G_{k, \sigma_1}$,
      $G_{k, \sigma_1} G_{k, \sigma_3}$ \\
    $H_{k, \sigma_1} H_{k, y \sigma_2} \tilde{H}^2$ & $1$, $G_{k, \sigma_3}$ \\
    $G_{k, \sigma_1} H_{k, y} H_{k, y \sigma_2} \tilde{H}^2$ &
      $G_{k, \sigma_1 \sigma_3}$, $G_{k, \sigma_1}$,
      $G_{k, \sigma_1} G_{k, \sigma_3}$ \\
    $G_{k, \sigma_1} G_{k, \sigma_2} H_{k, yy} \tilde{H}^2$ &
      $G_{k, \sigma_1 \sigma_3} G_{k, \sigma_2}$, $G_{k, \sigma_1}
      G_{k, \sigma_2 \sigma_3}$, $G_{k, \sigma_1} G_{k, \sigma_2}$,
      $G_{k, \sigma_1} G_{k, \sigma_2} G_{k, \sigma_3}$ \\
    $G_{k, \sigma_1} H_{k, \sigma_2} H_{k, yy} \tilde{H}^2$ &
      $G_{k, \sigma_1 \sigma_3}$, $G_{k, \sigma_1}$,
      $G_{k, \sigma_1} G_{k, \sigma_3}$ \\
    $G_{k, \sigma_2} H_{k, \sigma_1} H_{k, yy} \tilde{H}^3$ &
      $G_{k, \sigma_1 \sigma_3}$, $G_{k, \sigma_1}$,
      $G_{k, \sigma_1} G_{k, \sigma_3}$ \\
    $H_{k, \sigma_1} H_{k, \sigma_2} H_{k, yy} \tilde{H}^3$ &
      $1$, $G_{k, \sigma_3}$ \\
    $G_{k, \sigma_1} G_{k, \sigma_2} H_{k, y} H_{k, yy} \tilde{H}^3$ &
      $G_{k, \sigma_1 \sigma_3} G_{k, \sigma_2}$,
      $G_{k, \sigma_1} G_{k, \sigma_2 \sigma_3}$,
      $G_{k, \sigma_1} G_{k, \sigma_2}$,
      $G_{k, \sigma_1} G_{k, \sigma_2} G_{k, \sigma_3}$ \\
    $G_{k, \sigma_1} H_{k, \sigma_2} H_{k, y} H_{k, yy} \tilde{H}^3$ &
      $G_{k, \sigma_1 \sigma_3}$, $G_{k, \sigma_1}$,
      $G_{k, \sigma_1} G_{k, \sigma_3}$ \\
    \midrule
  \end{tabular}
\end{table}
This result yields the desired formula \eqref{eq-A3-calFIk3}
in a similar manner.
We complete the proof.
\end{proof}

\subsection{System of equations}
\label{sC3-sys}

In this appendix,  
as an application of the previous results, we introduce a method for solving  
a system of equations.  
We also provide an estimate of the partial derivatives of the solutions  
under certain conditions.  

\medskip

Let $\{ G_k^{(j)}: U \to \mathbb{R} \}_{k \in \mathbb{Z}_{> 0},
j \in \{1, 2, \cdots, m\}}$,
$m \in \mathbb{Z}_{> 0}$,
be $C^r$ functions from an open set $U \subset \mathbb{R}^n$ to $\mathbb{R}$.
Let $\{ H_k^{(j)}: U \times \mathbb{R}^m \to \mathbb{R}\}_{
  k \in \mathbb{Z}_{> 0}, j \in \{1, 2, \cdots, m\}}$
be $C^r$ functions.
For the above core objects, we set
\begin{align*}
  \mathbb{G}_2 := (\{ G_k^{(j)}: U \to \mathbb{R} \},
  \{ H_k^{(j)}: U \times \mathbb{R}^m \to \mathbb{R} \}).
\end{align*}
Let $U \times \mathbb{R}^m$ has the coordinates
$(\bm{x}, \bm{y})$, where $\bm{x} = (x_1, x_2, \cdots, x_n)$
and $\bm{y} = (y_1, y_2, \cdots, y_m)$.
Let $\Sigma_{\bm{x}} := \{x_1, x_2, \cdots, x_n\}$
and $\Sigma_{\bm{y}} := \{y_1, y_2, \cdots, y_m\}$.

\begin{proposition}[Solution method for a system of equations] \label{pc3-sys}
  We have the following two statements:
  \begin{enumerate}
    \item Assume
      \begin{gather*}
        \mathscr{H}_k^{(j)} = \mathscr{H}_k^{(j)}(\mathbb{G}_1)
          := \sup_{\bm{x} \in U, \, \bm{y} \in \mathbb{R}^m}
          |H_k^{(j)}(\bm{x}, \bm{y})|
        \to 0 \quad \text{as} \quad k \to \infty, \\
        \max_{\sigma \in \Sigma_{\bm{y}}}
          \sup_{\bm{x} \in U, \, \bm{y} \in \mathbb{R}^m}
          |H_{k, \sigma}^{(j)}(\bm{x}, \bm{y})|
        \to 0 \quad \text{as} \quad k \to \infty
      \end{gather*}
      for any $j \in \{1, 2, \cdots, m\}$.
      Then, there exists $\kappa = \kappa(\mathbb{G}_2) > 0$ such that
      the system of equations of $(\bm{x}, \bm{y}) \in U \times \mathbb{R}^m$
      \begin{align*}
        y_j = G_k^{(j)}(\bm{x}) + H_k^{(j)}(\bm{x}, \bm{y}), \quad
        j \in \{1, 2, \cdots, m\}
      \end{align*}
      has the solution
      \begin{align*}
        y_j = G_k^{(j)}(\bm{x}) + I_k^{(j)}(\bm{x}), \quad
        j \in \{1, 2, \cdots, m\}
      \end{align*}
      for any $k > \kappa$,
      where
      $I_k^{(j)}: D \to \mathbb{R}$
      are $C^r$ functions such that
      there exists a constant $C = C(\mathbb{G}_2) > 0$ satisfying
      \begin{align}
        |I_k^{(j)}(\bm{x})| \leq C \mathscr{H}_k^{(j)}
        \label{eq-A3-Ik_eval}
      \end{align}
      for any $\bm{x} \in D $, $j \in \{1, 2, \cdots, m\}$, and $k > \kappa$.
    \item We further assume
      $G_k^{(j)}(\bm{x})$ is constant and
      \begin{align*}
        \max_{\sigma \in \Sigma_{\bm{x}} \cup \Sigma_{\bm{y}}}
          \sup_{\bm{x} \in U, \, \bm{y} \in \mathbb{R}^m}
          |H_{k, \sigma}^{(j)}(\bm{x}, \bm{y})|
          \to 0 \quad \text{as} \quad k \to \infty
      \end{align*}
      for any $j \in \{1, 2, \cdots, m\}$.
      Then, there exists a constant $C = C(\mathbb{G}_2) > 0$ such that
      the first partial derivatives of the solution $I_k^{(j)}(\bm{x})$
      are estimated as
      \begin{align}
        |I_{k, \sigma}^{(j)}|
        \leq C \max_{\sigma' \in \{\sigma\} \cup \Sigma_{\bm{y}} \setminus \{y_j\}}
          \sup_{\bm{x} \in U, \, \bm{y} \in \mathbb{R}^m}
          |H_{k, \sigma'}^{(j)}(\bm{x}, \bm{y})|
        \label{eq-A3-Ikp_ev}
      \end{align}
      for any $j \in \{1, 2, \cdots, m\}$
      and $\sigma \in \Sigma_{\bm{x}}$.
  \end{enumerate}
\end{proposition}

\begin{proof}[Proof of Proposition~\ref{pc3-sys}]
We divide the proof into two parts, corresponding to the first and second items.  

\medskip

\textbf{(1) First item.}\quad
We prove the first item by mathematical induction by $m$.
The case $m = 1$ is proved from Proposition~\ref{pc1-sgl}.
We assume that the first item holds for $m$
and prove that the first item also holds for $m + 1$.

The equations
\begin{align}
  y_j = G_k^{(j)}(\bm{x}) + H_k^{(j)}(\bm{x}, \bm{y}), \quad
  j \in \{2, 3, \cdots, m + 1\}
  \label{eq-A3-msys}
\end{align}
can be solved by the assumption;
they have the solutions
\begin{align}
  y_j = G_k^{(j)}(\bm{x}) + \tilde{I}_k^{(j)}(\bm{x}, y_1), \quad
  j \in \{2, 3, \cdots, m + 1\}
  \label{eq-A3-solutions}
\end{align}
with the estimate as in \eqref{eq-A3-Ik_eval}.
For the remaining equation
\begin{align}
  y_1 = G_k^{(1)}(\bm{x}) + H_k^{(1)}(\bm{x}, \bm{y}),
  \label{eq-A3-rem}
\end{align}
we substitute \eqref{eq-A3-solutions} into the above equation
and get the equation of $(\bm{x}, y_1)$.
To apply Proposition~\ref{pc1-sgl} for the equation,
it suffices to check
\begin{align*}
  \left|\partial_{y_1} H_k^{(1)}(\bm{z})\right|
  \to 0 \quad \text{as} \quad k \to \infty,
\end{align*}
where $\bm{z} = (\bm{x}, y_1,
G_k^{(2)}(\bm{x}) + \tilde{I}_k^{(2)}(\bm{x}, y_1),
G_k^{(3)}(\bm{x}) + \tilde{I}_k^{(3)}(\bm{x}, y_1), \cdots,
G_k^{(m + 1)}(\bm{x}) + \tilde{I}_k^{(m + 1)}(\bm{x}, y_1))$.
By $\tilde{I}_{k, y_1}^{(j)} = H_{k, y_1}^{(j)}(1 - H_{k, y_j}^{(j)})^{-1}$,
we have
\begin{align*}
  \left|\partial_{y_1} H_k^{(1)}(\bm{z})\right|
  &= \left| H_{k, y_1}^{(1)}(\bm{z})
    + \sum_{j = 2}^{m + 1} \tilde{I}_{k, y_1}^{(j)}(\bm{x}, y_1)
      H_{k, y_j}^{(1)}(\bm{z}) \right| \\
  &\leq |H_{k, y_1}^{(1)}(\bm{z})|
    + \sum_{j = 2}^{m + 1}
      \frac{|H_{k, y_1}^{(j)}(\bm{z})| |H_{k, y_j}^{(1)}(\bm{z})|}
      {|1 - H_{k, y_j}^{(j)}(\bm{z})|}
  \to 0 \quad \text{as} \quad k \to \infty
\end{align*}
due to the assumptions.
Thus, we obtain the solution
\begin{align}
  y_1 = G_k^{(1)}(\bm{x}) + I_k^{(1)}(\bm{x})
  \label{eq-A3-y1_sol}
\end{align}
with the estimate in \eqref{eq-A3-Ik_eval}.
Putting
\begin{align}
  I_k^{(j)}(\bm{x})
    := \tilde{I}_k^{(j)}(\bm{x}, G_k^{(1)}(\bm{x}) + I_k^{(1)}(\bm{x}))
  \label{eq-A3-Ikrem}
\end{align}
for any $j \in \{2, 3, \cdots, m + 1\}$,
we complete the proof of the first item.

\medskip

\textbf{(2) Second item.}\quad
We prove this again by mathematical induction on $m$.  
The case $m = 1$ follows from Proposition~\ref{pc2-ev}.  
Assume that the second item holds for $m$.  
We will prove that the second item also holds for $m + 1$.  

Since the assumption holds,
the equations \eqref{eq-A3-msys} have the solutions \eqref{eq-A3-solutions} with
\begin{align}
  |\tilde{I}_{k, \sigma}^{(j)}(\bm{x}, y_1)|
  \leq C_1 \max_{\sigma' \in \{\sigma\} \cup \Sigma_{\bm{y}} \setminus \{y_1, y_j\}}
  \sup_{\bm{x} \in U, \, \bm{y} \in \mathbb{R}^m}
  |H_{k, \sigma'}^{(j)}(\bm{x}, \bm{y})|
  \label{eq-A3-IkpRem_ev}
\end{align}
for any $j \in \{2, 3, \cdots, m + 1\}$
and $\sigma \in \Sigma_{\bm{x}} \cup \{ y_1 \}$,
where $C_1 = C_1(\mathbb{G}_2) > 0$ is some constant.
We substitute \eqref{eq-A3-solutions} into \eqref{eq-A3-rem}
and get the equation of $(\bm{x}, y_1)$:
\begin{align*}
  y_1 = G_k^{(1)}(\bm{x}) + H_k^{(1)}(\bm{z}), \quad
  \bm{z} = (\bm{x}, y_1,
    \tilde{I}_k^{(2)}(\bm{x}, y_1),
    \tilde{I}_k^{(3)}(\bm{x}, y_1), \cdots,
    \tilde{I}_k^{(m + 1)}(\bm{x}, y_1)).
\end{align*}
Now, we pick $\sigma \in \Sigma_{\bm{x}}$.
Applying Proposition~\ref{pc2-ev} for the equation,
we get the solution \eqref{eq-A3-y1_sol} with the estimate
\begin{align*}
  |I_{k, \sigma}^{(1)}|
  \leq C_2 \sup_{\bm{x} \in U, \, y_1 \in \mathbb{R}}
    |\partial_\sigma H_k^{(1)}(\bm{z})|
\end{align*}
for some constant $C_2 = C_2(\mathbb{G}_2) > 0$.
By the chain rule, we have
\begin{align*}
  \left|\partial_\sigma H_k^{(1)}(\bm{z})\right|
  = \left| H_{k, \sigma}^{(1)}(\bm{z})
    + \sum_{j = 2}^{m + 1} \tilde{I}_{k, \sigma}^{(j)}(\bm{x}, y_1)
      H_{k, y_j}^{(1)}(\bm{z}) \right|.
\end{align*}
Hence, by \eqref{eq-A3-IkpRem_ev}
and the assumption of the second item, we obtain
\begin{align}
  |I_{k, \sigma}^{(1)}|
  \leq C_3 \max_{\sigma' \in \{\sigma\} \cup \Sigma_{\bm{y}} \setminus \{y_1\}}
    \sup_{\bm{x} \in U, \, \bm{y} \in \mathbb{R}^m}
    |H_{k, \sigma'}^{(1)}(\bm{x}, \bm{y})|
  \label{ec3-Ik1}
\end{align}
for some constant $C_3 = C_3(\mathbb{G}_2) > 0$.
On the other hand,
differentiating both sides of \eqref{eq-A3-Ikrem} with respect to $\sigma$,  
and using \eqref{eq-A3-IkpRem_ev}, we obtain  
\begin{align}
  |I_{k, \sigma}^{(j)}|
  = |\tilde{I}_{k, \sigma}^{(j)}
    + \tilde{I}_{k, y_1}^{(j)} I_{k, \sigma}^{(1)}|
  \leq C_4 \max_{\sigma' \in \{\sigma\} \cup \Sigma_{\bm{y}} \setminus \{y_j\}}
  \sup_{\bm{x} \in U, \, \bm{y} \in \mathbb{R}^m}
  |H_{k, \sigma'}^{(j)}(\bm{x}, \bm{y})|
  \label{ec3-Ikj}
\end{align}
for any $j \in \{2, 3, \cdots, m + 1\}$,
where $C_4 = C_4(\mathbb{G}_2) > 0$ is some constant.
The results \eqref{ec3-Ik1} and \eqref{ec3-Ikj}
complete the proof of the second item.
\end{proof}

\end{appendices}

\ifdefined\isMaster
\else
  \end{document}
\fi

\section*{Acknowledgments}
I would like to express my sincere gratitude to Shuhei Hayashi
for his invaluable guidance throughout the preparation of this paper.
I also wish to thank Shin Kiriki, Yushi Nakano, and Teruhiko Soma;
without their involvement,
this research would never have begun.
I am also grateful to Sogo Murakami for his insightful comments and continuous support.
I thank Katsutoshi Shinohara for arranging opportunities to connect
with researchers in related fields. I am deeply indebted to
Dmitry Turaev, Dongchen Li, Xiaolong Li, and Dimitrii Mints
for their professional feedback on the content of this work.

\bibliography{references}
\bibliographystyle{plain}

\end{document}